\documentclass[12pt]{amsart}
\font\emailfont=cmtt10

\usepackage{amsmath,amsthm,amsfonts,amscd,flafter,epsf}

\newcommand\commentable[1]{#1}

\newcommand\Id{\mathrm{Id}}

\newcommand{\Tors}{\mathrm{Tors}}

\newcommand{\HF}{HF}

\newtheorem{theorem}{Theorem}[section]
\newtheorem{prop}[theorem]{Proposition}

\newtheorem{conj}[theorem]{Conjecture}
\newtheorem{lemma}[theorem]{Lemma}

\newtheorem{example}[theorem]{Example}
\newtheorem{defn}[theorem]{Definition}

\newtheorem{remark}[theorem]{Remark}

\def\endproof{\relax\ifmmode\expandafter\endproofmath\else
  \unskip\nobreak\hfil\penalty50\hskip.75em\hbox{}\nobreak\hfil\bull
  {\parfillskip=0pt \finalhyphendemerits=0 \bigbreak}\fi}
\def\endproofmath$${\eqno\bull$$\bigbreak}
\def\bull{\vbox{\hrule\hbox{\vrule\kern3pt\vbox{\kern6pt}\kern3pt\vrule}\hrule}}

\newcounter{bean}
\newcommand{\Q}{\mathbb{Q}}

\newcommand{\Z}{\mathbb{Z}}

\newcommand{\OneHalf}{\frac{1}{2}}

\newcommand{\Zmod}[1]{\Z/{#1}\Z}

\newcommand{\Ker}{\mathrm{Ker}}

\newcommand{\Coker}{\mathrm{Coker}}

\newcommand{\cm}{\cdot}

\newcommand{\Nbd}[1]{{\mathrm{nd}}(#1)}
\newcommand{\nbd}[1]{\Nbd{#1}}

\newcommand{\ModSWfour}{\mathcal{M}}
\newcommand{\ModFlow}{\ModSWfour}

\newcommand{\SpinC}{{\mathrm{Spin}}^c}

\newcommand{\goesto}{\mapsto}

\newcommand\Wedge{\Lambda}

\newcommand\Hom{\mathrm{Hom}}

\newcommand\abuts\Rightarrow
\newcommand\Sym{\mathrm{Sym}}

\newcommand\mCP{{\overline{\mathbb{CP}}}^2}

\newcommand\HFpRed{\HFp_{\red}}
\newcommand\HFmRed{\HFm_{\red}}
\newcommand\mix{\mathrm{mix}}

\newcommand\liftGr{\widetilde\gr}

\newcommand\SpanA{{\mathrm{Span}}([\alpha_i]_{i=1}^g)}
\newcommand\SpanB{{\mathrm{Span}}([\beta_i]_{i=1}^g)}
\newcommand\SpanC{{\mathrm{Span}}([\gamma_i]_{i=1}^g)}

\newcommand\x{\mathbf x}
\newcommand\w{\mathbf w}

\newcommand\y{\mathbf y}

\newcommand\ModSphere{\ModFlow\left({\mathbb S}\longrightarrow 
\Sym^{g-1}(\Sigma_{1})\times \Sym^2(\Sigma_{2})\right)}
\newcommand\ModSpheres\ModSphere
\newcommand\CF{CF}

\newcommand\CFa{\widehat{CF}}
\newcommand\CFp{\CFb}
\newcommand\CFm{\CF^-}
\newcommand\CFleq{\CF^{\leq 0}}
\newcommand\HFleq{\HF^{\leq 0}}

\newcommand\HFpred{\HFp_{\rm red}}

\newcommand\HFmred{\HFm_{\rm red}}
\newcommand\HFred{\HF_{\rm red}}
\newcommand\coHFm{\HF_-}

\newcommand{\red}{\mathrm{red}}

\newcommand\HFp{\HFb}

\newcommand\HFm{\HF^-}
\newcommand\CFinf{CF^\infty}
\newcommand\HFinf{HF^\infty}
\newcommand\CFb{CF^+}
\newcommand\HFa{\widehat{HF}}
\newcommand\HFb{HF^+}
\newcommand\gr{\mathrm{gr}}
\newcommand\Mas{\mu}
\newcommand\UnparModSp{\widehat \ModSp}
\newcommand\UnparModFlow\UnparModSp
\newcommand\Mod\ModSp

\newcommand{\cald}{{\mathcal D}}

\newcommand{\spinc}{\mathfrak s}

\newcommand{\spinct}{\mathfrak t}

\newcommand\ModMaps{\mathcal M}
\newcommand\ModSp\ModMaps

\newcommand\Ta{{\mathbb T}_{\alpha}}
\newcommand\Tb{{\mathbb T}_{\beta}}
\newcommand\Tc{{\mathbb T}_{\gamma}}
\newcommand\Td{{\mathbb T}_{\delta}}

\newcommand\alphas{\mbox{\boldmath$\alpha$}}

\newcommand\betas{\mbox{\boldmath$\beta$}}
\newcommand\gammas{\mbox{\boldmath$\gamma$}}
\newcommand\deltas{\mbox{\boldmath$\delta$}}

\newcommand\PerDom{\mathcal P}

\newcommand\uCF{\underline\CF}
\newcommand\uCFinf{\uCF^\infty}
\newcommand\uHF{\underline{\HF}}
\newcommand\uHFp{\underline{\HF}^+}

\newcommand\uHFm{\underline{\HF}^-}
\newcommand\uHFa{\underline{\HFa}}

\newcommand\uHFinf{\uHF^\infty}

\newcommand\uFinf[1]{{\underline F}^\infty_{#1}}

\newcommand\Fm[1]{F^{-}_{#1}}
\newcommand\Fp[1]{F^{+}_{#1}}
\newcommand\Fa[1]{{\widehat F}_{#1}}
\newcommand\Finf[1]{F^{\infty}_{#1}}
\newcommand\Fleq[1]{F^{\leq 0}_{#1}}

\newcommand\sj{\mathfrak j}
\newcommand\Conjugate{\mathfrak J}

\newcommand\Fmix[1]{F^{\mix}_{#1}}
\newcommand\ugc[1]{{\underline g}^\circ_{#1}}
\newcommand\uGc[1]{{\underline G}^\circ_{#1}}
\newcommand\uginf[1]{{\underline g}^\infty_{#1}}

\newcommand\uEc[1]{{\underline E}^\circ_{#1}}
\newcommand\Ec[1]{{E}^\circ_{#1}}
\newcommand\Gc[1]{{G}^\circ_{#1}}
\newcommand\ueinf[1]{{\underline e}^\infty_{#1}}

\newcommand\Fc[1]{F^\circ_{#1}}

\newcommand\uFc{{\underline F}^\circ}
\newcommand\ufc{{\underline f}^\circ}

\newcommand\uHFc{\uHF^\circ}
\newcommand\uCFc{\uCF^\circ}
\newcommand\HFc{\HF^\circ}
\newcommand\CFc{\CF^\circ}
\newcommand\CFpm{\CF^\pm}

\newcommand\Dual{\mathcal D}
\newcommand\Duality\Dual

\newcommand\Tor{\mathrm{Tor}}
\newcommand\TaPr{\Ta'}
\newcommand\TbPr{\Tb'}
\newcommand\TcPr{\Tc'}

\newcommand\alphaprs{\alphas'}
\newcommand\betaprs{\betas'}
\newcommand\gammaprs{\gammas'}
\newcommand\Spider{\sigma}
\newcommand\EulerMeasure{\widehat\chi}

\newcommand\Knot{\mathbb K}
\newcommand\Link{\mathbb L}

\title[{Holomorphic triangles and four-manifold invariants}]
{Holomorphic triangles and invariants for smooth four-manifolds}

\author[Peter Ozsv{\'a}th]{Peter Ozsv\'ath}
\address{Department of
Mathematics,  Princeton University, New Jersey 08540 \newline
\indent{\emailfont{petero@math.princeton.edu}}}
\thanks{PSO was supported by NSF grant number DMS 9971950 and a Sloan 
Research Fellowship}

\author[Zolt{\'a}n Szab{\'o}]{Zolt{\'a}n Szab{\'o}} 
\address{Department of
Mathematics,  Princeton University, New Jersey 08540 \newline
\indent{\emailfont{szabo@math.princeton.edu}}}
\thanks{ZSz was supported by NSF grant number DMS 0107792
and a Packard Fellowship}

\begin{document}

\begin{abstract}  
	The aim of this article is to introduce invariants of
	oriented, smooth, closed four-manifolds, built using the Floer
	homology theories defined in~\cite{HolDisk} and
	\cite{HolDiskTwo}. This four-dimensional theory also endows
	the corresponding three-dimensional theories with additional
	structure: an absolute grading of certain of its Floer
	homology groups.  The cornerstone of these constructions is
	the study of holomorphic disks in the symmetric products of
	Riemann surfaces.  \end{abstract}

\maketitle
\section{Introduction}

In the first part of the paper, which comprises
Sections~\ref{sec:Cobordisms}-\ref{sec:BlowUp}, we set up the
invariant of a cobordism between two three-manifolds, which gives a
map between Floer homologies. In a related construction given in
Section~\ref{sec:AbsGrade}, we define an absolute $\Q$ lift of the
relative $\Z$-graded homology groups associated to a three-manifold
equipped with a torsion $\SpinC$ structure (by which we mean a
$\SpinC$ structure whose first Chern class is a torsion cohomology
class). In Section~\ref{sec:DefMixed}, we define a refined {\em mixed}
cobordism invariant which readily gives a smooth four-manifold
invariant (defined in Section~\ref{sec:ClosedInvariant}) for closed
four-manifolds $X$ with $b_2^+(X)>1$. We turn now to a more detailed
overview.


\subsection{Invariants of cobordisms}
In~\cite{HolDisk} and~\cite{HolDiskTwo}, we defined Floer homology
theories for oriented three-manifolds equipped with $\SpinC$ structures,
$\HFm(Y,\spinct)$, $\HFinf(Y,\spinct)$, and
$\HFp(Y,\spinct)$. (There is a fourth invariant, $\HFa(Y,\spinct)$, which
we will not discuss in this introduction.)  These invariants fit into a
long exact sequence
$$\begin{CD}
... @>>>\HFm(Y,\spinct)@>{\iota}>>\HFinf(Y,\spinct) @>{\pi}>> \HFp(Y,\spinct) @>{\delta}>>...
\end{CD}
$$
We abbreviate this long exact sequence $\HFc(Y,\spinct)$. (A rapid
overview of the general properties of these constructions, together
with a proof of their topological invariance in a stronger form, is
given in Section~\ref{sec:Naturality}.)
Recall that there is another associated three-manifold invariant, denoted
$$\HFpRed(Y,\spinc)=\Coker(\pi)\cong \Ker(\iota)=\HFmRed(Y,\spinc).$$
The isomorphism between $\HFpRed(Y,\spinc)\cong \HFmRed(Y,\spinc)$ is induced by
the coboundary map.

In the first part of the paper, we construct the invariant of a
connected cobordism $W$ between two connected three-manifolds $Y_1$
and $Y_2$, defined using the holomorphic triangle construction and a
handle-decomposition of $W$. Specifically, these constructions give
rise to a chain map between the chain complexes from $Y_1$ to $Y_2$,
whose induced maps on homology are invariants of $W$ (i.e. they are
independent of the handle decomposition).

\begin{theorem}
\label{intro:CobordismInvariant}
The maps on homology induced by a smooth, oriented
cobordism $W$ equipped with a
$\SpinC$ structure $\spinc\in\SpinC(W)$ are 
invariants of the cobordism, inducing a map of long exact
sequences:
$$
\begin{CD}
...
@>>>\HFm(Y_1,\spinct_1)@>{\iota_1}>>\HFinf(Y_1,\spinct_1) @>{\pi_1}>>\HFp(Y_1,\spinct_1)
@>>>...	\\
&& @V{F^-_{W,\spinc}}VV
@V{F^\infty_{W,\spinc}}VV
@V{F^+_{W,\spinc}}VV \\
... @>>>\HFm(Y_2,\spinct_2)@>{\iota_2}>>
\HFinf(Y_2,\spinct_2) @>{\pi_2}>>
\HFp(Y_2,\spinct_2)
@>>>...,
\end{CD}
$$ 
(where $\spinct_i\in\SpinC(Y_i)$ denotes the restriction of $\spinc$ to $Y_i$),
where the vertical maps are uniquely determined up to an overall sign, and
all squares are commutative.
\end{theorem}

The map on long exact sequences induced by the cobordism $W$ and
$\SpinC$ structure $\spinc$ is abbreviated $\Fc{W,\spinc}$.
Various refinements of the above map, including one using the action of first the homology of $W$,
and one using twisted coefficients, can be found in Section~\ref{sec:Cobordisms}.

The maps satisfy certain general properties: duality, conjugation
invariance, a blow-up formula, and composition properties. The duality
property relates the map induced by $F_{W,\spinc}$, thought of as a
cobordism from $Y_1$ to $Y_2$, with the induced map obtained by
thinking of $W$ as a cobordism from $-Y_2$ to $-Y_1$. Conjugation
invariance sets up an identification between the map induced by $W$
and the $\SpinC$ structure $\spinc$ with the same cobordism, equipped
with its conjugate $\SpinC$ structure ${\overline\spinc}$. The blow-up
formula relates the map induced by $W$, with that induced by the
(internal) connected sum of $W$ with $\mCP$ (the complex projective
plane, given the opposite of its complex orientation). 

The composition law states that if $W_1$ is a cobordism from $Y_1$ to
$Y_2$ and $W_2$ is a cobordism from $Y_2$ to $Y_3$, and we equip
$W_1$ and $W_2$ with $\SpinC$ structures $\spinc_1$ and $\spinc_2$ respectively (whose restrictions
agree over $Y_2$), then we have the following
relationship between the composition of $F_{W_1,\spinc_1}$ with $F_{W_2,\spinc_2}$, and the maps induced
by the composite cobordism $W=W_1\#_{Y_2} W_2$:
$$\Fc{W_2,\spinc_2}\circ
\Fc{W_1,\spinc_1} =
\sum_{\{\spinc\in\SpinC(W)\big|\spinc|W_1=\spinc_1,
\spinc|W_2=\spinc_2\}} 
\pm \Fc{W,\spinc}.$$

All of these properties (and various refinements) are stated precisely
in Section~\ref{sec:Cobordisms}.  The invariants are constructed in
Section~\ref{sec:ECobord}, and the properties are verified in
Sections~\ref{sec:ECobord}-\ref{sec:BlowUp}.

\subsection{Absolute gradings}

The principles used in construction of the cobordism invariant also
allow us to define an absolute $\Q$-lift of the relative $\Z$-gradings
of the Floer homology groups $\HFc(Y,\spinct)$ for a three-manifold
$Y$ equipped with a torsion $\SpinC$ structure $\spinct$. This lift is constructed in
Section~\ref{sec:AbsGrade}.

The relationship between the cobordism invariant and this absolute
grading is codified in the following formula, which holds for any cobordism $W$ from $Y_1$ to $Y_2$,
equipped with a $\SpinC$ structure $\spinc$ whose restriction to the two boundary components
$\spinct_1$ and $\spinct_2$ are both torsion:
$$\liftGr(\Fp{W,\spinc}(\xi))-\liftGr(\xi) = 
\frac{c_1(\spinc)^2-2\chi(W)-3\sigma(W)}{4}
$$
(where $\xi\in\HFp(\xi,\spinct_1)$ is any homogeneous cohomology class, and $\liftGr$ denotes the
absolute $\Q$ degree). 

\subsection{Closed four-manifold invariants}

There is a variant of the cobordism invariant which can be defined for
cobordisms $W$ with $b_2^+(W)>1$.

Observe first that if $W$ is a cobordism with $b_2^+(W)>0$, then the
map $\Finf{W,\spinc}$ induced on $\HFinf$ is trivial
(c.f. Lemma~\ref{lemma:BTwoPlusLemma}). If we have a cobordism with
$b_2^+(W)>1$, then we can cut $W$ along a three-manifold $N$, to
divide it into two cobordisms $W_1$ and $W_2$, both of which have
$b_2^+(W_i)>0$, in such a way that the map induced by restriction
$$\SpinC(W)\longrightarrow \SpinC(W_1)\times \SpinC(W_2)$$ is
injective (i.e. $\delta H^1(Y;\Z)\subset H^2(W)$ is trivial). 
In view of these remarks, if $\spinc$ is a $\SpinC$
structure and $\spinc_1$ and $\spinc_2$ denote its restriction to
$W_1$ and $W_2$ respectively, then $$\Fm{W_1,\spinc_1}\colon
\HFm(Y_1,\spinct_1) \longrightarrow
\HFm(N,\spinct)$$ factors through the inclusion of
$\HFmRed(N,\spinct)\hookrightarrow \HFm(N,\spinct)$, while
$$\Fp{W_2,\spinc_2}\colon \HFp(N,\spinct)\longrightarrow
\HFp(Y_2,\spinct_2)$$ factors through the projection
$\HFp(N,\spinct)\longrightarrow \HFpRed(N,\spinct)$.  Thus, by using
the identification of $\HFpRed(N,\spinct)\cong \HFmRed(N,\spinct)$ in
the middle, we can define the ``mixed invariant'' as a map
$$\Fmix{W,\spinc}\colon \HFm(Y_1,\spinct_1)\longrightarrow \HFp(Y_2,\spinct_2).$$

When $X$ is a closed four-manifold with $b_2^+(X)>0$, we can puncture
it in two points, and view the resulting object as a cobordism from
$S^3$ to $S^3$. By elaborating on the mixed invariant construction in
Section~\ref{sec:ClosedInvariant}, we obtain a map $$
\Phi_{\spinc} \colon \Z[U]\otimes \Wedge^* (H_1(X;\Z)/\Tors) \longrightarrow \Z,
$$ which is a homogeneous polynomial of degree $$
d(\spinc)=\frac{c_1(\spinc)^2-(2\chi(X)+3\sigma(X))}{4}, $$ and
well-defined up to sign.

It is, of course, interesting to compare this invariant with the
Seiberg-Witten invariant~\cite{Witten}. Indeed, we formulate the following:

\begin{conj}
Let $X$ be a closed, oriented, smooth four-manifold with $b_2^+(X)>1$.
Then, the invariant $\Phi_{X,\spinc}$ agrees with the Seiberg-Witten
invariant for $X$ in the $\SpinC$ structure $\spinc$.
\end{conj}

The invariant $\Phi$ can be used to give more directly topological
``gauge-theory free'' proofs of some facts about smooth four-manifolds
which have been previously established by means of Donaldson polynomials and
Seiberg-Witten invariants. We return to calculations of $\Phi$ in
\cite{HolDiskSymp}, and its non-vanishing properties
for symplectic manifolds, after 
calculating $\HFp$ for a number of three-manifolds
in~\cite{HolDiskGrad}.  These calculations are based on the surgery
long exact sequences from~\cite{HolDiskTwo}, combined with the
absolute $\Q$ grading defined in the present article. We content
ourselves here with some general properties and vanishing results, all
of which have natural analogues in Seiberg-Witten theory.

\subsection{Basic properties of the closed invariant}

The following is an analogue of Donaldson's connected sum theorem for
his polynomial invariants~\cite{DonaldsonPolynomials}. Unlike its gauge-theoretic
counterpart, the result follows rather directly from the definition of
the invariants.

\begin{theorem}
\label{intro:VanishingTheorem}
Let $X_1$ and $X_2$ be a pair of smooth, oriented four-manifolds with
$b_2^+(X_1)$, $b_2^+(X_2)>0$. Then, the invariants $\Phi_{X,\spinc}$
for the connected sum $X=X_1\# X_2$ vanish identically, for all
$\SpinC$ structures.
\end{theorem}

The blow-up formula for the cobordism invariant
translates directly into a corresponding 
blow-up formula for $\Phi_{X,\spinc}$
(compare~\cite{FintSternBlowup}):

\begin{theorem}
\label{intro:BlowUp}
Let $X$ be a closed, smooth, four-manifold with $b_2^+(X)>1$,
and let ${\widehat X}=X\#\mCP$ be its blowup.
Then, for each $\SpinC$ structure ${\widehat \spinc}\in\SpinC({\widehat X})$, with
$d({\widehat X},{\widehat \spinc})\geq 0$ we have the relation
$$\Phi_{{\widehat X},{\widehat \spinc}}(U^{\frac{\ell(\ell+1)}{2}}\cm \xi)=\Phi_{X,\spinc}(\xi),$$
where $\spinc$ is the $\SpinC$ structure over $X$ which agrees over $X-B^4$ with the restriction of
${\widehat\spinc}$, $\xi\in\Z[U]\otimes\Wedge^*(H_1(X)/\Tors)$ is any element of degree 
$d(X,\spinc)$, and $\ell$ is determined by 
$\langle c_1({\widehat \spinc}),[E]\rangle = \pm (2\ell+1)$,
where $E\subset{\widehat X}$ is the exceptional sphere.
\end{theorem}

For three-manifolds, Theorem~\ref{HolDiskTwo:thm:Adjunction}
of~\cite{HolDiskTwo} gives bounds on the Thurston norm in terms the
first Chern classes of $\SpinC$ structures for which $\HFp$ is
non-trivial. These ``adjunction inequalities'' have a straightforward
generalization in the four-dimensional context (compare~\cite{KMthom}, 
\cite{MSzT}, \cite{SympThom}).

\begin{theorem}
\label{intro:Adjunction}
Let $\Sigma\subset X$ be a homologically non-trivial embedded surface
with genus $g\geq 1$ and with non-negative self-intersection
number. Then, for each $\SpinC$ structure $\spinc\in\SpinC(X)$ for
which $\Phi_{X,\spinc}\not\equiv 0$, we have that $$\langle
c_1(\spinc),[\Sigma]\rangle + [\Sigma]\cdot[\Sigma] \leq 2g-2.$$
\end{theorem}

\vskip.2cm
\noindent{\bf Acknowledgements:} It is our pleasure to thank 
Tomasz Mrowka,, Robion Kirby, and Andr\'as Stipsicz for helpful
discussions during the preparation of this work.

\section{Preliminaries on Floer homology groups and
holomorphic triangles}
\label{sec:Naturality}

The aim of the present section is to collect some of the basic
properties of the three-manifold invariants introduced
in~\cite{HolDisk} and \cite{HolDiskTwo}, and to establish their
topological invariance in a form which will be useful to us in the
subsequent sections. After stating this result, a naturality result
for the isomorphism induced by a diffeomorphism, we turn to a general
discussion of the maps between the homology groups, as given using
Heegaard triples: recalling the four-dimensional interpretation
of Heegaard triples in Subsection~\ref{subsec:TopPrelim} 
(compare Section~\ref{HolDiskTwo:subsec:TriangleTopPrelim}
of~\cite{HolDiskTwo} for a more detailed discussion), and then in
Subsection~\ref{subsec:HolTriangles} stating the basic properties of
the induced maps. After recalling the calculations of the Floer
homology of $\#^n(S^1\times S^2)$, which is basic to much of the
present theory, we return to a proof of the naturality result in
Subsection~\ref{subsec:ProofNaturality}.

\subsection{Floer homology groups}

To fix terminology, a {\em set of attaching circles}
$\alphas=\{\alpha_1,...,\alpha_g\}$ in a genus $g$ surface
$\Sigma$ is a $g$-tuple of
pairwise disjoint, simple, closed curves in $\Sigma$ whose homology
classes are linearly independent in $H_1(\Sigma;\Z)$. A set of
attaching circles gives rise to a $g$-dimensional torus
$\Ta=\alpha_1\times...\times\alpha_g$ inside the $g$-fold symmetric
product of $\Sigma$, $\Sym^g(\Sigma)$.  A {\em pointed Heegaard
diagram} is a collection $(\Sigma,\alphas,\betas,z)$ where $\Sigma$ is an
oriented two-manifold,  $\alphas$ and $\betas$ are 
sets of attaching circles in $\Sigma$, and $z$ is a fixed reference
point in $\Sigma$ which is disjoint from all of the attaching
circles. A Heegaard diagram $(\Sigma,\alphas,\betas)$
specifies an oriented three-manifold $Y$ in a natural way.

In fact, given an oriented three-manifold equipped with a metric and
a self-indexing Morse function with unique index zero and three
critical points, there is an induced Heegaard diagram for $Y$ for
which the mid-level set is the Heegaard surface, points lying on
$\alpha_i$ (for any $i=1,...,g$) are the points in the mid-level which
flow out of the index one critical points under upward gradient flow,
and points lying on the $\beta_j$ are points in the mid-level which
flow into the index two critical points.  We say that two Heegaard
diagrams are {\em equivalent} if they are associated to two different
such Morse functions (and metrics) on the same three-manifold.  If two
Heegaard diagrams are equivalent, they can be connected by a sequence
of isotopies, handleslides, and stabilizations 
\footnote{For us, stabilization means
introducing a pair of canceling $\alpha$- and $\beta$-curves which are supported in
in a torus connected summand of the Heegaard surface. It is easy to see that more general stabilizations -- where the canceling curves are allowed to wander over $\Sigma$ -- differ from such stabilizations by a sequence of handleslides.}
As described
in~\cite{HolDisk}, the basepoint $z$ in the pointed Heegaard diagram
induces a map from intersection points of the tori $\Ta$ and $\Tb$ in
$\Sym^g(\Sigma)$ to $\SpinC$ structures over $Y$.

Suppose $Y$ is a three-manifold, equipped with a fixed $\SpinC$
structure $\spinct$. To a pointed Heegaard diagram
$(\Sigma,\alphas,\betas,z)$ for $Y$ (which satisfy certain additional
admissibility hypotheses when $b_1(Y)>0$, which we will recall
shortly) we can associate four chain complexes. These complexes are
$\CFa(\alphas,\betas,\spinct)$, which is freely generated by
intersection points between $\Ta$ and $\Tb$ representing the $\SpinC$
structure $\spinct$, $\CFinf(\alphas,\betas,\spinct)$, which is freely
generated by pairs consisting of such intersection and an integer, a
subcomplex $\CFm(\alphas,\betas,\spinct)$ where the integer is
required to be negative, and a quotient complex
$\CFp(\alphas,\betas,\spinct)$. The boundary maps are defined by
counting pseudo-holomorphic Whitney disks in $\Sym^g(\Sigma)$, and so
the definition of the boundary maps use some auxiliary data,
including a complex structure over $\Sigma$, and a one-parameter
variation in the induced almost-complex structure over
$\Sym^g(\Sigma)$, which we suppress from the notation for the chain complex
(indeed, as we
have suppressed the surface $\Sigma$ and its basepoint $z$). Each of
$\CFinf(\alphas,\betas,\spinct)$, $\CFm(\alphas,\betas,\spinct)$ and
$\CFp(\alphas,\betas,\spinct)$ comes equipped with the action of a
chain map $U$ which is induced from the map $U[\x,i]=[\x,i-1]$ on
$\CFinf(\alphas,\betas,\spinct)$ which decreases the (relative)
grading by two.

Recall that when $b_1(Y)>0$, we work with Heegaard
diagrams satisfying additional ``admissibility'' hypotheses.
To state these hypotheses, note that the two-dimensional homology
classes in $Y$ can be thought of as two-chains $\PerDom$ in $\Sigma$
whose local multiplicity at the reference point $z$
is zero, and 
whose boundary can be represented as a linear combination of
cycles representing elements of the $\alphas$ and
the $\betas$. Such two-chains are called {\em periodic domains}. A Heegaard
diagram is {\em strongly $\spinct$-admissible} for the $\SpinC$
structure $\spinct$ if for each (non-trivial) periodic domain $\PerDom$
for which $\langle c_1(\spinct), H(\PerDom)\rangle = 2n \geq 0$, some
local multiplicity in $\PerDom$ is greater than $n$. When we consider
$\CFm(\alphas,\betas,\spinct)$ and $\CFinf(\alphas,\betas,\spinct)$, we
will always be using strongly $\spinct$-admissible Heegaard
diagrams. When working with $\CFa(\alphas,\betas,\spinct)$ and
$\CFp(\alphas,\betas,\spinct)$, it suffices to work with {\em weakly
$\spinct$-admissible Heegaard diagrams} which are those for which each
non-trivial periodic domain $\PerDom$ with $\langle
c_1(\spinct),H(\PerDom)\rangle =0$ has both positive and negative
coefficients. Every Heegaard diagram for $Y$ is isotopic
to a strongly $\spinct$-admissible Heegaard diagram 
(see Section~\ref{HolDiskTwo:sec:Special} of~\cite{HolDiskTwo}).

It is shown in~\cite{HolDisk} and \cite{HolDiskTwo} that homology
groups of these chain complexes -- $\HFa(Y,\spinct)$,
$\HFinf(Y,\spinct)$, $\HFm(Y,\spinct)$, and $\HFp(Y,\spinct)$ -- are
topological invariants; where the latter three graded groups are
thought of as $\Z[U]$ modules. (In particular, they are independent of
the complex structure over $\Sigma$, the path of almost-complex
structures over $\Sym^g(\Sigma)$, and the Heegaard diagrams used in
their definition.)

When $b_1(Y)>0$, these groups have some extra structure. First of all,
the groups themselves depend on an additional choice of a coherent
orientation system, but there is a canonical such orientation system
fixed in Theorem~\ref{HolDiskTwo:thm:HFinfTwist} of~\cite{HolDiskTwo},
so that is the one we will use (unless otherwise specified). Second,
there is an action by the exterior algebra $\Wedge^* H_1(Y;\Z)/\Tors$,
graded so that $H_1(Y;\Z)/\Tors$ decreases degree by one. Third, there
are variants of these homology groups, but with coefficients twisted
by an arbitrary $\Z[H^1(Y;\Z)]$-module $M$, denoted
$\uHFa(Y,\spinct,M)$, $\uHFm(Y,\spinct,M)$, $\uHFinf(Y,\spinct,M)$,
and $\uHFp(Y,\spinct,M)$. These are all modules over the group-ring
$\Z[H^1(Y;\Z)]$. (When $M$ is not the trivial
module $\Z$ over $\Z[H^1(Y;\Z)]$, the corresponding Floer modules no
longer support an action by $\Wedge^* H_1(Y;\Z)/\Tors$.)

There are two natural long exact sequences whose existence is guaranteed
immediately
from the definitions of these groups 
\begin{equation}
\label{eq:HFaExactSequence}
\begin{CD}
... @>>>\HFa(Y,\spinct)@>{\widehat\iota}>>\HFp(Y,\spinct)@>{U}>> \HFp(Y,\spinct)@>>> ...
\end{CD}
\end{equation}
and 
\begin{equation}
\label{eq:HFinfExactSequence}
\begin{CD}
... @>>>\HFm(Y,\spinct)@>{\iota}>>\HFinf(Y,\spinct)@>{\pi}>> \HFp(Y,\spinct)@>>> ...
\end{CD}
\end{equation}
(with analogues in the twisted case, as well).

Sometimes, we abbreviate these long exact sequences simply by
$\HF^\circ(Y,\spinct)$. Both long exact sequences are also topological
invariants of $Y$. We make this topological invariance
statement more precise, by organizing the results
from~\cite{HolDisk} and \cite{HolDiskTwo} to prove the following:

\begin{theorem}
\label{thm:Naturality}
If $(\Sigma,\alphas,\betas,z)$ and $(\Sigma',\alphas',\betas',z')$ are
equivalent Heegaard diagrams which are strongly admissible for the $\SpinC$
structure $\spinct$, then there are induced isomorphisms of
corresponding long exact sequences: $$
\begin{CD}
...
@>>>\HFm(\alphas,\betas,\spinct)@>{\iota}>>
\HFinf(\alphas,\betas,\spinct) @>{\pi}>>
\HFp(\alphas,\betas,\spinct)
@>>>...	\\
&& @V{\Psi^-_{\spinct}}VV
@V{\Psi^\infty_{\spinct}}VV
@V{\Psi^+_{\spinct}}VV \\
...@>>>\HFm(\alphas',\betas',\spinct)@>{\iota'}>>
\HFinf(\alphas',\betas',\spinct) @>{\pi'}>>
\HFp(\alphas',\betas',\spinct)
@>>>...	
\end{CD}
$$
(i.e. where each square commutes), where the vertical maps commute with the 
actions of $\Z[U]\otimes \Wedge^*(H_1(Y;\Z)/\Tors)$;
and also
$$
\begin{CD}
...
@>>>\HFa(\alphas,\betas,\spinct)@>{\widehat\iota}>>
\HFp(\alphas,\betas,\spinct) @>{U}>>
\HFp(\alphas,\betas,\spinct)
@>>>...	\\
&& @V{{\widehat \Psi}_{\spinct}}VV
@V{\Psi^+_{\spinct}}VV
@V{\Psi^+_{\spinct}}VV \\
...@>>>\HFa(\alphas',\betas',\spinct)@>{\widehat\iota'}>>
\HFp(\alphas',\betas',\spinct) @>{U'}>>
\HFp(\alphas',\betas',\spinct)
@>>>...,
\end{CD}
$$
which commutes with the action of $\Wedge^* (H_1(Y;\Z)/\Tors)$.
Moreover,  the maps ${\widehat\Psi}$, $\Psi^-$,
$\Psi^\infty$, and $\Psi^+$ are uniquely determined up to an overall
factor of $\pm 1$.
\end{theorem}

We return to a proof of Theorem~\ref{thm:Naturality} in
Subsection~\ref{subsec:ProofNaturality}. (See also a twisted version
in Subsection~\ref{subsec:TwistedNaturality}.)

A {\em pointed Heegaard triple $(\Sigma,\alphas,\betas,\gammas,z)$} is
defined using three sets of attaching circles $\alphas$, $\betas$, and
$\gammas$, with a reference point $z\in\Sigma$ disjoint from all these
attaching circles. By counting holomorphic triangles in
$\Sym^g(\Sigma)$ with boundary conditions specified by the three
subspaces $\Ta$, $\Tb$, and $\Tc$, we get maps between tensor products
of Floer homology groups. On the other hand, the Heegaard triple and
the homotopy classes of triangles both can be interpreted in terms of
four-dimensional data. We briefly digress to recall this
four-dimensional interpretation (referring the reader to
Section~\ref{HolDiskTwo:sec:HolTriangles} of~\cite{HolDiskTwo} for
more details).

\subsection{Homotopy classes of triangles, and $\SpinC$ structures on four-manifolds}
\label{subsec:TopPrelim}

Given a pointed Heegaard triple we have three $g$-dimensional tori
$\Ta=\alpha_1\times...\alpha_g$, $\Tb=\beta_1\times...\times\beta_g$
and $\Tc=\gamma_1\times...\gamma_g$ which are embedded in the $g$-fold
symmetric product of $\Sigma$, $\Sym^g(\Sigma)$.

Clearly, the Heegaard triple specifies three three-manifolds
$Y_{\alpha,\beta}$, $Y_{\beta,\gamma}$, and $Y_{\alpha,\gamma}$ with
Heegaard diagrams given by using the corresponding pairs of $g$-tuples
of attaching circles. Indeed, the Heegaard triple specifies a
four-manifold $X_{\alpha,\beta,\gamma}$ which bounds these three three-manifolds.
We recall the construction presently.

Let $\Delta$ denote the two-simplex, with vertices $v_{\alpha},
v_{\beta}, v_{\gamma}$ labeled clockwise, and let $e_{i}$ denote the
edge $v_{j}$ to $v_{k}$, where $\{i,j,k\}=\{\alpha,\beta,\gamma\}$. 
Then, we form the identification space
$$X_{\alpha,\beta,\gamma}=\frac{\left(\Delta\times \Sigma\right)
\coprod \left(e_\alpha \times U_\alpha\right) \coprod
\left(e_\beta\times U_\beta \right)\coprod
\left(e_\gamma\times U_\gamma \right)}
{ \left(e_{\alpha}\times \Sigma\right) \sim \left(e_\alpha \times
\partial U_\alpha\right), \left(e_{\beta}\times \Sigma\right) \sim
\left(e_\beta \times \partial U_\beta\right), \left(e_{\gamma}\times
\Sigma\right) \sim \left(e_\gamma \times \partial U_\gamma\right) }.
$$ Over the vertices of $\Delta$, this space has corners, which can be
naturally smoothed out to obtain a smooth, oriented, four-dimensional
cobordism between the three-manifolds $Y_{\alpha,\beta}$,
$Y_{\beta,\gamma}$, and $Y_{\alpha,\gamma}$ as claimed.
Indeed, under the natural orientation conventions implicit in the above description,
$$\partial X_{\alpha,\beta,\gamma}=-Y_{\alpha,\beta}-Y_{\beta,\gamma}+Y_{\alpha,\gamma}.$$

Recall that the two-dimensional homology of $X$ corresponds to the set
of {\em triply-periodic domains} for the Heegaard triples -- these are
two-chains in $\Sigma$ (which can be represented by maps $\Phi\colon
F\longrightarrow \Sigma$ where $F$ is a two-dimensional
manifold-with-boundary) with multiplicity zero at the reference point
$z$, and whose boundaries are linear combinations of curves appearing in the
tuples $\alphas$, $\betas$, and $\gammas$. Such chains in turn
correspond uniquely to the integer relations they give amongst the homology
classes coming from the $\alphas$, $\betas$, and $\gammas$. Given such
a two-chain in $\Sigma$, the associated homology class $H(\PerDom)$ in
$H_2(X;\Z)$ can be concretely thought of as follows. We can extend
$\Phi\colon F\longrightarrow \Sigma$, to a map $${\widehat \Phi}\colon
{\widehat F}\longrightarrow X,$$ by viewing $\Phi$ as a map into some
copy of $\Sigma$ times an interior point $x\in\Delta$, and then
attaching cylinders connecting $\partial F$ to a collection of
$\alpha$, $\beta$, and $\gamma$-circles occuring in copies of $\Sigma$
on the $\alpha$, $\beta$ and $\gamma$-edges of $\Delta\times \Sigma$.
Then, we cap off these boundary components in the corresponding
$U_\alpha$, $U_\beta$ and $U_\gamma$ handlebodies, to get a map from a
closed surface into $X$. This closed surface represents the homology
class $H(\PerDom)$.

Fix intersection points $\x\in\Ta\cap \Tb$, 
$\y\in\Tb\cap \Tc$, $\w\in\Ta\cap\Tc$. A map
$$u \colon \Delta \longrightarrow \Sym^{g}(\Sigma) $$
with the boundary conditions that 
$u(v_{\gamma})=\x$, $u(v_{\alpha})=\y$, and $u(v_{\beta})=\w$, and 
$u(e_{\alpha})\subset \Ta$, $u(e_{\beta})\subset \Tb$, 
$u(e_{\gamma})\subset 
\Tc$ is called a {\em Whitney triangle connecting $\x$, 
$\y$, and $\w$}. There is a naturally defined map $$\epsilon\colon
(\Ta\cap\Tb)\times (\Tb\cap\Tc)\times (\Ta\cap\Tc)\longrightarrow
H_1(X;\Z)$$ with the property that $\epsilon(\x,\y,\w)=0$ if and only
if there is a Whitney triangle connecting $\x$, $\y$, and
$\w$. Two Whitney triangles
are said to be homotopic if the maps are homotopic through maps which
are all Whitney triangles. For fixed $\x$, $\y$, and $\w$, 
let $\pi_{2}(\x,\y,\w)$ denote the space of
homotopy classes of Whitney triangles connecting $\x$, $\y$, and $\w$.
Any two elements of this space $\pi_2(\x,\y,\w)$ have a naturally
associated difference, which is a periodic domain (after we subtract
off a sufficient multiple of $\Sigma$).  When $g>1$, this sets up an affine
isomorphism, whenever $\pi_2(\x,\y,\w)$ is non-empty, 
$$\pi_2(\x,\y,\w)\cong \Z \oplus
H_2(X_{\alpha,\beta,\gamma};\Z),$$ where the first factor is given by
the local multiplicity at the reference point $z$. 
When $g=1$, and $\pi_2(\x,\y,\w)$ is non-empty, then there is an isomorphism
$$\pi_2(\x,\y,\w)\cong H_2(X_{\alpha,\beta,\gamma};\Z).$$

As explained in Section~\ref{HolDiskTwo:sec:HolTriangles}
of~\cite{HolDiskTwo}, a Whitney triangle can be used to construct a
singular two-plane field in $X$ whose underlying $\SpinC$ structure is
independent of the choices made in its construction, giving rise to a
map $$\spinc_z\colon
\pi_2(\x,\y,\w)\longrightarrow
\SpinC(X),$$
where $\spinc_z(\psi)$ restricts to $\spinc_z(\x)$, $\spinc_z(\y)$,
and $\spinc_z(\w)$ at the three boundary components. 

There is a splicing map
$$\pi_2(\x,\y,\w)\times \pi_2(\x',\x)\times \pi_2(\y',\y)\times \pi_2(\w,\w')
\longrightarrow \pi_2(\x',\y',\w'),$$
which we denote simply by addition.
We have the following result from~\cite{HolDiskTwo}:

\begin{lemma}
Fix $\x,\x'\in \Ta\cap\Tb$, $\y,\y'\in\Tb\cap\Tc$, and
${\mathbf v},{\mathbf v}'\in\Ta\cap\Tc$.  
Two homotopy
classes $\psi\in
\pi_2(\x,\y,{\mathbf v})$ and 
$\psi'\in\pi_2(\x',\y',{\mathbf v}')$ induce the same $\SpinC$
structure if and only if there are homotopy classes
$\phi_1\in\pi_2(\x',\x)$, $\phi_2\in\pi_2(\y',\y)$, and
$\phi_3\in\pi_2({\mathbf v},{\mathbf v}')$ with
$$\psi'=\psi+\phi_1+\phi_2+\phi_3.$$
\end{lemma}

\subsection{Holomorphic triangles}
\label{subsec:HolTriangles}

We recall the maps induced by the holomorphic triangle construction.
Except, of course, for the interpretation of Whitney triangles using
four-manifolds, the results we recall here are modeled on
corresponding constructions in Lagrangian Floer homology.

We begin with a Heegaard triple $(\Sigma,\alphas,\betas,\gammas,z)$,
and fix a $\SpinC$ structure $\spinc$ over $X_{\alpha,\beta,\gamma}$. 
Under suitable admissibility hypotheses, there are chain maps:
$$f^\circ(\cdot,\spinc)\colon \CFc(\alphas,\betas,\spinc_{\alpha,\beta})
\otimes \CFleq(\betas,\gammas,\spinc_{\beta,\gamma}) \longrightarrow 
\CFc(\alphas,\gammas,\spinc_{\alpha,\gamma}),$$
where $\spinc_{\xi,\eta}$ is the restriction of $\spinc$ to $Y_{\xi,\eta}$
and where $\CFleq\cong\CFm$ is the chain complex generated by pairs $[\x,j]$ 
with $j\leq 0$, given by the formula:
$$
f^{\circ}([\x,i]\otimes[\y,j];\spinc)
=\sum_{\w\in\Ta\cap\Tc}
\sum_{\{\psi\in\pi_2(\x,\y,\w)\big|\spinc_z(\psi)=\spinc,\Mas(\psi)=0\}}
\Big(\#\ModFlow(\psi)\Big)\cm {[\w,i+j-n_z(\psi)]},
$$ where $\ModFlow(\psi)$ denotes the moduli space of
pseudo-holomorphic triangles in the homotopy class $\psi$, and
$\Mas(\psi)$ denotes its expected dimension.  There is also a variant
$$f^\leq(\cdot,\spinc)\colon
\CFleq(\alphas,\betas,\spinc_{\alpha,\beta})
\otimes \CFleq(\betas,\gammas,\spinc_{\beta,\gamma})
\longrightarrow \CFleq(\alphas,\gammas,\spinc_{\alpha,\gamma})$$
defined analogously.

The admissibility hypothesis relevant for the above maps is the
{\em $\spinc$-strong admissibility} for the Heegaard triple
(Definition~\ref{HolDiskTwo:def:AdmissibleTriple} of~\cite{HolDiskTwo}).
Specifically, we require that for each non-trivial triply-periodic
domain which can be written as a sum of doubly-periodic domains
$$\PerDom=\cald_{\alpha,\beta}+\cald_{\beta,\gamma}+\cald_{\alpha,\gamma}$$ with the property
that $$\langle c_1(\spinc_{\alpha,\beta}), H({\mathcal
D}_{\alpha,\beta})
\rangle + \langle c_1(\spinc_{\beta,\gamma}), H(\cald_{\beta,\gamma})
\rangle + \langle c_1(\spinc_{\alpha,\gamma}), H(\cald_{\alpha,\gamma})
\rangle =2n\geq 0,$$
we have some local multiplicity of $\PerDom$ which is strictly greater
than $n$.

Before stating the relevant invariance properties of these
constructions, we recall some of the notation involved for
isotopies. In Section~\ref{HolDiskTwo:sec:Special}
of~\cite{HolDiskTwo}, we achieved admissibility by using {\em special
isotopies} -- isotopies which never cross the basepoint $z$, and which
are realized as exact Hamiltonian isotopies in $\Sigma$. Now, if
$(\Sigma,\alphas,\betas,z)$ and $(\Sigma,\alphas',\betas,z)$ are
strongly $\spinc$-admissible Heegaard diagrams, and there is a special
isotopy connecting $\alphas$ to $\alphas'$, then the special isotopy
induces an isomorphism $$\Gamma_{\alpha',\alpha;\beta}\colon
\HFc(\alphas,\betas,\spinc)\longrightarrow
\HFc(\alphas',\betas,\spinc)$$
defined by counting holomorphic disks with time-dependent boundary
conditions; similarly, a special isotopy from $\betas$ to $\betas'$ induces an isomorphism
$$\Gamma_{\alpha;\beta,\beta'}\colon
\HFc(\alphas,\betas,\spinc)\longrightarrow
\HFc(\alphas,\betas',\spinc)$$

\begin{theorem}
\label{thm:Triangles}
Suppose that $(\Sigma,\alphas,\betas,\gammas,z)$ is
an admissible Heegaard triple for the $\SpinC$ structure
$\spinc$. Then, the
induced maps on homology $$F^\circ(\cdot,\spinc)\colon
\HF^\circ(\alphas,\betas,\spinc_{\alpha,\beta})\otimes
\HFleq(\betas,\gammas,\spinc_{\beta,\gamma})\longrightarrow
\HF^\circ(\alphas,\gammas,\spinc_{\alpha,\gamma}),$$
satisfy the following properties.
\begin{itemize}
\item $F^\circ(\cdot,\spinc)$ is independent of the choices of complex structures and perturbations used in the definition of the moduli spaces of triangles.
\item $F^\circ(\cdot,\spinc)$ 
commutes with the action by $\Z[U]$.
\item $F^\circ(\cdot,\spinc)$ is invariant under special isotopies, in the sense that
if $$\Gamma_{\alpha;\beta,\beta'}
\colon \HFc(\alphas,\betas,\spinc_{\alpha,\beta})\longrightarrow
\HFc(\alphas,\betas',\spinc_{\alpha,\beta'})$$
and $$\Gamma_{\alpha',\alpha;\beta}\colon
\HFc(\alphas,\betas,\spinc_{\alpha,\beta})\longrightarrow
\HFc(\alphas',\betas,\spinc_{\alpha',\beta})$$ 
(where $\spinc_{\alpha,\beta'}=\spinc_{\alpha,\beta}=\spinc_{\alpha',\beta}$
denotes the restrictions of $\spinc$ to the boundary component)
are the isomorphisms induced by
some isotopies from the $\alphas$ to the $\alphas'$, and $\betas$ to
$\betas'$, then the following diagrams commute: 
$$\begin{CD}
\HF^\circ(\alphas,\betas,\spinc_{\alpha,\beta})\otimes 
\HFleq(\betas,\gammas,\spinc_{\beta,\gamma}) 
@>{F^\circ(\cdot,\spinc)}>> 
\HF^\circ(\alphas,\gammas,\spinc_{\alpha,\gamma}) \\
@V{\Gamma_{\alpha',\alpha;\beta}\otimes\Id}VV
@VV{\Gamma_{\alpha',\alpha;\gamma}}V \\
\HF^\circ(\alphas',\betas,\spinc_{\alpha,\beta})\otimes 
\HFleq(\betas,\gammas,\spinc_{\beta,\gamma}) 
@>{F^\circ(\cdot,\spinc)}>> \HF^\circ(\alphas',\gammas,\spinc_{\alpha,\gamma})
\end{CD}$$
and
$$\begin{CD}
\HF^\circ(\alphas,\betas,\spinc_{\alpha,\beta})\otimes 
\HFleq(\betas,\gammas,\spinc_{\beta,\gamma}) 
@>{F^\circ(\cdot,\spinc)}>> \HF^\circ(\alphas,\gammas) \\
@V{\Gamma_{\alpha;\beta,\beta'}\otimes\Gamma_{\beta',\beta,\gamma}}VV
@VV{\Id}V \\
\HF^\circ(\alphas,\betas',\spinc_{\alpha,\beta'})\otimes 
\HFleq(\betas',\gammas,\spinc_{\beta',\gamma}) 
@>{F^\circ(\cdot,\spinc)}>> \HF^\circ(\alphas,\gammas).
\end{CD}
$$
\end{itemize}
\end{theorem}

The proof can be found in~\cite{HolDiskTwo}: the existence is shown in
Theorem~\ref{HolDiskTwo:thm:HolTriangles} of~\cite{HolDiskTwo};
independence of complex structures and isotopy invariance are respectively
Propositions~\ref{HolDiskTwo:prop:TrianglesJIndep} and
\ref{HolDiskTwo:prop:TriangleIsotopyInvariance} from that paper.

The holomorphic triangle construction also satisfies an associativity
property which we state separately.  
For the associativity, we use
pointed Heegaard quadruples
\begin{equation}
\label{eq:CoboundaryVanishes}
(\Sigma,\alphas,\betas,\gammas,\deltas,z),
\end{equation}
to which we can associate a 
four-manifold $X_{\alpha,\beta,\gamma,\delta}$ defined using a square
rather than a triangle. 
We assume that
the four-manifold $X_{\alpha,\beta,\gamma,\delta}$ satisfies the
additional hypothesis that 
\begin{eqnarray}
\label{eq:GoodDecomposition}
\delta H^1(Y_{\beta,\delta})|Y_{\alpha,\gamma}=0
&{\text{and}}&
\delta H^1(Y_{\alpha,\gamma})|Y_{\beta,\delta}=0.
\end{eqnarray}
In this case, one can formulate a strong admissibility
hypothesis, which depends on a 
$\delta H^1(Y_{\beta,\delta})+\delta H^1(Y_{\alpha,\gamma})$-orbit ${\mathfrak S}$ 
of a $\SpinC$ structure 
in $\SpinC(X_{\alpha,\beta,\gamma,\delta})$ (see
Section~\ref{HolDiskTwo:subsec:Assoc} of~\cite{HolDiskTwo} for the
definition). The topological hypotheses guarantee that strong
$\spinc$-admissibility can always be achieved for a Heegaard quadruple.

\begin{remark}
The topological hypothesis of Equation~\eqref{eq:CoboundaryVanishes}
is automatically satisfied if in our Heegaard quadruple, there are two
consecutive $g$-tuples of circles which span the same subspace of
$H_1(\Sigma;\Z)$. For instance, if the subspaces spanned by $\gammas$
and $\deltas$ coincide, then we can express each
$(\betas,\deltas)$-periodic domain as a sum of $(\gammas,\deltas)$ and
$(\betas,\gammas)$-periodic domains (which ensures that $\delta
H^1(Y_{\beta,\delta})=0$; and similarly, $(\alphas,\gammas)$-periodic
domains can be expressed as a sum of $(\gammas,\deltas)$- and
$(\deltas,\alphas)$-periodic domains. In fact, this guarantees that
$\delta H^1(Y_{\beta,\delta})+\delta H^1(Y_{\alpha,\gamma})=0$.
\end{remark}

\begin{theorem}
\label{thm:Associativity}
Let $(\Sigma,\alphas,\betas,\gammas,\deltas,z)$ be a
pointed Heegaard quadruple which is strongly
${\mathfrak S}$-admissible, where ${\mathfrak S}$ is
a $\delta H^1(Y_{\beta,\delta})+ \delta
H^1(Y_{\alpha,\gamma})$-orbit in
$\SpinC(X_{\alpha,\beta,\gamma,\delta})$.
Then, we have
\begin{eqnarray*}
\lefteqn{\sum_{\spinc\in{\mathfrak S}}
\Fc{\alpha,\gamma,\delta}
(\Fc{\alpha,\beta,\gamma}
(\xi_{\alpha,\beta}\otimes
\theta_{\beta,\gamma}; 
\spinc_{\alpha,\beta,\gamma})\otimes 
\theta_{\gamma,\delta};\spinc_{\alpha,\gamma,\delta})} \\
&=&\sum_{\spinc\in{\mathfrak S}}
\Fc{\alpha,\beta,\delta}(\xi_{\alpha,\beta}\otimes 
\Fleq{\beta,\gamma,\delta}
(\theta_{\beta,\gamma}\otimes\theta_{\gamma,\delta};\spinc_{\beta,\gamma,\delta});
\spinc_{\alpha,\beta,\delta}),
\end{eqnarray*}
where
$\xi_{\alpha,\beta}\in \HF^\circ(\alphas,\betas,\spinc_{\alpha,\beta})$, 
$\theta_{\beta,\gamma}$ and $\theta_{\gamma,\delta}$ lie in
$\HFleq(\betas,\gammas,\spinc_{\beta,\gamma})$ 
and $\HFleq(\gammas,\deltas,\spinc_{\gamma,\delta})$ 
respectively.
\end{theorem}

The above is a special case of
Theorem~\ref{HolDiskTwo:thm:Associativity} of~\cite{HolDiskTwo}.

There is a variant of the chain maps $\Fc{\alpha,\beta,\gamma}$ which incorporates an
action of the first homology of
$X_{\alpha,\beta,\gamma}$.

First of all, we claim that the map $$H_1(Y_{\alpha,\gamma}\coprod
Y_{\beta,\gamma}\coprod Y_{\alpha,\gamma};\Z) /\Tors \longrightarrow
H_1(X_{\alpha,\beta,\gamma};\Z)/\Tors$$ is surjective. Given $h\in
H_1(X_{\alpha,\beta,\gamma};\Z)/\Tors$, we let 
$$ (h_{\alpha,\beta},h_{\beta,\gamma},h_{\alpha,\gamma})\in 
H_1(Y_{\alpha,\gamma}\coprod
Y_{\beta,\gamma}\coprod Y_{\alpha,\gamma};\Z) /\Tors $$
be an element which maps to it. We then define
\begin{eqnarray}
\lefteqn{\Fc{\alpha,\beta,\gamma}(h\otimes \xi_{\alpha,\beta}\otimes \theta_{\beta,\gamma})
=}
\nonumber \\
&& 
\Fc{\alpha,\beta,\gamma}(\left(h_{\alpha,\beta}\cm \xi_{\alpha,\beta}\right)\otimes 
\theta_{\beta,\gamma}) + 
\Fc{\alpha,\beta,\gamma}(\xi_{\alpha,\beta}\otimes 
\left(h_{\beta,\gamma}\cm 
\theta_{\beta,\gamma}\right))
- h_{\alpha,\gamma} \cm \Fc{\alpha,\beta,\gamma}(\xi_{\alpha,\beta}\otimes \theta_{\beta,\gamma}),
\label{eq:DefHOneAction}
\end{eqnarray}
where the actions on the right-hand-side all represent the actions of $H_1$ 
of the three-manifolds on their corresponding Floer groups.

\begin{lemma}
\label{lemma:ExtendedTriangles}
The action defined in Equation~\eqref{eq:DefHOneAction} induces a map
$$
\Fc{\alpha,\beta,\gamma}(\cdot,\spinc)\colon \Wedge^* 
\left(H_1(X_{\alpha,\beta,\gamma};\Z)/\Tors\right)
\otimes \HFc(\alphas,\betas,\spinc_{\alpha,\beta})
\otimes \HFleq(\betas,\gammas,\spinc_{\beta,\gamma})
\longrightarrow
\HFc(\alphas,\gammas,\spinc_{\alpha,\gamma}).~\mbox{\hskip3in}
$$
\end{lemma}

\begin{proof}
Let $X=X_{\alpha,\beta,\gamma}$. We have natural isomorphisms
\begin{eqnarray*}
H_1(Y_{\alpha,\beta})/\Tors &\cong &
\Hom(H_2(Y_{\alpha,\beta}),\Z) \\
H_1(X)/\Tors &\cong &
\Hom(H^3(X,\partial X),\Z),
\end{eqnarray*}
where all homology groups above (and throughout this proof) are
understood to use integral coefficients.  Recall that $H_3(X)=0$, and
$H_2(X)$ has no torsion. Thus, we can dualize a portion of the
Mayer-Vietoris sequence for the pair $(X,\partial X)$ to give: $$
\begin{CD}
\Hom(H_2(X),\Z)@>>>\Hom(H_2(\partial X),\Z) 
@>>> \Hom(H_3(X,\partial X),\Z) @>>> 0,
\end{CD}
$$
in particular establishing surjectivity of the map on $H_1/\Tors$.

We must now verify that the action by
$(h_{\alpha,\beta},h_{\beta,\gamma},h_{\alpha,\gamma})$ is trivial on
the image of $\Hom(H_2(X),\Z)$.
Now, we have identifications:
\begin{eqnarray*}
\Hom(H_2(Y_{\alpha,\beta}),\Z)&\cong &
H^1(\Ta)\oplus H^1(\Tb)/H^1(\Sym^{g}(\Sigma))\\
\Hom(H_2(X),\Z) &\cong&
H^1(\Ta)\oplus H^1(\Tb) \oplus H^1(\Tc)/H^1(\Sym^g(\Sigma)).
\end{eqnarray*}
In particular, those elements of $\Hom(H_2(\partial X),\Z)$ which come
from $\Hom(H_2(X),\Z)$ can be represented by three constraints, one
in each of $H^1(\Ta)$, $H^1(\Tb)$, and $H^1(\Tc)$. 
By pushing the
constraints out to the ends of the triangles as in
Proposition~\ref{HolDiskTwo:prop:CommutesWithHOneAction}
of~\cite{HolDiskTwo}, it follows that the action 
triples
$(h_{\alpha,\beta},h_{\beta,\gamma},h_{\alpha,\gamma})$ coming from
$\Hom(H_2(X),\Z)$ is trivial. 

We have just seen how $\Wedge^1 (H_1(X)/\Tors)$ acts. This extends to
an action of the full exterior algebra: this follows from the
corresponding fact for the Floer homologies themselves.  Similarly,
independence of the action from the auxiliary data (paths of complex
structures, isotopies, etc.) follows from the corresponding facts for
the three bounding three-manifolds, as well.
\end{proof}

\subsection{Preliminaries on $\#^n(S^1\times S^2)$.}
Let $\spinc_0$ be the $\SpinC$ structure over
$\#^n(S^1\times S^2)$ whose first Chern class vanishes.
It is implicit in
Section~\ref{HolDiskOne:sec:HandleSlides} of~\cite{HolDisk}, and
spelled out in Section~\ref{HolDiskTwo:subsec:SOneTimesSTwo}
of~\cite{HolDiskTwo} that
\begin{eqnarray*}
\HFa(\#^n(S^1\times S^2),\spinc_0)&\cong& \wedge^* H^1(\#^n(S^1\times S^2);\Z),\\
\HFp(\#^n(S^1\times S^2),\spinc_0)&\cong& \wedge^* H^1(\#^n(S^1\times S^2);\Z)\otimes_\Z \Z[U^{-1}],\\
\HFm(\#^n(S^1\times S^2),\spinc_0)&\cong& \wedge^* H^1(\#^n(S^1\times S^2);\Z)\otimes_\Z \Z[U] \\
\HFinf(\#^n(S^1\times S^2),\spinc_0)&\cong& \wedge^* H^1(\#^n(S^1\times S^2);\Z)\otimes_\Z \Z[U,U^{-1}],
\end{eqnarray*}
where $U^{-1}$ has grading $2$.

In particular, there is a ``top-dimensional'' homology group of
$\HFa(\#^n(S^1\times S^2))$, which is isomorphic to $\Z$; let
${\widehat\Theta}_n$ denote one of its generators. The image of this
element in $\HFp(\#^n(S^1\times S^2))$ is denoted $\Theta_n^+$. A
generator for the top-dimensional homology group of
$\HFm(\#^n(S^1\times S^2),\spinc_0)\subset \HFinf(\#^n(S^1\times
S^2),\spinc_0)$ is denoted $\Theta^-_n$, while a generator for the
top-dimensional homology group of $\HFleq(\#^n(S^1\times S^2))$ is
denoted $\Theta_n$.  We drop the subscript $n$ when it is clear from
the context. Observe that all these elements are uniquely specified up
to sign.

\begin{defn}
\label{def:MaxInt}
Let $(\Sigma,\alphas,\betas,z)$ be a Heegaard diagram for
$\#^n(S^2\times S^1)$ which is strongly admissible for $\spinc_0$, the
$\SpinC$ structure whose first Chern class vanishes. We call an
intersection point $\x\in\Ta\cap\Tb$ {\em maximal} if 
$\spinc_z(\x)=\spinc_0$ and the relative grading of $\x$
(thought of as an element of
$\CFa(\#^n(S^2\times S^1))$) agrees with the relative grading of
${\widehat \Theta}_n\in\HFa(\#^n(S^2\times S^1))$.
\end{defn}

We can identify explicit representatives for the generators of
$\HFa(\#^n(S^1\times S^2),\spinc_0)$, as follows.

\begin{defn}
\label{def:StandardHeegaard}
Consider the pointed Heegaard diagram
$$(\Sigma,\{\alpha_1,...,\alpha_n\},\{\beta_1,...,\beta_n\},z),$$
where $\Sigma$ is an oriented two-manifold of genus $n$;
$\alpha_i$ is a small isotopic translate of $\beta_i$ (via
an isotopy supported in the complement of the basepoint $z$), meeting
it in two canceling transverse intersection points. This Heegaard
diagram is called a {\em standard Heegaard  diagram} for
$\#^n(S^2\times S^1)$.
\end{defn}

For a standard Heegaard decomposition of $(S^2\times S^1)$,
the tori $\Ta\cap\Tb$ meet in $2^n$ intersection points, corresponding
to the generators of 
$$\HFa(\#^n(S^1\times S^2),\spinc_0)\cong 
\wedge^* H^1(\#^n(S^1\times S^2);\Z).$$
In particular, there is a unique maximal intersection point, inducing
a representative of ${\widehat\Theta}_n$.

\subsection{Proof of Naturality}
\label{subsec:ProofNaturality}

We now turn to a proof of Theorem~\ref{thm:Naturality}.

\begin{defn}
Two Heegaard diagrams are said to be {\em strongly equivalent} if they
differ by a sequence of isotopies and handleslides.
\end{defn}

\begin{lemma}
\label{lemma:ConnectOrdered}
Let $(\Sigma_1,\alphas_1,\betas_1,z_1)$ and
$(\Sigma_2,\alphas_2,\betas_2,z_2)$ be two equivalent Heegaard
diagrams, then both can be stabilized to give Heegaard diagrams
$(\Sigma_1',\alphas_1',\betas_1',z_1')$ and
$(\Sigma_2',\alphas_2',\betas_2',z_2')$ which are strongly equivalent.
\end{lemma}

\begin{proof}
According to Proposition~\ref{HolDiskOne:prop:HeegaardMoves}
of~\cite{HolDisk}, it follows that
$(\Sigma_1,\alphas_1,\betas_1,z_1)$ and
$(\Sigma_2,\alphas_2,\betas_2,z_2)$ can be connected by a sequence of
Heegaard moves.

We claim that isotopies commute with stabilizations, in the following
sense.  Suppose that one can get from $(\Sigma,\alphas,\betas,z)$ to
$(\Sigma',\alphas',\betas',z')$ by a sequence of isotopies and handleslides
followed by
stabilizations, then one can also get from the first Heegaard diagram
by first a sequence of stabilizations, followed by isotopies and
handleslides. This follows from the simple observation that it is
equivalent to either
\begin{itemize}
\item isotope some $\alpha_i$ (resp. $\beta_i$) across a point
$w$ and then stabilize at $w$, 
\item or to first stabilize at $w$ and then
handleslide $\alpha_i$ (resp. $\beta_i$) twice over the newly
introduced $\alpha_g$ (resp. $\beta_g$). 
\end{itemize}

In this way, stabilizations can be moved before isotopies and
handleslides. They can also be moved before destabilizations:
destabilizing a canceling pair of curves $\alpha_g$ and $\beta_g$,
and then stabilizing at some point $w$ is the same as first
stabilizing at $w$ and then destabilizing the original pair of curves.

Thus, we can commute all stabilizations to the beginning of the sequence of moves, and
all destabilizations to the end.
This is equivalent to the assertion made in the lemma.
\end{proof}

Recall (c.f. Section~\ref{HolDiskOne:sec:DefHF} of~\cite{HolDisk}, see
especially Theorem~\ref{HolDiskOne:thm:IndepCxStruct}) that the Floer
homology groups depend on a choice of complex structure ${\mathfrak
j}$, and a generic path $J_s\subset {\mathcal U}$ of 
perturbations of the induced complex structure over
$\Sym^g(\Sigma)$, where ${\mathcal U}$ is a contractible
set of almost-complex structures, defined in~\cite{HolDisk}.

\begin{lemma}
Fix two different choices $({\mathfrak j}, J_s)$ and $(\mathfrak j',J_s')$
of complex structures and perturbations. Then, 
the induced isomorphism between Floer homologies 
$$\HF^\circ_{J_{s}}(\alphas,\betas,\spinc)\cong
\HF^\circ_{J_s'}(\alphas,\betas,\spinc)$$
are canonical, i.e. independent of paths connecting $J_s$ and $J_s'$. 
\end{lemma}

\begin{proof}
First, we suppose that ${\mathfrak j}={\mathfrak j}'$, writing
$J_s=J_s(0)$ and $J_s'=J_s(1)$. We can connect $J_s(0)$ and $J_s(1)$
by a one-parameter family $J_{s,t}$ as in the proof of
Theorem~\ref{HolDiskOne:thm:IndepCxStruct} of~\cite{HolDisk}
to construct a chain homotopy equivalence 
$$\Phi^\infty_{J_{s,t}}\colon
(\CFinf(\alphas,\betas,\spinc),\partial^\infty_{J_s(1)})
\longrightarrow 
(\CFinf(\alphas,\betas,\spinc),\partial^\infty_{J_s(0)}).$$ In fact, we claim
that the induced map on homology by $\Phi^{\infty}_{J_{s,t}}$ is
actually independent of the path of paths $J_{s,t}$ used in its
definition (i.e.depending only on its endpoints).  Since ${\mathcal U}$
is contractible, any two such paths $J_{s,t}(0)$ and $J_{s,t}(1)$ can
be connected by a homotopy $J_{s,t,\tau}$. We use this homotopy to
define a map $$H^\infty_{J_{s,t,\tau}}([\x,i]) =
\sum_\y\sum_{\{\phi\in\pi_2(\x,\y)|\Mas(\phi)=-1\}}
\# \left(\ModFlow_{J_{s,t,\tau}}(\phi)\right)\cdot[\y,i-n_z(\phi)].$$
which lowers degree by $1$. Now, counting the ends of the moduli spaces
$\ModFlow_{J_{s,t,\tau}}(\psi)$ with $\Mas(\psi)=0$, we see that
$H^\infty_{J_{s,t,\tau}}$ provides a chain homotopy between 
$\Phi^\infty_{J_{s,t}(0)}$ and $\Phi^\infty_{J_{s,t}(1)}$.

Since $\Phi^{\infty}_{J_{s,t}}$ respects the filtration induced by $n_z$, 
the analogous results hold for $\CFpm$ as well.

The boundary maps are invariant under small perturbations (provided
that the perturbations still give energy bounds), so -- in view of the
fact that the space of complex structures ${\mathfrak j}$ over $\Sigma$ is
contractible -- it follows also that the chain complexes are naturally
identified as ${\mathfrak j}$ is varied, as well. (Observe here that a codimension
two subset of complex structures was removed from consideration
in~\cite{HolDisk} to avoid difficulties occurring in moduli spaces with
Maslov index two, but these do not arise in the definition of the
boundary maps.)
\end{proof}

Suppose now that $(\Sigma,\alphas_1,\betas_1,z)$ and
$(\Sigma,\alphas_2,\betas_2,z)$ are strongly equivalent Heegaard diagrams both of which are
admissible for $\spinc$. Following
Section~\ref{HolDiskTwo:sec:Special} of~\cite{HolDiskTwo},
we can construct another pointed
Heegaard diagram $(\Sigma,\alphas_1',\betas_1',z)$ so that:
\begin{itemize}
\item the curves $\alphas_1'$ and $\betas_1'$ are connected to $\alphas_1$ and $\betas_1$ 
respectively by special isotopies
\item the Heegaard quadruple $(\Sigma,\alphas_1',\betas_1',\alphas_2,\betas_2,z)$ 
is strongly admissible,
for the unique $\SpinC$ structure over $X_{\alpha_1',\beta_1',\alpha_2,\beta_2}$
whose restriction to $(\Sigma,\alphas_1',\betas_2,z)\cong Y$ is $\spinc$ and 
whose restriction to $(\Sigma,\alphas_1',\alphas_2,z)\cong 
(\Sigma,\betas_1',\betas_2,z)\cong \#^{g}(S^1\times S^2)$ is $\spinc_0$.
\end{itemize}

We have an isomorphism induced by the special isotopies
$$\Gamma_{\alpha_1',\alpha_1;\beta_1,\beta_1'}\colon \CFc(\alphas_1,\betas_1,\spinc)
\longrightarrow \CFc(\alphas_1',\betas_1',\spinc)$$
defined by counting holomorphic disks with time-dependent boundary
conditions (analogous to those maps $\Gamma_{\alpha',\alpha;\beta}$
discussed earlier, only now we use time-dependent boundary conditions
on both boundaries of the strip).
We define the strong equivalence map $\Phi$ to be the
composite $$\begin{CD}
\CF^{\circ}(\alphas_1,\betas_1,\spinc) @>{\Gamma_{\alpha_1',\alpha_1;\beta_1,\beta_1'}}>>
\CF^{\circ}(\alphas_1',\betas_1',\spinc)
@>{\Theta_{\alpha_2,\alpha_1'}\otimes ~\cdot~ \otimes \Theta_{\beta_1',\beta_2}}>>
\CF^{\circ}(\alphas_2,\betas_2,\spinc)
\end{CD},$$
where the last map is shorthand for the map
$$\xi \mapsto \Fc{\alpha_2,\beta_1',\beta_2}
(\Fc{\alpha_2,\alpha_1',\beta_1'}(\Theta_{\alpha_2,\alpha_1'}\otimes \xi,\spinc)\otimes
\Theta_{\beta_1',\beta_2},\spinc),$$
where 
$$\Fc{\alpha_2,\alpha_1',\beta_1'}\colon \CFleq(\alphas_2,\alphas_1',\spinc_0)\otimes
\CFc(\alphas_1',\betas_1',\spinc)\longrightarrow \CFc(\alphas_2,\betas_1',\spinc)$$ 
is a version of the map associated to holomorphic triangles, where we use
$\CFleq$ on the first factor (rather than the second, as usual). 
Of course, the strong equivalence map $\Phi$ depends on the auxiliary Heegaard diagram
$(\Sigma,\alphas_1',\betas_1',z)$, as well as the isotopies connecting
$\alphas_1$ to $\alphas_1'$ and $\betas_1$ to $\betas_1'$.

To show that the map induced on homology is independent of the
various choices, we use some naturality properties of the maps induced
by isotopies. 

\begin{lemma}
\label{lemma:NaturalGammas}
The maps $\Gamma$ are natural under composition of
isotopies. Specifically, if we fix isotopies from $\alphas$ to
$\alphas'$, and isotopies from $\alphas'$ to $\alphas''$ (which never
cross the reference point $z$ in a Heegaard diagram
$(\Sigma,\alphas,\betas,z)$), then the maps on homology
$\Gamma_{\alpha'',\alpha;\beta}$ induced by the composite of the two
isotopies equals the composition
$\Gamma_{\alpha'',\alpha';\beta}\circ\Gamma_{\alpha',\alpha;\beta}$.
Similarly, if we have an isotopies of $\alphas$ to $\alphas'$ and
$\betas$ to $\betas'$, then (on homology),
$$\Gamma_{\alpha',\alpha;\beta,\beta'}=\Gamma_{\alpha',\alpha;\beta'}\circ
\Gamma_{\alpha;\beta,\beta'}=\Gamma_{\alpha';\beta,\beta'}\circ
\Gamma_{\alpha',\alpha;\beta}.$$
\end{lemma}

\begin{proof}
This result is a variation on the argument that the
$\Gamma_{\alpha;\beta,\beta'}$ is an isomorphism on homology, as
described in Theorem~\ref{HolDiskOne:thm:Isotopies} of~\cite{HolDisk} (which in
turn is a variation on familiar arguments from Floer theory). We focus
on the first statement: suppose that $\Psi_t'$ is an isotopy from
$\alphas$ to $\alphas'$ and $\Psi_t''$ is an isotopy from $\alphas'$ to
$\alphas''$. We consider holomorphic disks with time-dependent
boundary conditions, given by juxtaposing $\Psi_t$ and $\Psi_t'$, with 
a ``long gap'' in between: i.e. we define a one-parameter family of isotopies
$$\Psi_{\tau,t}=\left\{\begin{array}{ll}
\Psi'_{t-\tau} & {\text{if $t\geq 0$}} \\
\Psi''_{t+\tau} & {\text{if $t\leq 0$}}
\end{array}
\right.$$
We can count points in $\Psi_{\tau}$-time-dependent moduli spaces.
Taking the limit as $\tau\goesto\infty$, this can be used to construct
the chain homotopy between
$\Gamma_{\alpha'',\alpha';\beta}\circ\Gamma_{\alpha';,\alpha;\beta}$
and $\Gamma_{\alpha'',\alpha;\beta}$ as required. The other assertions
follow similarly.
\end{proof}

\begin{lemma}
\label{lemma:DefPhi}
Up to sign, the map on homology induced by a strong equivalence
$$\Phi\colon \HF^\circ(\alphas_1,\betas_1,\spinc)\longrightarrow
\HF^\circ(\alphas_2,\betas_2,\spinc)$$
is well-defined.
\end{lemma}

\begin{proof}
We must show that $\Phi$ is independent of the intermediate Heegaard
diagram and the isotopy. This follows from the following commutative
diagram:
$$\begin{CD}
\HF^{\circ}(\alphas_1,\betas_1,\spinc) 
@>{\Gamma_{\alpha_1',\alpha_1;\beta_1,\beta_1'}}>>
\HF^{\circ}(\alphas_1',\betas_1',\spinc)
@>{\Theta_{\alpha_2,\alpha_1'}\otimes ~\cdot~ \otimes \Theta_{\beta_1',\beta_2}}>>
\HF^{\circ}(\alphas_2,\betas_2,\spinc) \\
@V{=}VV @V{\Gamma_{\alpha_1'',\alpha_1';\beta_1',\beta_1''}}VV @V{=}VV \\
\HF^{\circ}(\alphas_1,\betas_1,\spinc) 
@>{\Gamma_{\alpha_1'',\alpha_1;\beta_1,\beta_1''}}>>
\HF^{\circ}(\alphas_1'',\betas_1'',\spinc)
@>{\Theta_{\alpha_2,\alpha_1''}\otimes ~\cdot~ \otimes \Theta_{\beta_1'',\beta_2}}>>
\HF^{\circ}(\alphas_2,\betas_2,\spinc)
\end{CD}.$$
Commutativity of the left-hand square follows from the naturality of
the $\Gamma$-maps under juxtaposition of isotopies
(Lemma~\ref{lemma:NaturalGammas} above). Commutativity of the
right-hand square up to sign follows from the isotopy invariance of
the triangle construction
(c.f. Lemma~\ref{HolDiskTwo:prop:TriangleIsotopyInvariance}
of~\cite{HolDiskTwo}), together with the observation that
\begin{eqnarray}
\Gamma_{\beta_1',\beta_1'';\beta_2}
(\Theta_{\beta_1'',\beta_2})=
\pm \Theta_{\beta_1',\beta_2}, &{\text{and}}&
\Gamma_{\alpha_2;\alpha_1'',\alpha_1'}(\Theta_{\alpha_2,\alpha_1''})
=\pm \Theta_{\alpha_2,\alpha_1'}.
\label{eq:Thetas}
\end{eqnarray}
since both sides of both equations represent
generators of a groups which are isomorphic to $\Z$. 
For the reader's convenience, we break 
the verification of the commutativity of this right-hand
square down into the following steps,
each of which is an application of either the isotopy invariance 
of the triangle construction (the last two commutative diagrams in 
Theorem~\ref{thm:Triangles}), the naturality of the isotopy maps 
(Lemma~\ref{lemma:NaturalGammas}), or Equations~\eqref{eq:Thetas}:
\begin{eqnarray*}
\lefteqn{F^\circ_{\alpha_2,\beta_1'',\beta_2}
(F^\circ_{\alpha_2,\alpha_1'',\beta_1''}(
\Theta_{\alpha_2,\alpha_1''}\otimes 
\Gamma_{\alpha_1'',\alpha_1';\beta_1',\beta_1''}(\xi))
\otimes \Theta_{\beta_1'',\beta_2})} \\
&=& \pm F^\circ_{\alpha_2,\beta_1',\beta_2}
(\Gamma_{\alpha_2;\beta_1'',\beta_1'}
(F^\circ_{\alpha_2,\alpha_1'',\beta_1''}(
\Theta_{\alpha_2,\alpha_1''}\otimes 
\Gamma_{\alpha_1'',\alpha_1';\beta_1',\beta_1''}(\xi)))
\otimes \Gamma_{\beta_1',\beta_1'';\beta_2}(\Theta_{\beta_1'',\beta_2}) \\
&=& 
\pm F^\circ_{\alpha_2,\beta_1',\beta_2}
(F^\circ_{\alpha_2,\alpha_1'',\beta_1'}(
\Theta_{\alpha_2,\alpha_1''}\otimes 
\Gamma_{\alpha_1'';\beta_1'',\beta_1'}\circ
\Gamma_{\alpha_1'',\alpha_1';\beta_1',\beta_1''}(\xi)))
\otimes \Gamma_{\beta_1',\beta_1'';\beta_2}
(\Theta_{\beta_1'',\beta_2})) \\
&=& 
\pm F^\circ_{\alpha_2,\beta_1',\beta_2}
(F^\circ_{\alpha_2,\alpha_1'',\beta_1'}(
\Theta_{\alpha_2,\alpha_1''}\otimes 
\Gamma_{\alpha_1'',\alpha_1';\beta_1'}(\xi))
\otimes \Gamma_{\beta_1',\beta_1'';\beta_2}
(\Theta_{\beta_1'',\beta_2})) \\
&=& 
\pm F^\circ_{\alpha_2,\beta_1',\beta_2}
(F^\circ_{\alpha_2,\alpha_1'',\beta_1'}(
\Theta_{\alpha_2,\alpha_1''}\otimes 
\Gamma_{\alpha_1'',\alpha_1';\beta_1'}(\xi))
\otimes \Theta_{\beta_1',\beta_2}) \\
&=& 
\pm F^\circ_{\alpha_2,\beta_1',\beta_2}
(F^\circ_{\alpha_2,\alpha_1',\beta_1'}(
\Gamma_{\alpha_2;\alpha_1'',\alpha_1'}(\Theta_{\alpha_2,\alpha_1''})\otimes 
\xi)
\otimes 
\Theta_{\beta_1',\beta_2}) \\
&=& 
\pm F^\circ_{\alpha_2,\beta_1',\beta_2}
(F^\circ_{\alpha_2,\alpha_1',\beta_1'}(
\Theta_{\alpha_2,\alpha_1'}\otimes 
\xi)
\otimes 
\Theta_{\beta_1',\beta_2}).
\end{eqnarray*}
(Here, when we go from the second to the third equation, we use a
variant of the first commutative diagram in Theorem~\ref{thm:Triangles},
only varying $\gamma$ rather than $\alpha$.)
\end{proof}

Suppose now that $(\Sigma^s,\alphas^s,\betas^s,z)$ is obtained from
$(\Sigma,\alphas,\betas,z)$ by a stabilization; i.e. we
let $E$ be an oriented two-manifold of genus one with a pair
$\alpha_{g+1}$ and $\beta_{g+1}$ of embedded curves meeting
transversally in a single intersection point $c$, and 
let $\Sigma^s=\Sigma\# E$, $\alphas^s=\alphas\cup\{\alpha_{g+1}\}$,
and $\betas^s=\betas\cup \{\beta_{g+1}\}$. Then we
have defined an isomorphism belonging to the stabilization
(c.f. Section~\ref{HolDiskOne:sec:Stabilization} of~\cite{HolDisk})
$$\sigma\colon \HF^\circ(\alphas,\betas,\spinc)\longrightarrow
\HF^\circ(\alphas^s,\betas^s,\spinc),$$
which (for appropriately chosen complex structures and perturbations) is 
given by the map which associates to each 
intersection point $\x\in\Sym^{g}(\Sigma)$, the stabilized
point 
$\x\times \{c\}\in\Sym^{g+1}(\Sigma\# E)$.

If we stabilize two strongly equivalent Heegaard diagrams, the resulting diagrams are
strongly equivalent, as well. We wish to show that the isomorphism induced by 
strong equivalence commutes with the isomorphism induced by stabilization. 
To this end, we shall employ the gluing theorem for holomorphic triangles,
Theorem~\ref{HolDiskTwo:thm:GlueTriangles} of~\cite{HolDiskTwo}, which
we state here for convenience:

\begin{theorem}
\label{thm:GlueTriangles}
Fix a pair of Heegaard diagrams 
\begin{eqnarray*}
(\Sigma,\alphas,\betas,\gammas,z) &{\text{and}}&
(E,\alpha_0,\beta_0,\gamma_0,z_0), 
\end{eqnarray*}
where $E$ is a Riemann surface of genus one. Consider the 
connected sum $\Sigma\# E$, where the connected sum points are near
the distinguished points $z$ and $z_0$ respectively. Fix intersection
points $\x,\y,\w$ for the first diagram and a class
$\psi\in\pi_2(\x,\y,\w)$, and intersection points $x_0$, $y_0$, and
$w_0$ for the second, with a triangle $\psi_0\in\pi_2(x_0,y_0,w_0)$
with $\Mas(\psi)=\Mas(\psi_0)=0$.  Suppose moreover that
$n_{z_0}(\psi_0)=0$.  Then, for a suitable choice of complex structures
and perturbations, we
have a diffeomorphism of moduli spaces: $$\ModFlow(\psi')\cong
\ModFlow(\psi)\times
\ModFlow(\psi_0),$$
where $\psi'\in \pi_2(\x\times x_0, \y\times y_0,\w\times w_0)$ is 
the triangle for $\Sigma\# E$ whose domain on the $\Sigma$-side 
agrees with $\cald(\psi)$, and whose domain on the $E$-side agrees with
$\cald(\psi_0)+n_{z}(\psi)[E]$. 
\end{theorem}

\begin{lemma}
\label{lemma:StabilizationEquivalence}
The isomorphism induced by stabilization commutes with that induced
by equivalence.
\end{lemma}

\begin{proof}
Let $\alpha_{g+1}'$ and $\beta_{g+1}'$ be small Hamiltonian isotopic
translates of $\alpha_{g+1}$ and $\beta_{g+1}$ in $E$. Also, let
$\alphas_1'$ and $\betas_1'$ be isotopic translates of the $\alphas_1$
and $\betas_1$ required in the definition of $\Phi$, chosen so that
the isotopies are supported away from the stabilization point. We claim then that
the following diagram commutes:
$$
\begin{CD}
\HF^\circ(\alphas_1,\betas_1)@>{\Gamma_{\alpha_1',\alpha_1;\beta_1}}>>
\HF^\circ(\alphas_1',\betas_1') \\
@V{\sigma_1}VV @VV{\sigma_1'}V \\
\HF^\circ(\alphas_1^s,\betas_1^s)
@>{\Gamma_{{\alpha_1'}^s,\alpha_1^s;\beta_1^s,{\beta_1'}^s}}>>
\HF^\circ({\alphas_1'}^s,{\betas_1'}^s),
\end{CD}
$$ where here $\alphas_1^s$, $\betas_1^s$, ${\alphas_1'}^s$, and
${\betas_1'}^s$ denote the stabilizations
$\alphas_1\cup\{\alpha_{g+1}\}$, $\betas_1\cup\{\beta_{g+1}\}$,
$\alphas_1'\cup\{\alpha_{g+1}'\}$, and $\betas_1'\cup\{\beta_{g+1}'\}$
respectively. This follows from the analogue of the gluing theorem for
stabilization invariance (Theorem~\ref{HolDiskOne:thm:StabilizeHFb}
of~\cite{HolDisk}), only using
time-dependent flow-lines.

Moreover, we have commutativity of 
$$
\begin{CD}
\HF^\circ(\alphas_1',\betas_1')@>{\Theta_{\alpha_2,\alpha_1'}\otimes~\cdot~\otimes
\Theta_{\beta_1',\beta_2}}>> 
\HF^\circ(\alphas_2,\betas_2) \\
@V{\sigma_1'}VV @VV{\sigma_2}V \\
\HF^\circ({\alphas_1'}^s,{\betas_1'}^s)
@>{\Theta_{\alpha_2^s,{\alpha_1'}^s}\otimes
~\cdot~\otimes\Theta_{{\beta_1'}^s,\beta_2^s}}>> 
\HF^\circ(\alphas_2^s,\betas_2^s) \\
\end{CD}
$$
This is follows from an application of Theorem~\ref{thm:GlueTriangles}
(see also Lemma~\ref{lemma:Stabilize} below).

Together, these two commutative squares show that stabilization commutes with the map
associated to strong equivalence.
\end{proof}

\vskip.2cm 
\noindent{\bf{Proof of Theorem~\ref{thm:Naturality}}}
Suppose that $(\Sigma_1,\alphas_1,\betas_1,z_1)$ and
$(\Sigma_2,\alphas_2,\betas_2,z_2)$ represent the same
three-manifold $Y$, and both are admissible for
$\spinc\in\SpinC(Y)$. Then, we connect them as in
Lemma~\ref{lemma:ConnectOrdered} and define $\Psi^\circ$ 
to be the composite
of the stabilization isomorphism composed with the equivalence
isomorphism, composed with the inverse of the final stabilization
isomorphism. 

We must argue that the map $\Psi^\circ$ is independent of the choice of
common stabilization $\Sigma_1'=\Sigma_2'$. Now, if $\Sigma_1''$ is
another intermediate choice, then there is a third stabilization
$\Sigma_1'''$ which is obtained by stabilizing both $\Sigma_1'$ and
$\Sigma_1''$. The map induced by factoring through $\Sigma_1'''$
agrees by that induced through $\Sigma_1'$ (and also $\Sigma_1''$) in view of 
Lemma~\ref{lemma:StabilizationEquivalence}.

Observe that the equivalence isomorphisms $\Psi^\circ$ are equivariant under
the action of $\Wedge^* H_1(Y,\spinc)/\Tors$, as was verified in 
Section~\ref{HolDiskTwo:sec:DefHF} of~\cite{HolDiskTwo}.

The fact that $\Psi^\circ$ induces a map of the long exact sequences
for $\HFa$ (Equation~\eqref{eq:HFaExactSequence}) is a formal
consequence of the fact that the isomorphism induced by isotopies
carries $\CFa(\alphas,\betas,\spinct)\subset
\CFp(\alphas,\betas,\spinct)$ to the corresponding $\CFa$ for the
isotopic Heegaard diagram, together with the fact that the map
$F^+(\cm,\spinc)$ induced by counting holomorphic triangles is
$U$-equivariant.
\vskip.2cm

\subsection{Twisted coefficients and naturality}
\label{subsec:TwistedNaturality}

Let $H=H^1(Y;\Z)$. There is a Floer homology with
twisted coefficients, denoted $\uHFc(Y,\spinct)$, which is a module
over the ring $\Z[H]$. We recall the construction briefly.

Recall that there is always a natural map from $\pi_2(\x,\x)$ to $H$,
which is obtained as follows. Each $\phi\in\pi_2(\x,\x)$ naturally
gives rise to an associated two-chain in $\Sigma$ whose boundary is a
collection of circles among the $\alphas$ and $\betas$. We can then
close off the two-chain to give a closed two-cycle in $Y$ by gluing on
copies of the attaching disks for the handlebodies in the Heegaard
diagram for $Y$. The Poincar\'e dual of this two-cycle is the
associated element of $H^1(Y;\Z)$. (When $g=1$, this map is an
isomorphism, when $g>1$, it, together with the multiplicity at $z$,
$n_z$, give a natural identification $\pi_2(\x,\x)\cong \Z\oplus H$.)

Fix a $\spinct$-admissible Heegaard diagram $(\Sigma,\alphas,\betas,z)$
for $Y$, and an {\em additive assignment} in the sense of
Section~\ref{HolDiskTwo:subsec:TwistedCoeffs} of~\cite{HolDiskTwo},
i.e. letting ${\mathfrak S}\subset \Ta\cap\Tb$ be the set of
intersection points $\x$ so that $\spinc_z(\x)=\spinct$, we fix a
collection of maps: $$A=\{A_{\x,\y}\colon
\pi_2(\x,\y)\longrightarrow H^1(Y;\Z)\}_{\x,\y\in{\mathfrak S}},$$
so that:
\begin{itemize}
\item
when $\x=\y$, $A_{\x,\x}$ is the canonical  map from 
$\pi_2(\x,\x)$ onto $H^1(Y;\Z)$ defined above
\item $A$ is
compatible with splicing in the sense that if
$\x,\y,\w\in{\mathfrak S}$, then for each $\phi_1\in\pi_2(\x,\y)$ and
$\phi_2\in\pi_2(\y,\w)$, we have that
$A(\phi_1*\phi_2)=A(\phi_1)+A(\phi_2)$.
\end{itemize}

We write elements of the group-ring $\Z[H]$ as finite sums $\sum_{h\in
H} n_h\cm e^h$ (with $n_h\in\Z$).  Consider the chain complex
$\uCFinf(Y,\spinct,A)$ which is freely generated as a
$\Z[H^1(Y;\Z)]$-module by elements $[\x,i]$ where $\x\in{\mathfrak S}$
and $i\in\Z$. The boundary map is defined by $$\partial^\infty[\x,i] =
\sum_{\y\in{\mathfrak S}} \sum_{\{\phi\in\pi_2(\x,\y)\big|\Mas(\phi)=1\}}
\left(\#\UnparModFlow(\phi)\right)\cm e^{A(\phi)} \otimes [\y,i-n_z(\phi)].$$
The homology groups of this complex are the completely twisted homology groups
$\uHFc(Y,\spinct)$. More generally, if $M$ is any module over
$\Z[H]$, then 
$\uHF(Y,\spinct,M)$ is defined to be the homology of 
$$\uCFc(Y,\spinct,M,A)\cong \uCFc(Y,\spinct,A)\otimes_{Z[H]}M.$$

Observe that the homology is independent of the additive assignment
$A$.  Specifically, if we have two additive assignments $A$ and $A'$,
then we can define an isomorphism of chain complexes $$\phi\colon
\uCFinf(Y,\spinct,M,A)\longrightarrow \uCFinf(Y,\spinct,M,A')$$ as
follows. Fix a point $\x_0\in{\mathfrak S}$, and let 
$$\phi(m \otimes[\x,i])
=(m\cm e^{A(\phi)-A'(\phi)})\otimes [\x,i],$$ where $\phi$ is
any element of $\pi_2(\x_0,\x)$. It is easy to see that $\phi$ is an
isomorphism of chain complexes (over $\Z[H]$); but the actual
isomorphism depends up to translation by an element of $H$
on the initial point
$\x_0\in{\mathfrak S}$. Indeed, the homology groups are topological
invariants, according to the following:

\begin{theorem}
\label{thm:TwistedNaturality}
Fix a module $M$ for $\Z[H]$, and suppose that 
$(\Sigma,\alphas,\betas,z)$ and $(\Sigma',\alphas',\betas',z')$ are
equivalent Heegaard diagrams which are admissible for
the $\SpinC$ structure $\spinct$, then there are induced isomorphisms
of the corresponding long exact sequences: $$
\begin{CD}
\uHFm(\alphas,\betas,\spinct,M)@>{\iota}>>
\uHFinf(\alphas,\betas,\spinct,M) @>{\pi}>>
\uHFp(\alphas,\betas,\spinct,M)
	\\
@V{{\underline\Psi}^-}VV
@V{{\underline \Psi}^\infty}VV
@V{{\underline \Psi}^+}VV \\
\uHFm(\alphas',\betas',\spinct,M)@>{\iota'}>>
\uHFinf(\alphas',\betas',\spinct,M) @>{\pi'}>>
\uHFp(\alphas',\betas',\spinct,M)
\end{CD}
$$ (i.e. where each square commutes), where the vertical maps commute
with the actions of $\Z[U]$.  Moreover, the maps ${\underline\Psi}^-$,
${\underline\Psi}^\infty$, and ${\underline\Psi}^+$ are uniquely
determined up to multiplication by $\pm 1$ and translation in $H$.
There is a similar canonically-induced map of the long exact sequence
for $\HFa$. 
\end{theorem}

The above theorem is proved the same way as Theorem~\ref{thm:Naturality} is proved; with 
suitable modifications made to the holomorphic triangle construction which allow for
twisted coefficients, and induced modules, as we describe in the next subsection.

Of course, the chain complex $\uCFc(Y,\spinct,M)$ is obtained from the
chain complex in the totally twisted case $\uCFc(Y,\spinct)$ by a
change of coefficients; thus, the corresponding homology groups are
related by a universal coefficients spectral sequence
(c.f.~\cite{EilenbergMacLane}) . In particular, when $M$ is the trivial
$\Z[H^1(Y;\Z)]$-module $\Z$, the $M$-twisted chain complex is the same
as the untwisted chain complex stated earlier. Moreover, it is easy to see that
$$\Tor^i_{\Z[H^1(Y;\Z)]}(\Z,\Z)\cong \Wedge^i(H_1(Y;\Z)/\Tors).$$
Thus, the universal
coefficients spectral sequence has $E_2$ term
$\Wedge^i(H_1(Y;\Z)/\Tors)\otimes_{\Z}\uHFc(Y,\spinct)$, and $E^\infty$
term $\HFc(Y,\spinct)$.

\subsection{Holomorphic triangles and twisted coefficients}

We set this up with slightly less generality than
in~\cite{HolDiskTwo}, but with sufficient generality for the
applications in the present paper.

Let $W$ be a cobordism from $Y_1$ to $Y_2$, and fix a module $M$ for $\Z[H^1(Y_1;\Z)]$.
The group
$$K(W)=\Ker\left(H^2(W,\partial W;\Z)\longrightarrow H^2(W;\Z)\right)$$
can be used to induce a module $M(W)$ for $\Z[H^1(Y_2;\Z)]$, defined by
$$M(W)=M\otimes_{\Z[H^1(Y_1;\Z)]}\Z[K(W)].$$

Suppose now that $(\Sigma,\alphas,\betas,\gammas,z)$ is a Heegaard
triple, and suppose that $(\Sigma,\betas,\gammas,z)$ represents
$\#^n(S^2\times S^1)$, so that when we fill it in with $\#^n(D^3\times
S^1)$, we obtain the cobordism $W$ above. Fix a $\SpinC$ structure
$\spinc$ over $W$, and let ${\mathfrak S}$ denote the set of homotopy
classes of triangles which induce the fixed $\SpinC$ structure
$\spinc$ over $X_{\alpha,\beta,\gamma}$. 

Now, two additive assignments $A_{\alpha,\beta}$, $A_{\beta,\gamma}$,
$A_{\alpha,\gamma}$ for the corresponding Heegaard diagrams on the
boundary and a single choice of $\psi_0\in{\mathfrak S}$ give rise to
a map $$A_W\colon {\mathfrak S}\longrightarrow K(W),$$ defined as
follows. If $\psi\in\pi_2(\x,\y,\w)\cap {\mathfrak S}$, then we can
find paths $\phi_{\alpha,\beta}\in\pi_2(\x_0,\x)$,
$\phi_{\beta,\gamma}\in\pi_2(\y_0,\y)$ and
$\phi_{\alpha,\gamma}\in\pi_2(\w_0,\w)$ with the property that
$$\psi=\psi_0+\phi_{\alpha,\beta}+\phi_{\beta,\gamma}+\phi_{\alpha,\gamma}.$$
Then, we define
$$A_W(\psi)=\delta(A(\phi_{\alpha,\beta})\oplus
A(\phi_{\alpha,\gamma})),$$ where $\delta$ denotes the coboundary map
$\delta\colon H^1(\partial W;\Z)\longrightarrow K(W)$.
It is easy to see that this is independent of the choices of 
$\phi_{\alpha,\beta}$, $\phi_{\beta,\gamma}$ and $\phi_{\alpha,\gamma}$ 
as above.

We now claim that the
holomorphic triangle construction gives a map
$$\uFc_{\alpha,\beta,\gamma}\colon
\uHFc(\alphas,\betas,\spinc|Y_1,M)\otimes
\HFleq(\betas,\gammas,\spinc_0) \longrightarrow \uHFc(\alphas,\gammas,\spinc|Y_2,M(W)),$$
to be the map on homology induced by the chain map $\ufc$ defined by:
\begin{eqnarray*}
\lefteqn{\ufc_{\alpha,\beta,\gamma}(m\otimes [\x,i]\otimes [\y,j];\spinc)
=}  \nonumber 
\\
&& \sum_{\w\in\Ta\cap\Tb} 
\sum_{\{\psi\in\pi_2(\x,\y,\w)\big| \spinc_z(\psi)=\spinc, \Mas(\psi)=0\}}
\Big(\#\ModFlow(\psi)\Big)
m\otimes e^{A_W(\psi)}
\otimes [\w,i+j-n_z(\psi)].
\label{eq:DefTriangleTwisted}
\end{eqnarray*}
This map commutes with the obvious $\Z[H^1(Y_1;\Z)]$ actions.

Although this construction depends on an initial triangle
$\psi_0\in{\mathfrak S}$, the following is clear:
If $\psi_1\in{\mathfrak S}$ is another such choice, then we can write
$\psi_1 = \psi_0+\phi_{\alpha,\beta}+\phi_{\beta,\gamma}+\phi_{\alpha,\gamma}$
with $A_{\beta,\gamma}(\phi_{\beta,\gamma})=0$. 
Then, if $\uFc_0$ and $\uFc_1$ are the two induced maps, then we have:
$$e^{A_{\alpha,\gamma}(\phi_{\alpha,\gamma})} \cm \uFc_0 = 
\uFc_1 \cm e^{A_{\alpha,\beta}(\phi_{\alpha,\beta})}.$$

Note that if $\gamma\colon M\longrightarrow M'$ is a map of
$\Z[H^1(Y;\Z)]$-modules, then there is an induced map 
$$H(\gamma)\colon \uHFc(Y,M)\longrightarrow \uHFc(Y,M')$$
defined in the obvious way.

The map $\uFc$ is a refinement of the map
$$\Fc{}\colon \HFc(\alphas,\betas,\spinc|Y_1)\otimes \HFleq(\betas,\gammas,\spinc_0)
\longrightarrow \HFc(\alphas,\gammas,\spinc|Y_2)$$
defined earlier, in the following sense.
Observe first that if $\Z$ is a trivial $\Z[H^1(Y;\Z)]$-module, then 
$\uHFc(Y,\spinc,\Z)=\HFc(Y,\spinc)$. Moreover, if $W$ is a cobordism from $Y_1$ to $Y_2$, then the trivial $H^1(Y_1;\Z)$-module $\Z$ induces the $H^1(Y_2;\Z)$-module 
$\Z[H^2(W,Y_2;\Z)]$. Now, letting $\epsilon\colon \Z[H^2(W,Y_2;\Z)]\longrightarrow \Z$
be the natural map (to the trivial $H^1(Y_2;\Z)$-module), we have that
\begin{equation}
\label{eq:TwistAndUntwist}
\Fc{}(\xi\otimes\eta,\spinc)=H(\epsilon)\circ \uFc(\xi\otimes\eta,\spinc).
\end{equation}

\section{Cobordisms}
\label{sec:Cobordisms}

In  Section~\ref{sec:ECobord}, we will use holomorphic triangles
to construct invariants (which are maps on the $\HF^\circ$)
induced by one-, two-, and three-handle addition
to a given three-manifold. By composing these invariants, we define 
the invariant of a cobordism, satisfying the following properties:

\begin{theorem}
\label{thm:CobordismInvariant}
Let $W$ be an oriented, smooth, connected, four-dimensional cobordism
with $\partial W= -Y_1\cup Y_2$. Fix a $\SpinC$ structure $\spinc\in\SpinC(W)$,
and let $\spinct_i$ denote its restriction to $Y_i$.
Then, by composing 
the maps associated to a handle decomposition of $W$, we obtain
a map
$$\Fc{W,\spinc}\colon
\HFc(Y_1,\spinct_1)\longrightarrow \HFc(Y_2,\spinct_2)$$ 
which is a
smooth oriented four-manifold invariant, uniquely defined up to sign.
More generally, $\Fc{W,\spinc}$ can be extended to a map
$$\Fc{W,\spinc}\colon \HFc(Y_1,\spinct_1)\otimes 
\Wedge^* \left(H_1(W,\Z)/\Tors\right)\longrightarrow
\HFc(Y_2,\spinct_2),$$
which is also a four-manifold invariant. More precisely, 
if $\phi\colon W \longrightarrow W'$ is an orientation-preserving  diffeomorphism, then we have a commutative diagram:
$$
\begin{CD}
\HFc(Y_1,\spinct_1)\otimes \Wedge^*(H_1(W;\Z)/\Tors)
@>{F_{W,\spinc}}>> \HFc(Y_2,\spinct_2) \\
@V{(\phi|Y_1)_*\otimes \phi_*}VV 
@V{(\phi|_{Y_2})_*}VV \\
\HFc(Y_1',\spinct_1')\otimes \Wedge^*(H_1(W';\Z)/\Tors)
@>{F_{W',\spinc'}}>> \HFc(Y_2,\spinct_2'),
\end{CD}
$$
where $\spinc=\phi^*(\spinc')$.
\end{theorem}

\begin{remark}
\label{rmk:Naturality}
When we write
$$\Fc{W,\spinc}\colon
\HFc(Y_1,\spinct_1)\longrightarrow \HFc(Y_2,\spinct_2),$$ 
we mean a map between the long exact sequences relating $\HFa$,
$\HFm$,
$\HFinf$, and $\HFp$. Specifically, the above theorem is saying that
there is a collection of maps $\Fa{W,\spinc}$, $\Fm{W,\spinc}$, $F^\infty_{W,\spinc}$
and $\Fp{W,\spinc}$ for which the following diagrams commute:
$$
\begin{CD}
... @>>>\HFm(Y_1,\spinct_1)@>{\iota}>>\HFinf(Y_1,\spinct_1)@>{\pi}>> \HFp(Y_1,\spinct_1)@>>> ... \\
&& @V{\Fm{W,\spinc}}VV @V{\Finf{W,\spinc}}VV @V{\Fp{W,\spinc}}VV \\
... @>>>\HFm(Y_2,\spinct_2)@>{\iota}>>\HFinf(Y_2,\spinct_2)@>{\pi}>> \HFp(Y_2,\spinct_2)@>>> ...\\
\end{CD}
$$
and also
$$
\begin{CD}
... @>>>\HFa(Y_1,\spinct_1)@>{\widehat\iota}>>
\HFp(Y_1,\spinct_1)@>{U}>> \HFp(Y_1,\spinct_1)@>>> ... \\
&& @V{\Fa{W,\spinc}}VV @V{\Fp{W,\spinc}}VV @V{\Fp{W,\spinc}}VV \\
... @>>>\HFa(Y_2,\spinct_2)
@>{\widehat\iota}>>
\HFp(Y_2,\spinct_2)@>{U}>> \HFp(Y_2,\spinct_2)@>>> ...\\
\end{CD}
$$
\end{remark}

This invariant enjoys a number of fundamental properties.

\begin{theorem}{\bf{(Finiteness)}}
\label{thm:Finiteness}
Let $W$ be a cobordism from $Y_1$ to $Y_2$, and fix $\SpinC$
structures $\spinct_1$ and $\spinct_2$ over $W$. Then, for each
$\xi\in\HFp(Y_1,\spinct_1)$, there are only finitely many
$\spinc\in\SpinC(W)$ for which $$ \Fp{W,\spinc}(\xi)\neq 0.$$
Moreover, for each integer $d$, and each element
$\eta\in\HFm(Y_1,\spinct_1)$, there are only finitely many $\SpinC$
structures $\spinc\in\SpinC(W)$ for which 
$$\Fm{W,\spinc}(\eta)\not\in U^d \HFm(Y_2,\spinct_2).$$
\end{theorem}

\begin{theorem} {\bf{(Composition Law)}}
\label{thm:Composition}
Let $W_1$ and $W_2$ be a pair of connected cobordisms with $\partial
W_1 = -Y_1\cup Y_2$, $\partial W_2 = -Y_2\cup Y_3$, and let
$W=W_1\cup_{Y_2} W_2$ be their composite. Fix $\SpinC$ structures
$\spinc_i\in\SpinC(W_i)$ for $i=1,2$ with $\spinc_1|_{Y_2}=\spinc_2|_{Y_2}$. 
Then, $\Fc{W_2,\spinc_2}\circ \Fc{W_1,\spinc_1}$ can be written as a sum:
$$\Fc{W_2,\spinc_2}\circ \Fc{W_1,\spinc_1}
=
\sum_{\{\spinc\in\SpinC(W)\big|\spinc|W_1=\spinc_1,
\spinc|W_2=\spinc_2\}} 
\pm \Fc{W,\spinc}.$$
\end{theorem}

There is a natural bilinear pairing   (c.f. Section~\ref{sec:Duality})
$$\langle,\rangle \colon
\CFinf(Y,\spinc)\otimes \CFinf(-Y,\spinc)\longrightarrow \Z$$ (c.f. Section~\ref{sec:Duality}) which
descends to give pairings on
$\HFinf(Y,\spinc)\otimes \HFinf(-Y,\spinc)$, 
$\HFp(Y,\spinc)\otimes \HFm(-Y,\spinc)$, and $\HFpRed(Y,\spinc)\otimes
\HFmRed(-Y,\spinc)$. 

\begin{theorem} 
\label{thm:Duality}
{\bf (Duality)} Let $W$ be a cobordism from $Y_1$ to $Y_2$. Then, the
map induced by $W$, thought of as a cobordism from $Y_1$ to $Y_2$, is dual 
to the map induced by $W$, thought of as a cobordism from $-Y_2$ to $-Y_1$;
i.e. for each  $\xi\in\HFp(Y_1,\spinct_1)$,
$\eta\in\HFm(-Y_2,\spinct_2)$, we have that $$\langle
F^+_{W,\spinc}(\xi),\eta\rangle_{\HFp(Y_2,\spinct_2)\otimes \HFm(-Y_2,\spinct_2)} = \langle
\xi,F^-_{W,\spinc}(\eta)\rangle_{\HFp(Y_1,\spinct_1)\otimes \HFm(-Y_1,\spinct_1)}.$$
\end{theorem}

In general, there is a $\Zmod{2}$ action on the set of $\SpinC$
structures over a given manifold. In three-dimensions, thinking of
$\SpinC$ structures as equivalence classes of nowhere-vanishing vector
fields $v$, this action is induced by the map $v\mapsto -v$. In four
dimensions, thinking of $\SpinC$ structures as equivalence classes of
almost-complex structures $J$ defined away from a finite collection of
points, the map is induced by the map $J\mapsto -J$.

In~\cite{HolDisk}, we defined an isomorphism of chain complexes, which
induces an identification $\HFc(Y,\spinc)$ with $\HFc(Y,{\overline\spinc})$.
We denote the isomorphism by
$\Conjugate_{Y}\colon \HFc(Y,\spinc)\longrightarrow\HFc(Y,{\overline \spinc})$.

\begin{theorem}\label{thm:Conjugation} {\bf{(Conjugation invariance)}}
Let $W$ be a cobordism as above. Then, 
$\Fc{W,\spinc}=\Conjugate_{Y_2}\circ \Fc{W,{\overline \spinc}}
\circ \Conjugate_{Y_1}$.
\end{theorem}

Proofs of duality and conjugation invariance will be given in
Section~\ref{sec:Duality}.

\begin{theorem}{\bf{(Blow-up formula)}}
\label{thm:BlowUp}
Let $W$ be the cobordism $\left([0,1]\times Y\right)\#\mCP$ (internal
connected sum). Then,
for each $\SpinC$ structure $\spinc$ over $Y$, we have that
$$\Fc{W,\spinc}\colon \HFc(Y,\spinct)\longrightarrow \HFc(Y,\spinct)$$
is multiplication by $U^{\frac{\ell\cm(\ell+1)}{2}}$ 
where $\spinc\in\SpinC W$ is characterized by 
$c_1(\spinc)|\{0\}\times Y=\spinct$ and
$\langle c_1(\spinc), E\rangle = \pm (2\ell+1)$, for $\ell\geq 0$.
\end{theorem}

\subsection{Twisted analogues}

We give the following twisted versions of
Theorems~\ref{thm:CobordismInvariant} and \ref{thm:Composition}.

Fix first a module $M$ over $\Z[H^1(Y_1;\Z)]$, and a cobordism $W$
from $Y_1$ to $Y_2$ as in Theorem~\ref{thm:CobordismInvariant}. 
Let $M(W)$ be the
induced module $M(W)$ over $\Z[H^1(Y_2;\Z)]$ induced by the cobordism
as in Subsection~\ref{subsec:TwistedNaturality}.

We have the following:

\begin{theorem}
\label{thm:CobordismInvariantTwist}
Fix a $\SpinC$ structure
$\spinc$ over $W$.
Then, there is an associated map
$$\Fc{W,\spinc,M}\colon 
\uHFc(Y_1,\spinc|Y_1,M)
\longrightarrow 
\uHFc(Y_2,\spinc|Y_2,M(W))
$$
which is uniquely defined up to multiplication by $\pm 1$,
left-translation by an element
of $H^1(Y_1;\Z)$, and right translation by an element of $H^1(Y_2;\Z)$
\end{theorem}

We let $[\Fc{W,\spinc}]$ be the $H^1(Y_1;\Z)\oplus H^1(Y_2;\Z)$-orbit 
of the map constructed in Theorem~\ref{thm:CobordismInvariantTwist}.

We state a refined version of the composition law. Let $W_1$ be a
cobordism from $Y_1$ to $Y_2$ and $W_2$ be a cobordism from $Y_2$ to
$Y_3$, and let $W=W_1\cup_{Y_2} W_2$ be their composite.  Then there
is a natural inclusion $$i_M\colon M(W)\longrightarrow M(W_1)(W_2), $$
which has a canonically defined right-inverse $$\Pi_M\colon
M(W_1)(W_2)\longrightarrow M(W).$$ To construct this, notice that
there is an $H^1(Y_2;\Z)$-equivariant inclusion $$
i\colon K(W)\longrightarrow
\frac{K(W_1)\oplus K(W_2)}{H^1(Y_2;\Z)},$$ 
as can be seen by inspecting the following commutative diagram:
$$
\begin{CD}
H^1(Y_2) @>{\delta}>> H^1(W_1,\partial W_1)\oplus
H^1(W_2,\partial W_2)
@>>>
H^2(W,\partial W) @>>> H^2(Y_2) \\
&&&& @VVV @AAA\\
&&&&H^2(W) @>>> H^2(W_1)\oplus H^2(W_2).
\end{CD}
$$
Thus, we can construct a  
right inverse 
$$\Pi\colon 
\Z[K(W_1)]\otimes_{\Z[H^1(Y_2;\Z)]} \otimes \Z[K(W_2)]
\cong \Z\Big[\frac{K(W_1)\oplus
K(W_2)}{H^1(Y_2;\Z)}\Big]\longrightarrow
\Z[K(W)],
$$
by defining 
$$\Pi\left(e^w\right)=\left\{\begin{array}{ll}
e^v & {\text{if $i(v)=w$}} \\
0 & {\text{otherwise}}
\end{array}\right.$$
Tensoring the projection with $M$ on the left, we obtain our required
projection map $$\Pi\colon M(W_1)(W_2)\longrightarrow
M(W).$$

\begin{theorem}
\label{thm:CompositionTwist}
Fix a module $M$ for $\Z[H^1(Y_1;\Z)]$, and let $W_1$ be a cobordism
from $Y_1$ to $Y_2$ and $W_2$ be a cobordism from $Y_2$ to $Y_3$. Let
$W$ be the composite cobordism $W=W_1\cup_{Y_2} W_2$. Fix a $\SpinC$
structure $\spinc$ over $W$ and let $\spinc_i=\spinc|W_i$ for $i=1,2$;
then there are choices of representatives $\uFc_{W_1,\spinc_1}\in
{[\uFc_{W_1,\spinc_1}]}$ and
$\uFc_{W_2,\spinc_2}\in{[\uFc_{W_2,\spinc_2}]}$, for which
$$[\uFc_{W,\spinc}]= {[\Pi \circ \uFc_{W_2,\spinc_2}\circ
\uFc_{W_1,\spinc_1}]}.$$ More generally, if $h\in H^1(Y_2;\Z)$, then
for the above choices of $\uFc_{W_1,\spinc_1}$ and
$\uFc_{W_2,\spinc_2}$, we have that $$[\uFc_{W,\spinc+\delta
h}]=[\Pi_M\circ\uFc_{W_2,\spinc_2}\circ e^h\cm \uFc_{W_1,\spinc_1}],$$
where $\delta h\in H^2(W;\Z)$ is the image of $h$ under the coboundary
homomorphism for the Mayer-Vietoris sequence for the decomposition of
$W$ into $W_1$ and $W_2$.
\end{theorem}

\section{Invariants of handles}
\label{sec:ECobord}

We build the invariant of a cobordism up from its handle
decomposition.  We describe first the most interesting part -- the
part belonging to the two-handles
(Subsection~\ref{subsec:TwoHandles}). After that, in
Subsection~\ref{subsec:OneThreeHandles}, we consider the one- and
three-handles. With these definitions in place, we define the
invariant of a cobordism, and verify its topological invariance
(Theorem~\ref{thm:CobordismInvariant} in the untwisted case, and
Theorem~\ref{thm:CobordismInvariantTwist}) in
Subsection~\ref{subsec:ProveCobordismInvariance}. Finally we give a
proof of the composition law (Theorems~\ref{thm:Composition} and
\ref{thm:CompositionTwist}), which follows readily from the
definition.

\subsection{Invariants of framed links}
\label{subsec:TwoHandles}

A framed link in a three-manifold $Y$ is a collection of $n$ disjoint,
embedded circles $K_1,...,K_n\subset Y$, together with a choice of
homology classes $\ell_i\in H_1(\partial \nbd{K_i})$, with $m_i\cm
\ell_i=1$ (where $m_i$ is the meridian of $K_i$). We will often abbreviate
the framed link 
$(K_i,\ell_i)_{i=1}^n$ by $\Link$.

By attaching two-handles along the framed link $\Link$, we naturally obtain a cobordism
$W(\Link)$ from 
$Y$ to the three-manifold $Y(\Link)$
(which is obtained by surgery along the link). 

Our aim of this section is to define and study an induced map
$$\uFc_{\Link,\spinc}\colon \uHFc(Y,\spinc|_Y,M) \longrightarrow 
\uHFc(Y({\mathbb L}),\spinc|_{Y(\Link)},M(W(\Link)))$$
associated to the framed link,
where $\spinc$ is a $\SpinC$ structure on the cobordism $W(\Link)$.

\begin{defn} 
A {\em bouquet for the link $\Link$} is a one-complex embedded in $Y$
which is the union of the link $\Link=\cup_{i=1}^n K_i$ with a
collection of paths connecting $K_i$ to a fixed reference point in
$Y$.
\end{defn}

The regular neighborhood of a bouquet $B(\Link)$ is a genus $n$ handlebody $V$.
There is a subset of its
boundary which is identified with 
a disjoint union of $n$ punctured tori $F_i$. 

\begin{defn}
\label{def:HeegTrip}
A {\em Heegaard triple subordinate to the bouquet $B(\Link)$} is a
Heegaard triple $(\Sigma,\alphas,\betas,\gammas)$ with the following
properties:
\begin{list}
	{(\arabic{bean})}{\usecounter{bean}\setlength{\rightmargin}{\leftmargin}}
\item $(\Sigma, \{\alpha_1,...,\alpha_g\}, \{\beta_{n+1},...\beta_{g}\})$
describes the complement of the bouquet $B(\Link)$, 
\item $\gamma_{n+1},...\gamma_{g}$ are small
isotopic translates of the $\beta_{n+1},...,\beta_{g}$, 
\item after surgering out the $\{\beta_{n+1},...,\beta_g\}$,
the induced curves $\beta_i$ and $\gamma_i$ (for $i=1,...,n$)
lie in the punctured torus $F_i\subset \partial V$,
\item for
$i=1,...,n$, the curves $\beta_i$ represent
meridians for the $K_i$ which are disjoint from all the $\gamma_j$ for
$i\neq j$, and meet $\gamma_i$ in a single transverse intersection point
\item for $i=1,...,n$ the homology classes of the $\gamma_i$ 
correspond to the framings $\ell_i$ under the natural identification
$$H_1(\partial\nbd{\bigcup K_i})\cong H_1(\partial V).$$
\end{list}
\end{defn}

Let $X_{\alpha,\beta,\gamma}$ be the cobordism specified by the
Heegaard triple $(\Sigma,\alphas,\betas,\gammas)$.

\begin{prop}
The manifold $X_{\alpha,\beta,\gamma}$ has three boundary components,
$-Y$, $Y({\mathbb L})$, and $\#^{g-n}(S^1\times S^2)$. Furthermore,
after filling in the third boundary by the boundary connected sum
$\#^{g-n}(S^1\times B^3)$, we obtain the standard cobordism
$W(Y,{\mathbb L})$.
\end{prop}

\begin{proof} 
As a warm-up, consider the trivial case where $n=0$ -- i.e. the $\gamma_i$ are
all small translates of the $\beta_i$. In this case, it is easy to see
that $X_{\alpha,\beta,\gamma}$ is diffeomorphic to
$Y_{\alpha,\beta}\times [-1,1]$, with a regular neighborhood of the
$U_{\beta}\times \{0\}$ deleted. The boundary is, of course,
$\#^g(S^1\times S^2)$, so when we fill that back in, we obtain $Y\times I$.

Next, consider the case where $n=1$. In this case,
$$Y_{\beta, \gamma}=S^3\# \left(\#^{g-1}(S^2\times S^1)\right),$$ where the $S^3$
factor is obtained by filling in $\beta_1$ and $\gamma_1$. Now,
$$W(Y,\Link)-
\Big([0,1]\times (D^2\times \gamma_1)\cup (D^2\times D^2)\Big)\cong Y\times [-1,1].$$
Observe that $[0,1]\times (D^2\times \gamma_1)\cup (D^2\times D^2)$ is a
four-cell, attached to $Y\times I$ along a neighborhood $K$ with the
specified framing. The case for arbitrary $n$ follows in the same manner.
\end{proof}

Fix a $\SpinC$ structure $\spinc$ over $W=W(Y,{\mathbb L})$.  We
define the map $$ \uFc_{\Link,\spinc}\colon \uHFc(Y,\spinc|_Y,M)
\longrightarrow \uHFc(Y({\mathbb L}),\spinc|_{Y(\Link)},M(W)), $$ by
$$\uFc_{\Link,\spinc}(\xi)=\uFc(\xi\otimes \Theta,\spinc),$$ where
$\Theta$ is a generator for the top-dimensional homology
$\HFleq(\#^{g-n}(S^1\times S^2),\spinc_0)$. Since there are two
possible generators for this latter group, the map
$\uFc_{\Link,\spinc}$ is defined only up to sign.  Moreover, the maps
themselves are defined only up to translations by $H^1(Y;\Z)$ and $H^1(Y(\Link),\Z)$.
But more importantly, the map appears to depend on the bouquet and
the admissible triple. However, we have the following:

\begin{theorem}
\label{thm:LinkInvariant}
The map $\uFc_{\Link,\spinc}$ depends only on the three-manifold $Y$
and the framed link ${\mathbb L}$ and the $\SpinC$ structure $\spinc$
over $W(\Link)$ in the following sense.  Let
$(\Sigma_1,\alphas_1,\betas_1,\gammas_1,z_1)$ and
$(\Sigma_2,\alphas_2,\betas_2,\gammas_2,z_2)$ be a pair of Heegaard
triples which are subordinate to two bouquets for $\Link\subset Y$. Then, 
we have a commutative diagram:
$$\begin{CD}
\uHFc(\alphas_1,\betas_1,\spinc|Y,M)@>{{\uFc_{W,\spinc}}^{1}}>>
\uHFc(\alphas_1,\gammas_1,\spinc|Y(\Link),M(W)) \\
@V{\Psi^\circ_1}VV @V{\Psi^\circ_2}VV \\
\uHFc(\alphas_2,\betas_2,\spinc|Y,M)@>{{\uFc_{W,\spinc}}^{2}}>>
\uHFc(\alphas_2,\gammas_2,\spinc|Y(\Link),M(W)) \\
\end{CD}
$$
where $\Psi^\circ_1$ and $\Psi^\circ_2$ are isomorphisms 
induced by equivalences between the Heegaard diagrams 
(c.f. Theorem~\ref{thm:TwistedNaturality}).
\end{theorem}

The proof of this theorem occupies the rest of this section. 
We first show that the invariant
$\uFc_{\Link,\spinc}$ is unchanged under certain operations on the
Heegaard diagram (which leave the bouquet unchanged). 
These operations are given in the following:

\begin{lemma}
\label{lemma:LinkInvariance}
Let $Y$ be a closed, oriented three-manifold equipped
with an framed link $\Link\subset Y$ and associated bouquet $B$. Then, there 
is a Heegaard triple subordinate to $B$, and indeed
any two subordinate pointed Heegaard triples
can be connected by a sequence of the following moves:
\begin{list}
	{(\arabic{bean})}{\usecounter{bean}\setlength{\rightmargin}{\leftmargin}}
\item 
\label{item:FirstMove}
handleslides and isotopies amongst the $\{\alpha_1,...,\alpha_g\}$
\item handleslides and isotopies
amongst the $\beta_i$ for $i=n+1,...,g$, carrying along the $\gamma_i$
for $i=n+1,...,g$, as well
\item isotopies and possible handleslides of some $\beta_i$ for $i\in 1,...,n$ across
$\beta_j$ for $j\in n+1,...,g$
\item 
\label{item:PenultimateMove}
isotopies and possible handleslides of the $\gamma_i$ for $i\in 1,...,n$ across
the $\gamma_j$ for $j\in n+1,...,g$
\item stabilizations (introducing $\alpha_{g+1}$, $\beta_{g+1}$, and
$\gamma_{g+1}$).
\end{list}
As usual, isotopies and handleslides here take place in the complement
of the basepoint $z$.
\end{lemma}

\begin{proof}
Fix a Morse function on $B(\Link)$ with one critical point of index
three and $n$ index two critical points. Let $\{\beta_1,...,\beta_n\}$
be the attaching circles for these index two critical points. 

We complete the above to a 
triple extending this Morse function over 
all of $Y-B(\Link)$ (introducing only one 
index zero,  and no new index three critical points). This gives rise to
$\Sigma$, $\{\alpha_1,...,\alpha_g\}$ and
$\{\beta_{n+1},...,\beta_g\}$ in the usual manner. Let $z$ then be any
basepoint in the complement of these curves. Curves
$\{\gamma_{n+1},...,\gamma_g\}$ can be found then by taking small
translates of the corresponding $\beta_i$.

Now, if we have two such Morse functions (giving rise to Heegaard
triples subordinate to a fixed bouquet), we can connect the Morse
functions through a generic one-parameter family. This allows us to
connect the two handle decompositions for $Y-B(\Link)$ through a
sequence of handle-slides, pair creations and annihilations. 

In this process, we possibly introduce new index zero and index three
critical points which cancel new attaching circles $\alpha_i$ and
$\beta_j$. It is straightforward to see that no matter how these
cancellations occur, the two handle decompositions differ by only
isotopies and handleslides amongst the $\{\alpha_1,...,\alpha_g\}$ and 
$\{\beta_{n+1},...,\beta_g\}$, and also 
stabilizations (c.f. Proposition~\ref{HolDiskOne:prop:HeegaardMoves}
of~\cite{HolDisk}). 

Observe that if we surger $\{\beta_{n+1},...,\beta_g\}$ out of
$\Sigma$, the remaining two-manifold is identified with $\partial
B(\Link)$; in particular, we have candidates for the curves
$\{\gamma_1,...,\gamma_n\}$ in $\partial B(\Link)$, which we can
perturb slightly so that they induce curves in $\Sigma$.

Next, we consider the dependence of the construction on the $\beta_i$
and $\gamma_i$ for $i=1,...,n$.  The  link specifies the
homology classes of the $\beta_i$ in $F_i$, while the framing
specifies the homology classes of the $\gamma_i$. Different choices
$\beta_i$ and $\beta_i'$ can be connected by an isotopy in $F_i$. It
follows that in $\Sigma$ they can be connected by a sequence of
isotopies and handleslides across the
$\{\beta_{n+1},...,\beta_g\}$. The $\gamma_i$ for $i=1,...,n$ follow
in the same manner.

Since $\Sigma-\beta_1-...-\beta_g-\gamma_1-...-\gamma_n$ is connected,
we can realize moving the basepoint as a sequence of handleslides
amongst the $\alpha_i$ (c.f. Section~\ref{HolDiskOne:sec:HandleSlides}
of~\cite{HolDisk}).
\end{proof}

We have the following special case of Theorem~\ref{thm:LinkInvariant}:

\begin{prop}
\label{prop:IndepOfHeegaard}
Fix a framed link $\Link\subset Y$, a $\SpinC$ structure $\spinc\in
W(\Link)$, and a bouquet $B$ for $\Link$.  Then, the map
$\uFc_{\Link,\spinc}$ is independent of the underlying Heegaard triple
subordinate to the bouquet, in the sense that if 
$(\Sigma_1,\alphas_1,\betas_1,\gammas_1,z_1)$ and
$(\Sigma_2,\alphas_2,\betas_2,\gammas_2,z_2)$ are a pair of Heegaard
triples which are subordinate to a fixed  bouquet $B(\Link)$ 
for $\Link\subset Y$. Then, 
we have a commutative diagram:
$$\begin{CD}
\uHFc(\alphas_1,\betas_1,\spinc|Y,M)@>{{\uFc_{W,\spinc}}^{1}}>>
\uHFc(\alphas_1,\gammas_1,\spinc|Y(\Link),M(W)) \\
@V{\Psi^\circ_1}VV @V{\Psi^\circ_2}VV \\
\uHFc(\alphas_2,\betas_2,\spinc|Y,M)@>{{\uFc_{W,\spinc}}^{2}}>>
\uHFc(\alphas_2,\gammas_2,\spinc|Y(\Link),M(W)) \\
\end{CD}
$$
where $\Psi^\circ_1$ and $\Psi^\circ_2$ are isomorphisms 
induced by equivalences between the Heegaard diagrams 
\end{prop}

This proposition is divided into steps, the first of which is
a stabilization invariance for holomorphic triples, which in turn
follows from the gluing result for holomorphic triangles,
Theorem~\ref{thm:GlueTriangles}.

\begin{lemma}
\label{lemma:Stabilize}
The map $\uFc_{\Link,\spinc}$ commutes with the isomorphisms
induced by stabilizations
(introducing $\alpha_{g+1}$, $\beta_{g+1}$, and $\gamma_{g+1}$). 
\end{lemma}

\begin{proof}
Start with a Heegaard triple $(\Sigma,\alphas,\betas,\gammas,z)$.
The stabilization involves forming the connected sum of 
$\Sigma$  with a
genus one surface $E$ containing: a pair of curves $\alpha_{g+1}$ and
$\beta_{g+1}$ (which meet in a single, transverse intersection point),
and a new curve $\gamma_{g+1}$ which is a small isotopic copy of
$\beta_{g+1}$ meeting it in two transverse intersection points (and
meeting $\alpha_g$ in a single transverse intersection point). 
Then, the new Heegaard triple 
$$(\Sigma,\{\alpha_1,...,\alpha_{g+1}\},\{\beta_1,...,\beta_{g+1}\},
\{\gamma_1,...\gamma_{g+1}\},z)$$
represents a stabilization of the original Heegaard triple.

Let $x_0=\alpha_{g+1}\cap\beta_{g+1}$,
$w_0=\alpha_{g+1}\cap\gamma_{g+1}$, and let $y_0$ be the intersection
point of $\beta_{g+1}$ with $\gamma_{g+1}$ with higher relative
degree. It is easy to see that there is a unique homotopy class of
triangle $\psi_0\in\pi_2(x_0,y_0,w_0)$ with nowhere negative
coefficients, and moreover $\Mas(\psi_0)=0$, and $\psi_0$ has a
unique, smooth holomorphic representative (if we take a constant
complex structure on $E$), see Figure~\ref{fig:Stabilize}.  
Clearly, if the $\Theta\in
\HFleq(\Tb,\Tc,\spinc_0)$ represents a top-dimensional generator,
where $\Tb=\beta_1\times...\times\beta_g$ and
$\Tc=\gamma_1\times...\times\gamma_g$, then
$\Theta'=\Theta\times\{y_0\}\in\HFp(\Tb\times\beta_{g+1},\Tc\times\gamma_{g+1},\spinc_0)$.

\begin{figure}
\mbox{\vbox{\epsfbox{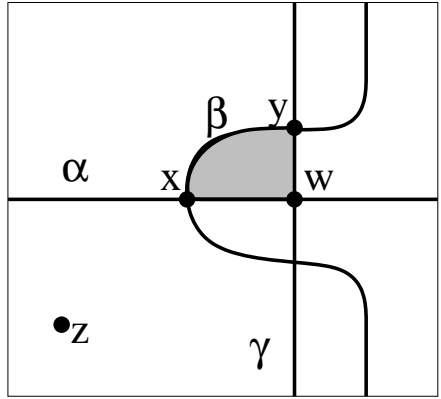}}}
\caption{\label{fig:Stabilize}
{\bf{Stabilization.}}  We have pictured here the torus $E$ used in
stabilizing. (We have dropped all the subscripts from the picture.) 
The connected sum point is labelled $z$. The domain of $\psi_0$ is
lightly shaded in the picture.}
\end{figure}

Consider the square
$$
\begin{CD}
\CFp(\Ta,\Tb,\spinc|_Y) @>{F^+_1}>> \CFp(\Ta,\Tc,\spinc|_{Y(\Link)}) \\
@V{\cong}VV @VV{\cong}V \\
\CFp(\Ta\times\alpha_{g+1},\Tb\times\beta_{g+1})
@>{F^+_2}>>
\CFp(\Ta\times\alpha_{g+1},\Tb\times\gamma_{g+1}),
\end{CD}
$$ where the vertical isomorphisms are induced by the maps
$\x\in\Sym^g(\Sigma-\nbd{z})\mapsto \x\times x_0\in
\Sym^{g+1}(\Sigma\# E)$ and $\y\mapsto \y\times y_0$. For 
all sufficiently long connected sum neck, these induce are
isomorphisms of chain complexes, according to
Theorem~\ref{HolDiskOne:thm:StabilizeHFb} of~\cite{HolDisk}.
Theorem~\ref{thm:GlueTriangles} gives a map
$$\pi_2(\x,\Theta,\w)\longrightarrow \pi_2(\x\times
x_0,\Theta',\w\times w_0)$$ and an identification of formally
zero-dimensional moduli spaces $\ModFlow(\psi)\cong\ModFlow(\psi')$,
which shows that the square exhibited above commutes.
\end{proof}

\vskip.2cm
\noindent{\bf{Proof of Proposition~\ref{prop:IndepOfHeegaard}.}}
First, we organize Lemma~\ref{lemma:LinkInvariance} in the same manner
as Lemma~\ref{lemma:ConnectOrdered}: if
$(\Sigma_1,\alphas_1,\betas_1,\gammas_1,z_1)$ and
$(\Sigma_2,\alphas_2,\betas_2,\gammas_2,z_2)$ are a pair of Heegaard
triples, both of which are subordinate to a fixed bouquet $B(\Link)$
for a framed link $\Link\subset \Knot$, then we can go from one to the
other by first a sequence of stabilizations, then a strong equivalence
-- i.e. using only
Moves~(\ref{item:FirstMove})-(\ref{item:PenultimateMove}) of
Lemma~\ref{lemma:LinkInvariance} (in particular leaving the underlying
two-manifold $\Sigma$ unchanged) -- and then a sequence of
destabilizations.

Having shown (in Lemma~\ref{lemma:Stabilize}) that $\uFc_{\Link,\spinc}$ 
commutes
with the isomorphism induced by stabilization, it suffices to show
that $\uFc_{\Link,\spinc}$ commutes with the isomorphism induced by a strong
equivalence. 

Suppose that $(\Sigma,\alphas_1,\betas_1,\gammas_1,z)$ and
$(\Sigma,\alphas_2,\betas_2,\gammas_2)$ are strongly equivalent. Then
it is easy to see that $\alphas_2$, $\betas_2$, and $\gammas_2$ are
obtained by handleslides and isotopies amongst the $\alphas_1$,
$\betas_1$, and $\gammas_1$ respectively. Following the scheme from
Subsection~\ref{subsec:ProofNaturality}, we choose isotopic copies
$\alphas_1'$, $\betas_1'$ and $\gammas_1'$ of the $\alphas_1$,
$\betas_1$, and $\gammas_1$, so that each of
$(\Sigma,\alphas_1,\alphas_1',z)$, $(\Sigma,\betas_1,\betas_1',z)$ and
$(\Sigma,\gammas,\gammas_1',z)$ is admissible. Now the fact that
$F_{\Knot,\spinc}$ commutes with the isomorphism of the strong
equivalence follows from the commutativity of the following two squares.

First, we have
$$
\begin{CD}
\HF^\circ(\alphas_1,\betas_1)@>{\Gamma_{\alpha_1',\alpha_1;\beta_1,\beta_1'}}>>
\HF^\circ(\alphas_1',\betas_1') \\
@V{F^\circ(\cdot{\otimes\Theta_{\beta_1,\gamma_1}},\spinc)}VV 
@VV{F^\circ(\cdot\otimes \Theta_{\beta_1',\gamma_1'},\spinc)}V \\
\HF^\circ(\alphas_1,\gamma_1)
@>{\Gamma_{\alpha_1',\alpha_1;\gamma_1,\gamma_1'}}>>
\HF^\circ(\alphas_1',\gamma_1').
\end{CD}
$$

This square commutes up to sign since 
isotopy invariance of the triangle construction gives that
\begin{eqnarray*}
\Gamma_{\alpha_1',\alpha_1;\gamma_1,\gamma_1'}
\circ F^\circ_{\alpha_1,\beta_1,\gamma_1}(\xi\otimes\Theta_{\beta_1,\gamma_1},\spinc)
&=& F^\circ_{\alpha_1',\beta_1,\gamma_1'}(\Gamma_{\alpha_1',\alpha_1;\beta_1}(\xi)\otimes
\Gamma_{\beta_1;\gamma_1,\gamma_1'}(\Theta_{\beta_1,\gamma_1}),\spinc) \\
&=& \pm
F^\circ_{\alpha_1',\beta_1,\gamma_1'}(\Gamma_{\alpha_1',\alpha_1;\beta_1}(\xi)\otimes
\Theta_{\beta_1,\gamma_1'}) \\
&=& 
\pm F^\circ_{\alpha_1',\beta_1,\gamma_1'}(\Gamma_{\alpha_1',\alpha_1;\beta_1,\beta_1'}(\xi)\otimes
\Gamma_{\beta_1',\beta_1;\gamma_1'}(\Theta_{\beta_1,\gamma_1'})) \\
&=& 
\pm F^\circ_{\alpha_1',\beta_1,\gamma_1'}(\Gamma_{\alpha_1',\alpha_1;\beta_1,\beta_1'}(\xi)\otimes
\Theta_{\beta_1',\gamma_1'})),
\end{eqnarray*}
in view of the fact that
$\Gamma_{\beta_1;\gamma_1,\gamma_1'}(\Theta_{\beta_1,\gamma_1'})=\pm
\Theta_{\beta_1,\gamma_1'}$, and 
$\Gamma_{\beta_1',\beta_1;\gamma_1'}(\Theta_{\beta_1,\gamma_1'})=
\pm \Theta_{\beta_1',\gamma_1'}$ (since both sides of both equations are
generators for a group  -- the top-dimensional homology of 
$\HFleq(\#^{g-1}(S^1\times S^2),\spinc_0)$ -- which is isomorphic to $\Z$).

The second commutative square is 
$$
\begin{CD}
\HF^\circ(\alphas_1',\betas_1')@>{\Theta_{\alpha_2,\alpha_1'}\otimes~\cdot~\otimes
\Theta_{\beta_1',\beta_2}}>> 
\HF^\circ(\alphas_2,\betas_2) \\
@V{\otimes\Theta_{\beta_1',\gamma_1'}}VV 
@VV{\otimes\Theta_{\beta_2,\gamma_2}}V \\
\HF^\circ(\alphas_1',\gammas_1')
@>{\Theta_{\alpha_2,\alpha_1'}\otimes
~\cdot~\otimes\Theta_{\gamma_1',\gamma_2}}>> 
\HF^\circ(\alphas_2,\gammas_2) \\
\end{CD}
$$
which now commutes up to sign, owing to the associativity 
of the triangle construction (Theorem~\ref{thm:Associativity}), and now the fact that
$$F^{\leq 0}(\Theta_{\beta_1',\gamma_1'}\otimes\Theta_{\gamma_1',\gamma_2},\spinc_0)
=\pm \Theta_{\beta_1',\gamma_2}
= \pm F^{\leq 0}(\Theta_{\beta_1',\beta_2}\otimes\Theta_{\beta_2,\gamma_2},\spinc_0).$$
\qed
\vskip.2cm

\begin{lemma}
\label{lemma:IndepOfBouquet}
The map $\uFc_{\Link,\spinc}$ is independent of the bouquet.
\end{lemma}

\begin{proof}
Suppose that $B$ and $B'$ are a pair of bouquets which differ in the
choice of one path $\sigma_1$ and $\sigma_1'$. We construct two Heegaard
triples 
\begin{eqnarray*}
(\Sigma,\alphas,\betas,\gammas,z) &{\text{and}}&
(\Sigma,\alphas,\betaprs,\gammaprs,z)
\end{eqnarray*}
in such a manner that the sets
$\betaprs$ are obtained from the $\betas$ by a sequence of
handleslides, and the $\gammaprs$ are obtained from the $\gammas$ by a
sequence of handleslides.

To this end, consider the regular neighborhood $M$ of $B\cup B'$, This
is a handlebody of genus $n+1$.  On the complement $Y-M$, we have a
Morse function with one index zero critical point, $g$ index one, and
$g-n-1$ index two critical points. Let $\{\alpha_1,...,\alpha_g\}$ and
$\{\beta_{n+2},...,\beta_g\}$ be the corresponding attaching
circles. Let $\{\beta_1,...,\beta_n\}$ be attaching circles in
$\partial M$ for the index two critical points corresponding to
reattaching the components $K_1,...,K_n$ of the link. There is one
additional index two critical point on $M$. Let $\beta_{n+1}$ be the
attaching circle dual to $\sigma_1$: this is the attaching circle for
a Morse function on $M-\nbd{B}$ with a single index two critical
point. Similarly, let $\beta_{n+1}'$ be the attaching circle dual to
$\sigma'$ (while all other $\beta_i'=\beta_i$). We can arrange for
$\beta_{n+1}'$ to be disjoint from $\beta_{n+1}$. Let $\gamma_i$ for
$i=1,...,n$ correspond to the given framing of the link, and
$\gamma_i$ be small isotopic translates of the $\beta_i$ for
$i=n+1,...,g$, while $\gamma_i'$ is a small isotopic translate of
$\beta_i'$ (while all other $\gamma_i'=\gamma_i$), then,
$(\Sigma,\alphas,\betas,\gammas,z)$ is subordinate to $B(\Link)$,
while $(\Sigma,\alphas,\betaprs,\gammaprs,z)$ is subordinate to
$B(\Link')$.

Now, if we surger
$\{\beta_1,...,\beta_{n-1},\beta_{n+1},...,\beta_{g}\}$ out of
$\Sigma$, we obtain a surface of genus one, with two disjoint,
embedded, homologically non-trivial curves induced by $\beta_{n}$ and
$\beta_n'$. These curves must then be isotopic in the torus: thus,
$\beta_{n}'$ can be obtained by handlesliding $\beta_n$ over some collection 
of the
$\{\beta_1,...,\beta_{n-1},\beta_{n+1},...,\beta_{g}\}$. The analogous remark
applies to obtaining $\gamma_n'$ from $\gamma_n$.
Thus, $\betas'$ are obtained from $\betas$ by handleslides and isotopies, and 
$\gammas'$ are obtained from $\gammas$ by handleslides and isotopies.

The result then follows from the usual commutative diagram
$$
\begin{CD}
\uHFc(\alphas,\betas,\spinc|_Y,M)
@>{\uFc(\cdot\otimes \Theta_{\beta,\gamma},\spinc)}>> 
\uHFc(\alphas,\gammas,\spinc|_{Y(\Link)},M(W)) \\
@V{\uFc(\cdot\otimes \Theta_{\beta,\beta'},\spinc)}VV 
@VV{\uFc(\cdot\otimes \Theta_{\gamma,\gamma'},\spinc)}V \\
\uHFc(\alphas,\betas',\spinc|_Y,M)
@>{\uFc(\cdot\otimes \Theta_{\beta',\gamma'},\spinc)}>> 
\uHFc(\betas',\gammas',\spinc|_{Y(\Link)},M(W)), 
\end{CD}
$$ (where the vertical maps are isomorphisms induced by strong
equivalence, and the two horizontal maps are the two candidates for
$\uFc_{\Link,\spinc}$ belonging to the two bouquets), and 
use of associativity, according to which we have the
following equation up to sign 
$$\Fleq{}(\Theta_{\beta,\gamma}\otimes
\Theta_{\gamma,\gamma'},\spinc_0)=\pm \Theta_{\beta,\gamma'}=
\pm\Fleq{}(\Theta_{\beta,\beta'}\otimes\Theta_{\beta',\gamma'},\spinc_0).$$
\end{proof}

\vskip.2cm
\noindent{\bf{Proof of Theorem~\ref{thm:LinkInvariant}}.}
This is now an immediate consequence of Proposition~\ref{prop:IndepOfHeegaard},
together with Lemma~\ref{lemma:IndepOfBouquet}.
\qed
\vskip.2cm

\subsection{Compositions of link invariants}

The invariant $\uFc_{\Link,\spinc}$ satisfies the following analogue of 
Theorem~\ref{thm:Composition}.

If we partition the link $\Link=\Link_1\cup\Link_2$, this gives 
a decomposition of 
the cobordism $W(\Link)$ as a union $$ W(\Link)=W_1\cup_{Y(\Link_1)}
W_2,$$ where $W_1=W(Y,\Link_1)$, and $W_2$ is the cobordism from
$Y(\Link_1)$ to $Y(\Link)$ associated to $\Link_2$, thought of as a
framed link in $Y(\Link_1)$. Recall that we have a  projection map
$$\Pi\colon \uHFc(Y(\Link),M(W_1)(W_2))\longrightarrow \uHFc(Y(\Link),M(W)).$$

\begin{prop}
\label{prop:CompositeLink}
There are choices $\uFc_{W_1,\spinc_1}$ and $\uFc_{W_2,\spinc_2}$
from the equivalence class of
maps  $\uFc_{Y,\Link_1,\spinc|W_1}$ and $\uFc_{Y(\Link_1),\Link_2}$
respectively, with the property that the projection of the
composite map $\Pi\circ \uFc_{W_2}\circ \uFc_{W_1}$ 
represents $\uFc_{Y,\Link}$.
\end{prop}

\begin{proof}
Let $\Link_1$ and $\Link_2$ have $m$ and $n$ components, respectively.
Let $\{\beta_1,...,\beta_{n+m}\}$ correspond to the link
$\Link_1\cup\Link_2$. Let $\delta_i$ be small isotopic translates of
the $\gamma_i$ for $i=1,...,m$ correspond to the framings of
$\Link_1$, and let $\delta_i$  be small isotopic translates of
$\beta_i$ for $i=m+1,...,m+n$. Clearly, 
$$\uFc_{Y,\Link_1}=
\uFc(\cdot\otimes\Theta_{\beta,\delta},\spinc,M)
\colon \uHFc(\alphas,\betas,\spinc|_Y,M) \longrightarrow
\uHFc(\alphas,\deltas,\spinc|_{Y(\Link_1)},M(W_1)),$$
while 
$$\uFc_{Y(\Link_1);\Link_2}=
\uFc(\cdot\otimes\Theta_{\delta,\gamma},\spinc,M(W_1))
\colon \uHFc(\alphas,\deltas,\spinc|_{Y(\Link_1)},M(W_1)(W_2)) 
\longrightarrow
\uHFc(\alphas,\gammas,\spinc|_{Y(\Link_1\cup\Link_2)}),$$

Next, we claim that 
$$\Fleq{}(\Theta_{\beta,\delta}\otimes
\Theta_{\delta,\gamma},\spinc_0)=\pm\Theta_{\beta,\gamma}.$$ This follows from
Theorem~\ref{thm:GlueTriangles}, and two model calculations in the
genus one surface $E$. In one case, we have three curves $\beta$,
$\gamma$, and $\delta$ which are all three Hamiltonian isotopic
translates of one another, as in Figure~\ref{fig:OneHandle} below
(with different labellings); in the other case, $\beta$ and $\delta$
are small isotopic translates of one another, and $\gamma$ intersects
$\beta$ (and also $\delta$) in a single, transverse intersection point
(we met this configuration already in the proof of stabilization
invariance, see Figure~\ref{fig:Stabilize}, only with different
notation). In each case, there is a unique homotopy class with nowhere
negative coefficients connecting the relevant intersection points, and
the homotopy class supports a single smooth solution.

The theorem now follows from associativity in the twisted case. Note
that in the present application of associativity, we consider the
pointed Heegaard quadruple
$(\Sigma,\alphas,\betas,\deltas,\gammas,z)$, and as such we must
verify that 
$$
\delta H^1(Y_{\alpha,\delta})|Y_{\beta,\gamma}=0=
\delta H^1(Y_{\beta,\gamma})|Y_{\alpha,\delta}.
$$
Both follow from the fact that the map on $H_2(Y_{\beta,\gamma})\longrightarrow H_2(X_{\alpha,\beta,\delta,\gamma})$ is trivial
(to see this, observe that the $(\betas,\gammas)$-periodic domains,
which evidently involve relations amongst the $\beta_i$ and $\gamma_j$
for $i,j\in \{m+1,...,m+n\}$, can be expressed as sums of
$(\betas,\gammas)$- and $(\betas,\deltas)$-periodic domains). 
\end{proof}

\subsection{One- and three-handles}
\label{subsec:OneThreeHandles}

Let $U$ be the cobordism obtained by adding a single one-handle to a
connected three-manifold $Y$. This is a cobordism between $Y$ and
$Y'=Y\#(S^1\times S^2)$. Clearly, $\SpinC(U)\cong \SpinC(Y)$;
moreover, the $\SpinC$ structures over $Y'$ which extends over $U$ are
those which have trivial first Chern class on the $S^1\times S^2$
factor.

We define the invariant for the one-handle addition: as follows.  Fix
a standard Heegaard diagram $(E,\alpha,\beta,z_0)$ for $S^2\times S^1$
(in the sense of Definition~\ref{def:StandardHeegaard}) so that
$\alpha$ and $\beta$ meet in a pair of intersection points; and let
$\theta$ be the maximal one.  Given a Heeegaard decomposition
$(\Sigma,\alphas,\betas,z)$ for $Y$, let
$(\Sigma',\alphas',\betas',z')=(\Sigma,\alphas,\betas,z)\#(E,\alpha,\beta,z_0)$
be an associated Heegaard decomposition for $Y\#(S^2\times
S^1)$. Then, we define $$\ugc{U,\spinc}\colon
\uCFinf(\alphas,\betas,\spinc,M)
\longrightarrow \uCFinf(\alphas',\betas',
\spinc\#\spinc_0,M)$$ be the map induced from
$$\uginf{U,\spinc}(m\otimes[\x,i])=m\otimes [\x\times\{\theta\},i].$$
(Observe that $M(U)\cong M$.)
When the complex structure on $\Sigma\# E$ is sufficiently
stretched out, this is a chain map
(c.f. Proposition~\ref{HolDiskTwo:intro:ConnectedSumSphereProd}
of~\cite{HolDiskTwo}), and we define the map associated to the one-handle
$$\uGc{U,\spinc}\colon \uHFc(Y,\spinc|Y,M)\longrightarrow
\uHFc(Y',\spinc|Y',M)$$
be the induced map on homology.

\begin{theorem}
\label{thm:OneHandles}
The maps $\uGc{U,\spinc}$ depends only on the three-manifold $Y$
and the $\SpinC$ structure $\spinc$ in the following sense. If
$(\Sigma_1,\alphas_1,\betas_1,z_1)$ and
$(\Sigma_2,\alphas_2,\betas_2,z_2)$ are equivalent Heeegaard diagrams,
then the corresponding diagrams
$(\Sigma_1',\alphas_1',\betas_1',z_1')$ and
$(\Sigma_2',\alphas_2',\betas_2',z_2')$ are equivalent, and indeed we
have a commutative diagram 
$$ \begin{CD}
\HFc(\alphas_1,\betas_1,\spinct) @>{\uGc{U,\spinc}}>>
\HFc(\alphas_1',\betas_1',\spinct\#\spinc_0) \\
@V{\Psi^\circ_Y}VV	@VV{\Psi^\circ_{Y\# (S^2\times S^1)}}V \\
\HFc(\alphas_2,\betas_2,\spinct)
@>{\uGc{U,\spinc}}>>
\HFc(\alphas_2',\betas_2',\spinct\#\spinc_0),
\end{CD}
$$
where the vertical maps are the isomorphisms induced from the
equivalences of the Heegaard diagrams
(c.f. Theorem~\ref{thm:Naturality}).
\end{theorem}

\begin{proof}
The equivalence of $(\Sigma,\alphas_1,\betas_1,z_1)$ and
$(\Sigma,\alphas_2,\betas_2,z_2)$ induces the equivalence between the
corresponding diagrams $(\Sigma_1,\alphas_1',\betas_1',z_1')$ and
$(\Sigma_2,\alphas_2',\betas_2',z_2')$ in view of the fact that an
isotopy of an attaching circle in a Heegaard diagram for $Y$
corresponds to a pair of handleslides across $\alpha_{g+1}$ or
$\beta_{g+1}$ supported in $E$ for the induced Heegaard diagram of
$Y\#(S^2\times S^1)$.

The fact that $G^\infty_{U,\spinc}$ commutes with the isomorphism
induced by equivalence, follows from a gluing analogous to
Lemma~\ref{lemma:Stabilize}, only now the corresponding picture in the
torus $E$ is different; see Figure~\ref{fig:OneHandle}. Here, $\alpha$,
$\beta$, and $\gamma$ are all isotopic translates of the same curve
(the curve $\alpha_{g+1}$ in $E$).

\begin{figure}
\mbox{\vbox{\epsfbox{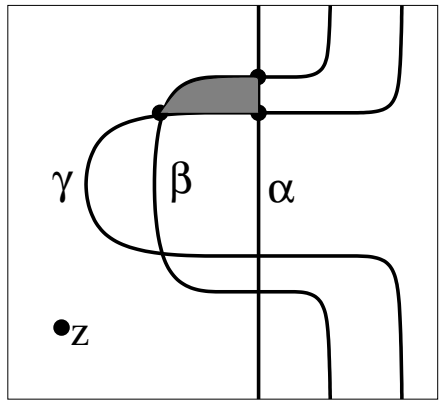}}}
\caption{\label{fig:OneHandle}
{\bf{One handles.}}  Three curves in the torus, obtained as small
isotopic translates of one another. The domain of holomorphic triangle
connecting the top-dimensional intersection points is lightly shaded.}
\end{figure}

\end{proof}

Dually, if $V$ is the cobordism obtained by adding a single
three-handle along a non-separating two-sphere $S\subset Y'$, then 
$Y'$ can be written as $Y\#(S^2\times S^1)$ and $V$
is a cobordism from $Y'=Y \# (S^2\times S^1)$ to $Y$.
There is a special kind of compatible Heegaard diagram induced by the
embedded sphere.

\begin{lemma}
\label{lemma:SphereHeegaard}
Let $S\subset Y'$ be a non-separating embedded sphere in a
three-manifold, then there is an induced split Heegaard diagram
$(\Sigma',\alphas',\betas',z')$ for $Y'$ of the form
$$(\Sigma',\alphas',\betas',z')=(\Sigma,\alphas,\betas,z)\#(E,\alpha,\beta,z_0),$$
where $(E,\alpha,\beta,z_0)$ is a standard Heegaard diagram of
$S^2\times S^1$ (and where the sphere $S$ is represented in this
factor). Moreover, if we have two such split diagrams which 
are equivalent,
\begin{eqnarray*}
(\Sigma_1,\alphas_1,\betas_1,z_1)\#(E,\alpha,\beta,z_0)&{\text{and}}&
(\Sigma_2,\alphas_2,\betas_2,z_2)\#(E,\alpha,\beta,z_0),
\end{eqnarray*} then $(\Sigma_1,\alphas_1,\betas_1,z_1)$
and $(\Sigma_2,\alphas_2,\betas_2,z_2)$ are equivalent Heegaard diagrams.
\end{lemma}

\begin{proof}
We build the Heegaard diagram starting from a handle decomposition of
the neighborhood of the two-sphere, as in
Lemma~\ref{HolDiskTwo:lemma:SpecialHeegaard} of~\cite{HolDiskTwo}
(where we were interested in the case of embedded surfaces of genus
$g>0$, though the construction works the same way when
$g=0$). Specifically, we start with Morse function on the neighborhood
of the two-sphere $S$, and consider extensions to the complement
$Y-S$.  In this manner, we obtain a Heegaard diagram
$(\Sigma',\alphas',\betas',z')$, with a periodic domain representing
$S$ which is bounded by curves $\alpha_{g+1}$ and $\beta_{g+1}$ which
are small translates of one another. Indeed, in the present case, we
arrange that $\beta_{g+1}$ is a small Hamiltonian translate of
$\alpha_{g+1}$. Now, any two Heegaard
diagrams which arise in this manner are equivalent, through an
equivalence which leaves $\alpha_{g+1}$ and $\beta_{g+1}$ unchanged.

Pick any curve $\delta$ dual to $\alpha_{g+1}$ (and
$\beta_{g+1}$). It is easy to see that after handleslides across the
$\alpha_{g+1}$, we can arrange for all the remaining
$\{\alpha_1,...,\alpha_g\}$ to be disjoint from $\gamma$. Similarly,
handlesliding across $\beta_{g+1}$, we arrange for the
$\{\beta_1,...,\beta_{g}\}$ to be disjoint from $\delta$. The curves
$\delta$ and $\alpha_{g+1}$ and $\beta_{g+1}$ lie in a torus summand
of $\Sigma'$, giving us the required splitting. 

To verify uniqueness of the splitting, suppose we have two different
such dual curves $\delta$ and $\delta'$. 
In the first case, we destabilize along the
neighborhoods $\alpha_{g+1}\cup \delta$ (i.e. surgering out a regular
neighborhood), and in the second we destabilize along
$\alpha_{g+1}\cup\delta'$. Both two-manifolds are identified
with the two-manifold obtained by surgering out $\alpha_{g+1}$. To see
that the two induced Heegaard diagrams are equivalent, we make the
following observation. Suppose that $\alpha_i$ is a curve which meets
$\delta$, but misses $\delta'$. By handlesliding repeatedly across
$\alpha_{g+1}$, we obtain a new curve $\alpha_i'$ which meets meets
$\delta'$, but is disjoint from $\delta$. However, if we then surger
out $\alpha_{g+1}$, the induced curves (images of $\alpha_i$ and
$\alpha_i'$ are isotopic. (See Figure~\ref{fig:ThreeHandles}.)

\begin{figure}
\mbox{\vbox{\epsfbox{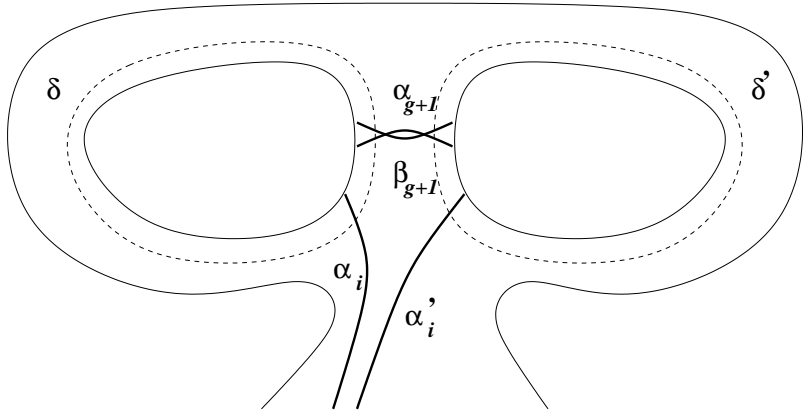}}}
\caption{\label{fig:ThreeHandles}
{\bf{Two-spheres splitting Heegaard diagrams.}}  This is an illustration of the uniqueness claim in
Lemma~\ref{lemma:SphereHeegaard}.  The curves $\alpha_{g+1}$ and
$\beta_{g+1}$ bound a periodic domain representing the embedded
two-sphere, and the curves $\delta$ and $\delta'$ give two candidates
for the splitting of $E$. Handlesliding $\alpha_i$ across $\alpha_{g+1}$, we get
a curve $\alpha_i'$ which is isotopic to $\alpha_i$, after we surger out
$\alpha_{g+1}$.}
\end{figure}

\end{proof}

We now take a split Heegaard diagram 
$$(\Sigma',\alphas',\betas',z')=
(\Sigma,\alphas,\betas,z)\#(E,\alpha,\beta,z_0)$$ 
as in the above lemma.
Then, we define $$\ueinf{V,\spinc}\colon
\CFinf(\alphas',\betas',
\spinc|Y') \longrightarrow \CFinf(\alphas,\betas,\spinc|Y)$$
to be the map
$$\ueinf{V,\spinc}[\x\times \{y\},i]
=\left\{\begin{array}{ll}
[\x,i] & {\text{if $y$ is the minimal intersection point
for the $S^2\times S^1$ factor}} \\
0 & {\text{otherwise}}
\end{array}
\right.$$
Again, if the complex structure on $\Sigma\#E$ is sufficiently
stretched, out, $\ueinf{V,\spinc}$ gives a chain map, inducing the
required map 
on homology
$$\uEc{V,\spinc}\colon
\uHFc(Y',\spinc|Y',M)\longrightarrow \uHFc(Y,\spinc|Y,M).$$

The following analogue of Theorem~\ref{thm:OneHandles} holds for
three-handles.

\begin{theorem}
\label{thm:ThreeHandles}
Fix a non-separating embedded sphere $S\subset Y$, and let
$$(\Sigma_1',\alphas_1',\betas_1',z_1')
=(\Sigma_1,\alpha_1,\beta_1,z_1)\#(E,\alpha,\beta,z_0)$$
and 
$$(\Sigma_2',\alphas_2',\betas_2',z_2')
=(\Sigma_2,\alpha_2,\beta_2,z_2)\#(E,\alpha,\beta,z_0)$$
be a pair of equivalent split Heegaard diagrams, then 
the following square commutes:
$$ \begin{CD}
\HFc(\alphas_1',\betas_1',\spinct\#\spinc_0) @>{\uEc{V,\spinc}}>>
\HFc(\alphas_1,\betas_1,\spinct) \\
@V{\Psi^\circ_{Y'}}VV	@VV{\Psi^\circ_{Y}}V \\
\HFc(\alphas_2',\betas_2',\spinct\#\spinc_0)
@>{\uEc{V,\spinc}}>>
\HFc(\alphas_2,\betas_2,\spinct),
\end{CD}
$$ where the vertical maps are the isomorphisms induced by the
equivalences.
\end{theorem}

\begin{proof}
This follows exactly as in the case for one-handles,
Theorem~\ref{thm:OneHandles}: i.e. we apply the gluing of
Theorem~\ref{thm:GlueTriangles} for the appropriate triangle in the torus
with three curves illustrated in Figure~\ref{fig:OneHandle}.
\end{proof}

\subsection{Invariance of the maps associated to cobordisms}
\label{subsec:ProveCobordismInvariance}

We prove
Theorems~\ref{thm:CobordismInvariantTwist} and
\ref{thm:CobordismInvariant}.

Let $W$ be a cobordism from $Y_1$ to $Y_2$. We decompose $W=W_1\#_{Y_1'} W_2
\#_{Y_2'} W_3$, where $W_1$ is a collection of one-handles, $W_2$ is a a
collection of two-handles, and $W_3$ is a collection of three-handles.
Concretely, $W_2$ can be represented as a framed link $\Link\subset
Y_1'=Y_1\#\left(\#^{\ell}(S^2\times S^1)\right)$. We then define
$$\uFc_{W,\spinc}=\uEc{W_3,\spinc}\circ \uFc_{\Link,\spinc}\circ
\uGc{W_1,\spinc},$$
where $\uEc{W_3,\spinc}$ and $\uGc{W_1,\spinc}$ are defined to be composites 
of the maps $\uEc{}$ and $\uGc{}$ induced by the the various one- and three-handles.

To see that this depends on the underlying four-manifold only (not on
the particular handle decomposition), we proceed in several steps.

\begin{lemma}
The maps $\uGc{W_1,\spinc}$ are invariant under 
the ordering of the one-handles, and handleslides amongst them.
\end{lemma}

\begin{proof}
We verify that $\uGc{W_1,\spinc}$ is
independent of the ordering of the one-handles. Addition of a pair of
two one-handles to a Heegaard diagram corresponds to adding two
handles to the Heegaard surface to obtain $\Sigma'$, and placing a
pair of new curves $\alpha_1$ and $\beta_1$ and $\alpha_2$ and
$\beta_2$, where $\alpha_i$ is a small isotopic translate of
$\beta_i$, and both pairs are supported inside the two new
handles. From this description, it is clear that the composite of the
two one-handle additions is independent of their order. Handleslide
invariance follows from this. 
\end{proof}

\begin{lemma}
\label{lemma:FourDHS}
Fix a framed link $\Link$ in $Y$, and let $\Link'$ be the framed link
obtained by handleslides amongst the components of $\Link$. Then, the
maps $\uFc_{Y;\Link}$ and $\uFc_{Y;\Link'}$ are equal.
\end{lemma}

\begin{proof}
Let $K_1'$ be obtained from $K_1$ by a handleslide over $K_2$. To
perform this handleslide, one needs a path $\sigma$ joining $K_1$ to $K_2$. 
After the handleslide, there is a natural path $\sigma'$ joining $K_1'$ and $K_2$.

Complete $K_1\cup \sigma\cup K_2$ to a bouquet for $\Link$. Let
$(\Sigma,\alphas,\betas,\gammas,z)$ be the associated triple, where
$\beta_1$ is dual to $K_1$, $\beta_2$ is dual to $K_2$.  We can
complete $K_1'\cup \sigma'\cup K_2$ to a bouquet for $\Link'$ so that
there is a subordinate Heegaard triple $$(\Sigma',\alphas,\betaprs,\gammaprs,z)$$ with the property that
$\beta_2'$ is obtained as a handleslide of $\beta_2$ over $\beta_1$,
$\gamma_1'$ is obtained as a handleslide of $\gamma_1$ over
$\gamma_2$, and all other $\beta_i'=\beta_i$, $\gamma_i'=\gamma_i$ 
(c.f. Figure~\ref{fig:HandleSlides}).

\begin{figure}
\mbox{\vbox{\epsfbox{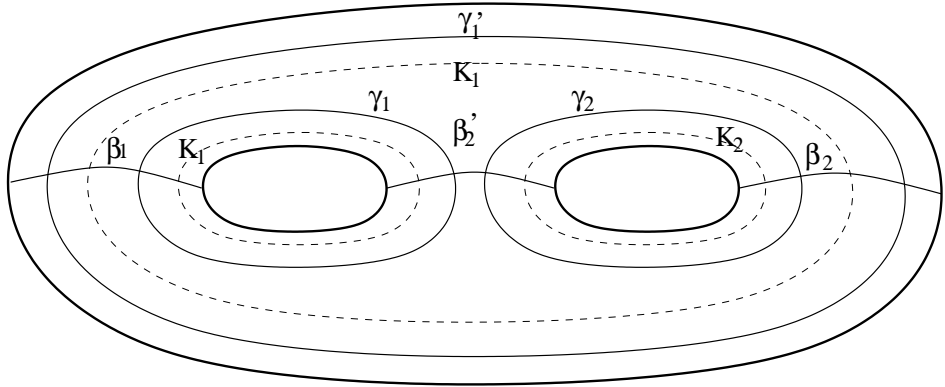}}}
\caption{\label{fig:HandleSlides}
{\bf{Four-dimensional handleslides.}}  An illustration of
four-dimensional handleslides, in the case where the knot has two
components.  The handleslide takes place inside the solid torus
pictured; the dashed lines $K_1$ and $K_2$ constitute the first link,
and $K_1'$ is obtained by sliding $K_1$ over $K_2$. The darker circles
($\beta_1$, $\gamma_1$, $\gamma_2$, $\gamma_1$ and $\beta_2'$)
take place on the boundary of the handlebody, in the Heegaard diagram.
The dual to $K_1'$, $\beta_1'$ is not pictured: it is an isotopic copy
of $\beta_1$; similarly, $\gamma_2'$ is not pictured, as it is an
isotopic translate of $\gamma_2$.}
\end{figure}

We then have the following commutative diagram:
$$
\begin{CD}
\uHFc(Y,\spinc|Y) @>{\otimes \Theta_{\beta,\gamma}}>> \uHFc(Y(\Link),
\spinc|Y(\Link)) \\
@V{\otimes \Theta_{\beta,\beta'}}VV
@V{\otimes \Theta_{\gamma,\gamma}}VV \\
\uHFc(Y,\spinc|Y) @>{\otimes \Theta_{\beta',\gamma'}}>>
\uHFc(Y(\Link'),\spinc|Y(\Link'))
\end{CD}
$$

Commutativity follows from associativity, and the observation that
$$\Fleq{\beta,\beta',\gamma'}
(\Theta_{\beta,\beta'}\otimes \Theta_{\beta',\gamma'},\spinc_0)
=\pm \Theta_{\beta,\gamma'}
= \pm \Fleq{\beta,\gamma,\gamma'}(\Theta_{\beta,\gamma}\otimes 
\Theta_{\gamma,\gamma'},\spinc_0),$$
according to the handleslide invariance of the homology groups.
\end{proof}

\begin{lemma}
\label{lemma:ThreeHandleSlides}
The maps $\uEc{W_3,\spinc}$ are invariant under 
the ordering of the three-handles, and handleslides amongst them.
\end{lemma}

\begin{proof}
Again, we verify independence of the ordering, i.e.
$$\uEc{S_2,\spinc}\circ
\uEc{S_1,\spinc}=\uEc{S_1,\spinc}\circ\uEc{S_2,\spinc}.$$
To see this, we can equip $Y$ with a Heegaard decomposition which
splits off a standard decomposition of $\#^2(S^2\times S^1)$: $(E\#
E,\{\alpha_1,\alpha_2\},\{\beta_1,\beta_2\},z)$. It is clear from the
definitions now that the two composite maps agree.
\end{proof}

\begin{lemma}
\label{lemma:CancelOneTwo}
Let $W_1$ be the cobordism obtained by adding a one-handle to $Y$ and
$W_2$ be the two-handle attached along any knot $\Knot$ which cancels
the one-handle. Then, the induced map 
$$\uFc_{\Knot,\spinc}\circ \uGc{W_1,\spinc}\colon
\uHFc(Y,\spinc,M)\longrightarrow
\uHFc(Y,\spinc,M)$$ 
corresponds to the identity map.
\end{lemma}

\begin{proof}
Fix a Heegaard diagram $(\Sigma,\alphas,\betas,z)$ for $Y_1$. Let
$(E,\alpha,\beta,\gamma,z_0)$ be a Heegaard triple where $E$ is a
genus one surface, $\alpha$ and $\beta$ are exact Hamiltonian
translates of one another, and $\gamma$ is a curve meeting both
transversally in a single intersection point. 
In fact, if we consider the Heegaard triple
$$(\Sigma\#E,\alphas',\betas',\gammas',z')=
(\Sigma,\alphas,\betas,\gammas,z)\#(E,\alpha,\beta,\gamma,z_0),$$
where the $\gammas$ are small
Hamiltonian translates of the $\betas$, then
this  Heegaard triple represents a two-handle $\Knot_0$ which cancels the
one-handle in $Y\times (S^2\times S^1)$  which was just
introduced. The model case calculation in $E$ now shows that
the map induced by the composite induces the same map in homology as
the map 
$$ \CFinf(\alphas,\betas,\spinc)\longrightarrow 
\CFinf(\alphas\cup\{\alpha\},\betas'\cup\{\gamma\},\spinc).$$
given by $$[\x,i]\mapsto [\x'\times \{c\},i],$$ where $c\in
\alpha\cap\gamma$ is the unique intersection point in $E$, and $\x'$ is the intersection point of $\Ta\cap \Tb'$ closest to $\x\in\Ta\cap\Tb$: thus, 
this is equivalent to the map induced by stabilization.

More generally, if we choose another two-handle $\Knot$ which cancels
the given one-handle, we claim that the induced composite is the same.
This can be shown by constructing a Heegaard diagram for which 
the induced $\delta$ belonging to
the framing of $\Knot$ differs from $\gamma$ by a sequence of
isotopies and handleslides across the
$\{\beta_1',...,\beta_g'\}$. This follows from the fact that, if we
surger out the $\{\beta_1',...,\beta_g'\}$ from $\Sigma\# E$, we are
left with a torus, in which $\gamma$ and $\delta$ are two curves which
meet the remaining $\beta$ transversally in a single transverse
intersection point (after isotopies, if necessary). Thus, in this
torus, $\gamma$ and $\delta$ are isotopic.  It follows then that the
two maps $\otimes
\Theta_{\beta,\gamma}$ and $\otimes \Theta_{\beta,\delta}$ (where here
$\deltas=\{\delta_1,...,\delta_g,\delta\}$, and $\delta_i$ for
$i=1,...,g$ are small, Hamiltonian translates of the corresponding
$\beta_i$) which correspond to $F^\circ_{\Knot}$ and
$F^\circ_{\Knot_0}$ agree after we post-compose the second with
multiply the the second with $\otimes \Theta_{\delta,\gamma}$, which
in turn corresponds to a strong equivalence between two Heegaard diagrams
belonging to the same three-manifold $Y_1=Y_3$.
\end{proof}

Dually, we have the following:

\begin{lemma}
\label{lemma:CancelTwoThree}
Let $W_1$ be the cobordism obtained by attaching a two-handle to $Y$
along a framed knot $\Knot$, and let $W_2$ be a three-handle attached along a two-sphere which cancels the knot. Then, the composite
$$\uEc{W_2,\spinc}\circ\uFc_{\Knot,\spinc}\colon\uHFc(Y,\spinc,M)
\longrightarrow \uHFc(Y,\spinc,M)$$
corresponds to the identity map.
\end{lemma}

\begin{proof}
The proof follows from ``turning around'' the proof of 
Lemma~\ref{lemma:CancelOneTwo}.
\end{proof}

\vskip.2cm\noindent{\bf{Proof of Theorem~\ref{thm:CobordismInvariantTwist}}}
Consider the Kirby calculus picture for the cobordism $W$
(see~\cite{KirbyCalculus}, \cite{GompfStipsicz}). Any two such pictures can be
connected by a sequence of pair cancellations and additions, and a
sequence of handleslides (Kirby moves). In fact, since $W$ is a
connected, manifold-with-boundary, $W$ has such a description with
no zero- or four-handles; moreover, any two such
descriptions can be connected through a sequence of Kirby moves which
never introduce new zero- or four-handles. On the other hand, the above
lemmas ensure that the map $\uFc_{W,\spinc}$ is invariant under all the Kirby
moves, and hence, it is a four-manifold invariant.
\qed
\vskip.2cm

For the untwisted case (Theorem~\ref{thm:CobordismInvariant}), 
when we wish to construct the map
$\Fc{W,\spinc}$ without the action of the exterior algebra
of
$H_1(W;\Z)/\Tors$, we can appeal directly to 
Theorem~\ref{thm:CobordismInvariantTwist} above.
Specifically, suppose that $W$ is a cobordism from $W_1$ to $W_2$. 
Taking trivial coefficients $M=\Z$ for $H^1(Y;\Z)$, then we have that the induced module 
$M(W)\cong \Z[H^2(W,Y_2)]$. There is a canonical map of
$\Z[H^1(Y_2;\Z)]$-modules $$\epsilon\colon
\Z[H^2(W,Y_2)]\longrightarrow \Z;$$
so we can define $\Fc{W,\spinc}$ to be the composite
$$\Fc{W,\spinc}= H(\epsilon)\circ\uFc_{W,\spinc}.$$

Equivalently (c.f. Equation~\eqref{eq:TwistAndUntwist}), we can
construct the untwisted maps for cobordisms using the untwisted
triangle construction, and observe that the proof of
Theorem~\ref{thm:CobordismInvariantTwist} (and all its lemmas) adapts
with notational changes to prove
Theorem~\ref{thm:CobordismInvariant}. We adopt this point of view,
when including the action of the one-dimensional homology.

\vskip.2cm
\noindent{\bf{Proof of Theorem~\ref{thm:CobordismInvariant}}}
To construct the invariant as a map
$$\Fc{W,\spinc}\colon \HFc(Y_1,\spinct_1)\otimes \Wedge^* (H_1(W,\Z)/\Tors)\longrightarrow
\HFc(Y_2,\spinct_2),$$
in general, we proceed as follows. 

As before, we split $W=W_1\cup_{Y_1'} W_2\cup_{Y_2'} W_3$, where $W_1$ and
$W_3$ are collections of one- and three-handles respectively, and
$W_2$ is the cobordism induced by the two-handles. Let
$X\subset W_2$ be the four-manifold obtained from a Heegaard triple
for $W_2$ (i.e. $X$ is the four-manifold underlying a
Heegaard triple subordinate to some bouquet for the
link $\Link$ describing $W_2$). 
Clearly, $W$ is obtained from $X$ by adding three- and
four-handles, so in particular the inclusion determines an isomorphism
$H_1(X)\cong H_1(W)$. Thus, we define for each $\gamma\in \Wedge^*
(H_1(W)/\Tors) \cong
\Wedge^* (H_1(X)/\Tors)$, 
$$\Fc{W,\spinc}(\gamma\otimes\xi)
=\Ec{V,\spinc}\circ \Fc{\Link,\spinc}(\gamma\otimes \Gc{W_1,\spinc}),$$
where we now use the extended triangle map for
$\Fc{}$ using the exterior algebra on the one-dimensional homology
(c.f. Lemma~\ref{lemma:ExtendedTriangles}).

The proof of Theorem~\ref{thm:CobordismInvariantTwist} now adapts with
the following observations. First, note that equivalences $\Psi$ of
Heegaard diagrams induce maps on Floer homology which are equivariant
under $H_1(Y)/\Tors$-actions. Thus, we can include the action in the
analogue of 
Theorem~\ref{thm:LinkInvariant}: if $(\Sigma_1,\alphas_1,\betas_1,\gammas_1,z_1)$
and $(\Sigma_2,\alphas_2,\betas_2,\gammas_2,z_2)$ are a pair of Heegaard triples 
subordinate to two bouquets for a link in $Y$, then the following diagram commutes:
$$\begin{CD}
\HFc(\alphas_1,\betas_1,\spinc|Y)\otimes \Wedge^* (H_1(X)/\Tors)
@>{{\Fc{W,\spinc}}^{1}}>>
\HFc(\alphas_1,\gammas_1,\spinc|Y(\Link)) \\
@V{\Psi^\circ_1}VV @V{\Psi^\circ_2}VV \\
\HFc(\alphas_2,\betas_2,\spinc|Y)\otimes \Wedge^* (H_1(X)/\Tors)
@>{{\Fc{W,\spinc}}^{2}}>>
\HFc(\alphas_2,\gammas_2,\spinc|Y(\Link)) \\
\end{CD}
$$
\qed

The treatment of one- and three-handles goes through again with only
notational changes. We observe that the map $$\Gc{W_1,\spinc}\colon
\HFc(Y)\longrightarrow \HFc(Y\#(S^2\times S^1)) $$ is equivariant
under the action of $H_1(Y)/\Tors$. (This was established in the the
precise formulation of
Proposition~\ref{HolDiskTwo:intro:ConnectedSumSphereProd}
of~\cite{HolDiskTwo}, which is
Proposition~\ref{HolDiskTwo:prop:ConnSum} of~\cite{HolDiskTwo}.)
It follows easily that if $W$ and $W'$ differ by the addition of a canceling pair of
one- and two-handles, then the induced maps $\Fc{W,\spinc}$ and $\Fc{W',\spinc}$,
thought of as maps including the exterior algebra action,
agree. The result now follows.
\vskip.2cm

\subsection{Finiteness}

We turn now to Theorem~\ref{thm:Finiteness}.

\vskip.2cm
\noindent{\bf{Proof of Theorem~\ref{thm:Finiteness}.}}
It suffices to consider the case where $W$ consists only of
two-handles. Let $(\Sigma,\alphas,\betas,\gammas,z)$ be a subordinate
Heegaard triple, which we can assume is strongly $\spinc$-admissible
(in the sense of Definition~\ref{HolDiskTwo:def:AdmissibleTriple}
of~\cite{HolDiskTwo}, see also Subsection~\ref{subsec:HolTriangles} of
the present paper). Note that this admissibility condition depends on
$\spinc$ only through its restriction to the boundary.  We wind
further to achieve the following additional admissibility condition:
we arrange that each triply-periodic domain with $n_z(\PerDom)=0$ has
both positive and negative coefficients. It is easy to see that under
this hypothesis, for each $\x\in\Ta\cap\Tb$, $\Theta\in\Tb\cap\Tc$ and
$\y\in\Ta\cap\Tc$ and for each integer $j$, there are only finitely
many homotopy classes of maps $\psi\in\pi_2(\x,\y,\w)$ with
$n_z(\psi)=j$ which support holomorphic representatives. Both
finiteness conditions are now clearly satisfied.
\qed
\vskip.2cm

\subsection{Composition laws revisited}

The composition laws for cobordism invariants follows 
readily from the
definitions, together with the version of the composition law already
stated for links (c.f. Proposition~\ref{prop:CompositeLink}).
But first, 
we must verify that we can commute two-handle additions past
one-handle additions; and similarly, three-handle additions commute
past two-handle additions.

\begin{prop}
\label{prop:CommuteOneTwo}
Let $\Link\subset Y$ be a link in a three-manifold, and let $\Link'$
be the corresponding link in the three-manifold $Y'=
Y\#(\#^n(S^2\times S^1))$ obtained by attaching one-handles. Then, the following diagram commutes:
$$
\begin{CD}
\uHFc(Y,\spinc|Y) @>{\Fc{\Link,\spinc}}>> \HFc(Y(\Link),\spinc|Y(\Link)) \\
@V{\uGc{V,\spinc}}VV @V{\uGc{V(\Link),\spinc}}VV \\
\uHFc(Y',\spinc|Y') @>{\Fc{\Link',\spinc}}>> \HFc(Y'(\Link'),\spinc|Y'(\Link')).
\end{CD}
$$
\end{prop}

\begin{proof}
This follows immediately from the gluing result for holomorphic
triangles, Theorem~\ref{thm:GlueTriangles}, together with the model
case of a torus with three (homologically non-trivial, simple, closed)
curves in the torus which are exact Hamiltonian translates of one
another. (In fact, we have verified this proposition in a special case
when we showed in Theorem~\ref{thm:OneHandles} that $\Gc{}$ commutes
with the maps induced by equivalences of Heegaard diagrams;
but the case where the horizontal map is induced 
more generally by a Heegaard triple, the proof is
no different.)
\end{proof}

Similarly, we have the following:

\begin{prop}
\label{prop:CommuteTwoThree}
Let $\Link'\subset Y'$ be a link in a three-manifold, which is
disjoint from a non-separating two-sphere $S\subset Y'$, then the following diagram
commutes: 
$$
\begin{CD}
\uHFc(Y',\spinc|Y') @>{\Fc{\Link',\spinc}}>> \HFc(Y'(\Link),
\spinc|Y'(\Link')) \\
@V{\uEc{U,\spinc}}VV @V{\uEc{U(\Link),\spinc}}VV \\
\uHFc(Y,\spinc|Y) @>{\Fc{\Link,\spinc}}>> \HFc(Y(\Link),
\spinc|Y(\Link)),
\end{CD}
$$
where $Y$ is the result of adding a three-handle to $Y'$ along
$S$, and $\Link\subset Y$ is the link induced from $\Link'$.
\end{prop}

\begin{proof}
We construct the Heegaard triple
$(\Sigma',\alphas',\betas',\gammas',z')$ for $Y'$ so that the last
curves are $\alpha_{g+1}$, $\beta_{g+1}$ and $\gamma_{g+1}$ are small
Hamiltonian isotopic translates of one another, all supported inside a
summand $E$ of $\Sigma'=\Sigma\# E$, so that the
$\alphas=\{\alpha_1,...,\alpha_g\}$, $\betas=\{\beta_1,...,\beta_g\}$
and $\gammas=\{\gamma_1,...,\gamma_g\}$ are supported inside
$\Sigma$. Then, $(\Sigma,\alphas,\betas,\gammas,z)$ is a Heegaard
triple for $\Link\subset Y$, and the result follows from the gluing
for holomorphic triangles with the model calculation of
$(E,\alpha_{g+1},\beta_{g+1},\gamma_{g+1},z_0)$ as before.
\end{proof}

\vskip.2cm
\noindent{\bf{Proof of Theorem~\ref{thm:CompositionTwist}.}}
The proof is now follows easily from 
Propositions~\ref{prop:CompositeLink}, \ref{prop:CommuteOneTwo}, and
\ref{prop:CommuteTwoThree}.
\qed

The need to sum over $\SpinC$ structures for the composition law in
the untwisted case follows from the corresponding fact in the
associativity of the triangle construction. This results in the following 
analogue of Proposition~\ref{prop:CompositeLink} in the untwisted case:

\begin{prop}
\label{prop:CompositeLinkUntwist}
Suppose that $W_1$ and $W_2$ are composed entirely of two-handles, then
there are choices of sign for $\Fc{W_1,\spinc_1}$ and
$\Fc{W_2,\spinc_2}$ with the property
that
$$\Fc{W_2,\spinc_2}\circ \Fc{W_1,\spinc_1}
=
\sum_{\{\spinc\in\SpinC(W_1\# W_2)\big|\spinc|W_1=\spinc_1,
\spinc|W_2=\spinc_2\}} 
\pm \Fc{W_1\# W_2,\spinc}.$$
More generally, given $\zeta_1\in \Wedge^*(H_1(W_1;\Z)/\Tors)$ and $\zeta_2\in \Wedge^*(H_1(W_2,\Z)\Tors)$, we have that
$$\Fc{W_2,\spinc_w}(\zeta_2\otimes \Fc{W_1,\spinc_1}(\zeta_1\otimes \cdot))
= 
\sum_{\{\spinc\in\SpinC(W_1\# W_2)\big|\spinc|W_1=\spinc_1,
\spinc|W_2=\spinc_2\}} 
\Fc{W_1\# W_2,\spinc}(\zeta_3\otimes \cdot),$$
where $\zeta_3\in\Wedge^*(H_1(W_1\# W_2)/\Tors)$ is the image of $\zeta_1\otimes \zeta_2$ under the natural map. 
\end{prop}

\begin{proof}
Apply the proof of Proposition~\ref{prop:CompositeLink}, only using
associativity in the untwisted case. Let
$(\Sigma,\alphas,\betas,\gammas,z)$ and
$(\Sigma,\alphas,\gammas,\deltas,z)$ be the two Heegaard triples
constructed in the proof.  To verify that the action of $H_1$ is
respected as claimed, we represent the action of $\zeta_1\in
H_1(X_{\alpha,\beta,\gamma};\Z)\cong H_1(W_1)$ by the action of a class from
$H_1(Y_{\alpha,\beta})$, which in turn is a constraint from
$H^1(\Ta)$. It follows immediately that: $$\Fc{W_2,\spinc_w}(1\otimes
\Fc{W_1,\spinc_1}(\zeta_1\otimes \cdot)) =
\sum_{\{\spinc\in\SpinC(W_1\# W_2)\big|\spinc|W_1=\spinc_1,
\spinc|W_2=\spinc_2\}} 
\Fc{W_1\# W_2,\spinc}(\zeta_1\otimes \cdot)$$
(where the element $\zeta_1$ appearing on the right-hand-side is
to be interpreted as the image of $\zeta_1$ inside $H_1(W_1\# W_2;\Z)$).
Classes $\zeta_2$ coming from $H_1(X_{\alpha,\gamma,\delta})\cong
H_1(W_2)$ can be represented by classes from $H_1(Y_{\alpha,\delta})$, 
to get the corresponding equation that:
$$\Fc{W_2,\spinc_w}(\zeta_2\otimes
\Fc{W_1,\spinc_1}(1\otimes \cdot)) =
\sum_{\{\spinc\in\SpinC(W_1\# W_2)\big|\spinc|W_1=\spinc_1,
\spinc|W_2=\spinc_2\}} 
\Fc{W_1\# W_2,\spinc}(\zeta_2\otimes \cdot).$$
\end{proof}

\vskip.2cm
\noindent{\bf{Proof of Theorem~\ref{thm:Composition}.}}
With Proposition~\ref{prop:CompositeLinkUntwist} replacing
Proposition~\ref{prop:CompositeLink}, the proof of
Theorem~\ref{thm:CompositionTwist} applies.
\qed

\section{Duality and Conjugation invariance}
\label{sec:Duality}

In the present section, we describe duality for the maps, and
the closely related conjugation invariance.

\subsection{Duality}
We discuss duality for the maps associated to cobordisms, proving
Theorem~\ref{thm:Duality} and \ref{thm:Conjugation}. But first, we
recall some aspects of duality for the homology groups of a
three-manifold.

If $Y$ is a three-manifold, then we can think of $\SpinC(Y)$
structures as equivalence classes of nowhere vanishing vector
fields. Since a nowhere vanishing vector field over $Y$ can also be
viewed as a vector field over $-Y$, we get a naturally induced bijection
$$\SpinC(Y)\cong \SpinC(-Y).$$ 

In terms of Heegaard diagrams, if the diagram
$(\Sigma,\alphas,\betas,z)$ describes $-Y$, then
$(-\Sigma,\alphas,\betas,z)$ describes $Y$. Under the above
identification, if $\x\in\Ta\cap\Tb$ is an intersection point, then
the $\SpinC$ structure induced by $\x$ and the basepoint $z$ for the
first Heegaard diagram agrees with that induced by $\x$ and $z$ for
the second diagram.

We define a pairing 
$$\langle,\rangle \colon
\CFinf(Y,\spinc)\otimes \CFinf(-Y,\spinc)\longrightarrow \Z,$$ 
by the formula
$$\langle [\x,i],[\y,j]\rangle = \left\{\begin{array}{ll}
1 & {\text{if $\x=\y$ and $i+j+1=0$}} \\
0 & {\text{otherwise}}
\end{array}\right.,$$
where we are using $(\Sigma,\alphas,\betas,z)$ for $Y$ and
$(-\Sigma,\alphas,\betas,z)$ for $-Y$ as above.

\begin{lemma}
Under the above pairing, we have the identities:
\begin{eqnarray}
\langle \xi,\partial^\infty_{-Y} \eta \rangle
&=& \langle \partial^\infty_Y \xi,\eta \rangle, 
\label{eq:HomologyDual}
\\
\langle \xi, U\eta \rangle &=&  \langle U \xi,\eta \rangle.
\label{eq:UDual}
\end{eqnarray}
\end{lemma}

\begin{proof}
Now, precomposition of the reflection of the disk across the real axis
gives a map from $\pi_2(\x,\y)$ to $\pi_2(\y,\x)$ which identifies the
$\sj$-holomorphic maps in $\Sigma$ with the $(-\sj)$-holomorphic maps in
$-\Sigma$ (indeed, if $J_t$ is a one-parameter family of
almost-complex structures over $\Sym^g(\Sigma)$, then this reflection
identifies the $J_t$-pseudo-holomorphic maps with the
$(-J_t)$-pseudo-holomorphic maps in $\Sym^g(-\Sigma)$). The first
formula follows.
The second formula is immediate.
\end{proof}

According to Equation~\eqref{eq:HomologyDual}, the pairing we have
defined descends to give pairings
\begin{eqnarray*}
\langle,\rangle \colon \HFinf(Y,\spinc)\otimes \HFinf(-Y,\spinc)&\longrightarrow& \Z, \\
\langle,\rangle \colon \HFp(Y,\spinc)\otimes \HFm(-Y,\spinc)&\longrightarrow& \Z.
\end{eqnarray*}
Indeed, according to Equation~\eqref{eq:UDual}, for this second
pairing, $\HFmred(-Y,\spinc)$ pairs trivially with the image of
$\HFinf(Y,\spinc)$ inside $\HFpred(Y,\spinc)$; thus, we obtain an
induced pairing $$\langle,\rangle \colon
\HFpred(Y,\spinc)\otimes \HFmred(-Y,\spinc)\longrightarrow \Z.$$

In the case where $c_1(\spinc)$ is torsion, the pairing (on the level
of chain complexes) induces an isomorphism between the chain complex
for $\CFp(Y,\spinc)$ and the cochain complex for
$\CFm(-Y,\spinc)$. 

If $W$ is a cobordism from $Y_1$ to $Y_2$, then, $W$ can also be
thought of as a cobordism from $-Y_2$ to $-Y_1$. We show that
$$F^+_{W,\spinc}\colon \HFp(Y_1,\spinc|Y_1)\longrightarrow
\HFp(Y_2,\spinc|Y_2)$$ is adjoint (under the above pairing) to the map
$$F^-_{W,\spinc}\colon \HFm(-Y_2,\spinc|Y_2)\longrightarrow
\HF(-Y_1,\spinc|Y_1),$$ as stated in Theorem~\ref{thm:Duality}.

\vskip.3cm
\noindent{\bf{Proof of Theorem~\ref{thm:Duality}}.}
When viewing $W$ as a cobordism from $-Y_2$ to $-Y_1$, three-handles are viewed
as one-handles, the one-handles are viewed as three-handles, and the collection of two handles are viewed as another dual collection.

In particular, we consider the case where $W$ consists of a single
one-handle, attached to $Y_1=Y$ to obtain $Y_2=Y\#(S^1\times S^2)$.
Suppose that $Y_1$ is represented by $(\Sigma,\alphas,\betas,z)$, and
then that $Y_2=(\Sigma,\alphas,\betas,z)\#(E,\alpha_0,\beta_0,z_0)$
where $(E_0,\alpha_0,\beta_0,z_0)$ is a standard Heegaard splitting of
$S^1\times S^2$. Then, we can alternatively think of the cobordism as
going from $-Y_2 =
(-\Sigma,\alphas,\betas,z)\#(-E,\alpha_0,\beta_0,z_0)$ to
$(-\Sigma,\alphas,\betas,z)$ by attaching a three-handle to the
two-sphere specified by the pair $\alpha_0$ and $\beta_0$. Now, as
defined in Subsection~\ref{subsec:OneThreeHandles}, $F^+_{W}$ is the map
taking $[\x,i]$ to $[\x\times\{\theta\},i]$, where $\theta$ is the
intersection point of $\alpha_0\cap\beta_0$ which has higher Maslov
index. But viewed as an intersection point for
$(-E,\alpha_0,\beta_0,z_0)$, $\theta$ has lower Maslov index. It
follows then immediately from the definition of the map on
three-handles (c.f. Subsection~\ref{subsec:OneThreeHandles}), that 
$F^+_{W,\spinc}$ as a map from the homology $Y_1$ to $Y_2$
is dual to $F^-_{W,\spinc}$ as a map from the homology of $-Y_2$ to $-Y_1$.
The result follows with notational changes 
when $W$ is obtained as a union of one-handles (or three-handles).

We now focus on the case of two-handles. Let
$(\Sigma,\alphas,\betas,\gammas,z)$ be a Heegaard triple subordinate
to some bouquet for the link $\Link\subset Y$. Then, the Heegaard
triple $(-\Sigma,\alphas,\gammas,\betas,z)$ is a Heegaard triple
subordinate some link $\Link'\subset -Y(\Link)$, on which a surgery
gives $-Y$. We must compare holomorphic triangles for these two Heegaard triples.

Observe that there is a unique antiholomorphic involution $R$ of the
triangle $\Delta$ with edges $\alpha$, $\beta$, and $\gamma$, which
switches the $\beta$ and $\gamma$ edges and maps the $\alpha$ edge to
itself. Precomposing $\psi'$ with $R$ induces an identification for
each $\x\in\Ta\cap\Tb$, $\w\in\Tb\cap\Tc$, $\y\in\Ta\cap\Tc$ of
classes $$\pi_2(\x,\w,\y)\cong\pi_2(\y,\w,\x)$$ which carries all the
$\Sigma$-holomorphic to $-\Sigma$ holomorphic data:
$\Mas_{\Sigma}(\psi)=\Mas_{-\Sigma}(\psi\circ R)$ and
$\ModFlow_{\Sigma}(\psi)=\ModFlow_{-\Sigma}(\psi\circ R)$. Since 
$R$ reverses orientation, we have that
$n^{\Sigma}_z(\psi\circ R)=n_z^{-\Sigma}(\psi)$.

Thus,
\begin{eqnarray*}
\langle F^+_{\Link,\spinc}([\x,i]\otimes [\w,j]),[\y,k]\rangle
&=&\sum_{\{\psi\in\pi_2(\x,\w,\y)\big|\Mas_{\Sigma}(\psi)=0,\spinc_z(\psi)=\spinc,
n_z(\psi)=i+j+k+1\}}\#\ModFlow_\Sigma(\psi) \\
&=& 
\sum_{\{\psi\in\pi_2(\y,\w,\x)\big|\Mas_{-\Sigma}(\psi)=0,\spinc_z(\psi)=\spinc,
n_z(\psi)=i+j+k+1\}}\#\ModFlow_{-\Sigma}(\psi) \\
&=& \langle [\x,i], F^+_{\Link',\spinc}([\y,k]\otimes[\w,j]).
\end{eqnarray*}

Observe next that if $\Theta$ is some chain in
$\CFleq(\betas,\gammas,z)$ representing a generator for the highest
dimensional non-trivial homology of $\HFleq(\#^n(S^1\times
S^2),\spinc_0)$ (where we are using the Heegaard diagram
$(\Sigma,\betas,\gammas,z)$), then $\Theta$ can be viewed as a chain
for $\CFleq(\gammas,\betas,z)$ (using $-\Sigma$ as the dividing
surface); as such it will also represent a generator for the highest
non-trivial homology of $\HFleq(\#^n(S^1\times S^2),\spinc_0)$ (this
follows from the conjugation invariance of the three-dimensional
groups, together with the observation that $\spinc_0$ is fixed under
conjugation). The theorem follows similarly when $W$ is composed of a
collection of two-handles. The general case follows, combined with the
remarks on one- and three-handles at the beginning of the proof.
\qed
\vskip.3cm

\subsection{Conjugation invariance}

Recall that if $(\Sigma,\alphas,\betas,z)$ is a Heegaard decomposition
for $Y$, then $(-\Sigma,\betas,\alphas,z)$ also represents $Y$, and if
$\x\in\Ta\cap\Tb$ is an intersection point, then the two $\SpinC$
structure associated to $\x$ using the two Heegaard diagrams are
conjugates of one another. This induces isomorphisms $$\Conjugate_Y^\circ
\colon \HFc(Y,\spinct)\longrightarrow \HFc(Y,{\overline \spinct}).$$

We wish to prove the corresponding conjugation invariance of maps
associated to cobordisms. But, first we prove a technical lemma,
regarding the invariant of maps associated to framed links.

Consider a framed link $\Link\subset Y$. We
call a Heegaard triple $(\Sigma,\deltas,\alphas,\betas,z)$ {\em
left-subordinate} to a bouquet for the link if $Y_{\delta,\alpha}\cong
\#^\ell(S^1\times S^2)$, $Y_{\alpha,\beta}\cong Y$, and
$Y_{\delta,\alpha}\cong Y(\Link)$. We define a link invariant 
$$K^{\circ}_{\Link,\spinc}\colon
\HFc(Y,\spinc|_Y)
\longrightarrow \HFc(Y(\Link),\spinc|_{Y(\Link)})$$
using a left-subordinate Heegaard triple, using the formula
$$K^{\circ}_{\Link,\spinc}(\xi)=f^\circ(\Theta_{\delta,\alpha}\otimes
\xi,\spinc),$$ where $\Theta_{\delta,\alpha}$ is, once again, a top-dimensional generator
of $\HFleq(\#^\ell(S^2\times S^1),\spinc_0)$.
This agrees with $F^+_{\Link,\spinc}$, 
as defined in
Subsection~\ref{subsec:TwoHandles}, according to the following:

\begin{lemma}
\label{lemma:LeftTheta}
The invariant of the link $K_{\Link,\spinc}^+$ agrees with
$F_{\Link,\spinc}^+$ up to sign.
\end{lemma}

\begin{proof}
For each link, we can find a Heegaard quadruple
$$(\Sigma,\deltas,\alphas,\betas,\gammas,z)$$ with the property that
$(\Sigma,\deltas,\alphas,\betas,z)$ is left-subordinate to
$\Link\subset Y$, $(\Sigma,\alphas,\betas,\gammas,z)$ is subordinate
(in the usual sense) to the link $\Link\subset Y$,
$(\Sigma,\deltas,\betas,\gammas,z)$ is subordinate to a cobordism from
$Y$ to $Y\#^n(S^1\times S^2)$ (obtained as zero-surgery on an
$n$-component framed link), and $(\Sigma,\deltas,\alphas,\gammas)$ is
left-subordinate to the same cobordism. Moreover, the Heegaard diagram
$(\Sigma,\deltas,\gammas,z)$ is equipped with $n$ distinguished two-spheres
which can be cancelled (by adding three-handles) to obtain a Heegaard diagram for $Y$ back.
We then consider the diagram $$\begin{CD}
\HFp(\alphas,\betas,\spinc|Y_{\alpha,\beta})
@>{f^+(\cdot\otimes \Theta_{\beta,\gamma},\spinc)}>>
\HFp(\alphas,\gammas,\spinc|Y_{\alpha,\gamma}) \\
@V{f^+(\Theta_{\delta,\alpha}\otimes \cdot,\spinc)}VV 
@VV{f^+(\Theta_{\delta,\alpha}\otimes\cdot,\spinc)}V \\
\HFp(\deltas,\betas,\spinc|Y_{\delta,\beta}) 
@>{f^+(\cdot\otimes \Theta_{\beta,\gamma},\spinc)}>>
\HFp(\deltas,\gammas,\spinc|Y_{\delta,\gamma}),
\end{CD}$$
which commutes by the usual associativity. 
Post-composing with the three-handles, (which cancel the two-handles whether
they are added ``from the left'' according to Lemma~\ref{lemma:CancelTwoThree}, 
as in the bottom row above,
or ``from the right'' according to a suitably modified version of that lemma, 
as in the right column), we see that the result follows.

We have pictured the Heegaard quadruple in the case where the link has
a single component, c.f. Figure~\ref{fig:JInvariance}.

\begin{figure}
\mbox{\vbox{\epsfbox{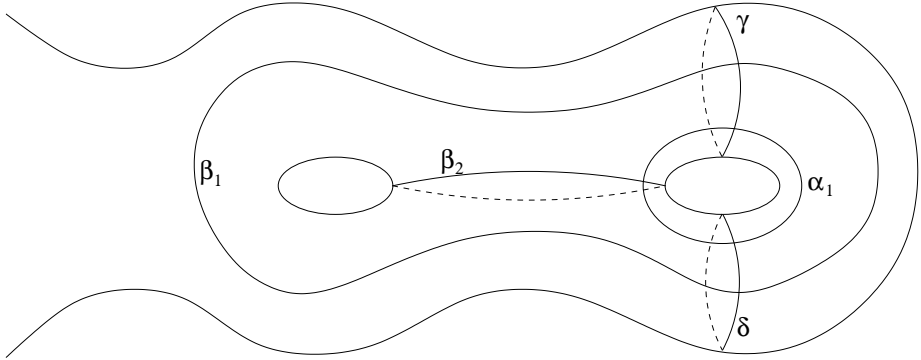}}}
\caption{\label{fig:JInvariance}
{\bf{Heegaard quadruple from Lemma~\ref{lemma:LeftTheta}.} } Let
$\delta_1=\delta$ and $\delta_i$ be a small isotopic translate of
$\alpha_i$ for $i>1$; let $\gamma_1=\gamma$ and $\gamma_i$ be a small
isotopic translate of $\beta_i$ for $i>1$.  Then,
$(\Sigma,\deltas,\alphas,\betas,\gammas,z)$ is a Heegaard quadruple as
required by Lemma~\ref{lemma:LeftTheta}. The curves $\gamma$ and
$\delta$ give the two-sphere along which the three-handle can be added.
(Here, our link has one component; in the more general case, we graft
$n$ copies of this picture.)}
\end{figure}
\end{proof}

We can now prove the conjugation invariance of the maps associated to cobordisms.

\vskip.2cm
\noindent{\bf{Proof of Theorem~\ref{thm:Conjugation}.}}
The result is obvious for cobordisms consisting only of one- and
three-handles.  Thus, we focus on the case of two-handles:
i.e. $W=W(\Link)$ for some framed link $\Link\subset Y$ inside
$Y=Y_1$. Let $(\Sigma,\alphas,\betas,\gammas,z)$ be a Heegaard triple
subordinate to some bouquet for the link. Then,
$(-\Sigma,\gammas,\betas,\alphas)$ represents the same oriented
four-manifold, and indeed, it is left-subordinate to the link
$\Link\subset Y$.  We claim that for each chain $\xi\in\CFc(\alphas,\betas,\spinc|Y)$,
$$f^\circ_{\alpha,\beta,\gamma}(\xi\otimes \Theta_{\beta,\gamma})
=\Conjugate_{Y(\Link)}\circ
f^\circ_{\gamma,\beta,\alpha}(\Theta_{\gamma,\beta}\otimes
\Conjugate_{Y}(\xi)).$$  
(Here, as usual, $\Theta_{\beta,\gamma}$ is a generator for the top
non-zero homology group of $\HFleq(\#^n(S^2\times S^1),\spinc_0)$, and
so is $\Conjugate_{\#^n(S^2\times
S^1)}(\Theta_{\beta,\gamma})=\Theta_{\gamma,\beta}$.)  This comes from
an identification between holomorphic triangles for the two Heegaard
triples which is defined by pre-composing triangles
$\psi\in\pi_2(\x,\w,\y)$ by the anti-holomorphic involution of the
$\rho$ triangle which preserves the $\beta$ edge and switches the $\alpha$-
and $\gamma$-edges.  

We claim that this pre-composition induces conjugation on the level of
(four-dimensional) $\SpinC$ structures. More precisely, if
$$\Phi\colon X(\Sigma,\alphas,\betas,\gammas,z)\longrightarrow
X(-\Sigma,\gammas\betas,\alphas,z)$$ is the obvious diffeomorphism of
four-manifolds associated to Heegaard triples, then
$$\Phi^*(\spinc_z(u\circ\rho))\cong \Conjugate(\spinc_z(u)).$$ This
can be seen on the level of the two-plane fields which determine the
$\SpinC$ structures
(see~\ref{HolDiskTwo:subsec:TrianglesAndSpinCStructures} of~\cite{HolDiskTwo}): the (singular) two-plane fields
$\Phi^*(\spinc_z(u\circ\rho))$ and $\spinc_z(u)$ agree, only they have
opposite orientations.

Conjugation invariance now follows from these observations, together 
with Lemma~\ref{lemma:LeftTheta}.
\qed

\section{Blowing up}
\label{sec:BlowUp}

The ``blow-up formula'' stated in Theorem~\ref{thm:BlowUp} can be
reduced to a calculation in 
a genus one surface.

Consider a Heegaard triple $(E,\{\alpha_0\},\{\beta_0\},\{\gamma_0\},z_0)$ where
$\alpha_0$, $\beta_0$, $\gamma_0$ are three curves each meeting
pairwise in a single intersection point, as pictured in
Figure~\ref{fig:cptwo} (with subscripts dropped). This Heegaard triple
is subordinate to the unknot in the three-sphere with framing $-1$:
i.e. the underlying four-manifold is  $\mCP$ punctured in three points.

Let  $x_0=\alpha_0\cap\beta_0$,
$y_0=\beta_0\cap\gamma_0$, and
$w_0=\alpha_0\cap\gamma_0$, we have that
$\pi_2(x_0,y_0,w_0)\cong \Z$, and each homotopy class represents a
different $\SpinC$ structure $\spinc$. Choosing the basepoint as in
Figure~\ref{fig:cptwo}, we have a family of homotopy classes $\{\psi^\pm_k\}$
indexed by a sign and a positive integer $k$, with
\begin{eqnarray*}
\Mas(\psi^\pm_k)&=&0, \\
n_z(\psi^\pm_k)&=&\frac{k(k-1)}{2} \\
\end{eqnarray*}
(compare Proposition~\ref{HolDiskTwo:prop:HoClassesCancel}
of~\cite{HolDiskTwo}).  Each homotopy class has a unique, smooth
holomorphic representative, so that $\#\ModFlow(\psi^\pm_k)=\pm 1$.

We will use the following:

\begin{lemma}
\label{lemma:CPtwoBar}
Let $\{\psi^\pm_k\}_{k=1}^\infty$ be the family of triangles described above, 
and let $E\in H_2(\mCP;\Z)$ be a generator. Then, 
$$\Big|\langle c_1(\spinc_z(\psi_k)),[E]\rangle\Big| = 2k+1.$$
\end{lemma}

The proof is deferred to Subsection~\ref{subsec:CPtwoBar}, after we
establish  a more general formula
for the first Chern class of the $\SpinC$ structure underlying a
Whitney triangle, in terms of combinatorial data in a Heegaard triple
(Proposition~\ref{prop:COneFormula}).

\vskip.3cm
\noindent{\bf{Proof of Theorem~\ref{thm:BlowUp}.}}
The cobordism $W$ can be expressed by a single two-handle (an 
unknot in $Y$ with framing $-1$). 
Let $(\Sigma,\alphas,\betas,z)$ be a genus $g$ Heegaard diagram for 
$Y$. Let $\gammas$ be small isotopic translates of the $\betas$. 
Then, the connected sum of $(\Sigma,\alphas,\betas,\gammas,z)$ with
$(E,\{\alpha_0\},\{\beta_0\},\{\gamma_0\},z_0)$
is a Heegaard triple which is subordinate to the unknot with framing $-1$. 
The theorem then follows immediately from the above observations, together 
with the gluing theorem for holomorphic triangles
Theorem~\ref{thm:GlueTriangles}.
\qed
\vskip.3cm

\subsection{The first Chern class formula}
\label{subsec:COneFormula}

We give a formula for calculating the evaluation of the first Chern
class of a $\SpinC$ structure over $X$, as specified by a Whitney
triangle and a base-point, on a two-dimensional homology class in $X$,
as specified by a triply-periodic domain. This quantity is expressed
purely in terms of data on $\Sigma$ (compare~\cite{HolDiskTwo}, for an
analogous calculation in dimension three). This formula is used to
establish Lemma~\ref{lemma:CPtwoBar} above.

To state the formula, we define certain quantities associated to
triply-periodic domains and classes in $\pi_2(\x,\y,\w)$: the Euler
measure, and the dual spider number.

Fix a triply-periodic domain ${\mathcal P}$, thought of as a two-chain
in $\Sigma$ which spans some triple
$(a,b,c)\in\SpanA\oplus\SpanB\oplus\SpanC$. We define the {\em Euler
measure} $\EulerMeasure({\mathcal P})$ as follows. 
Find a representative
$$\Phi\colon F\longrightarrow \Sigma$$
representing $\mathcal P$, where $F$ is a two-manifold with boundary,
and where $\Phi$ which is an immersion in a neighborhood of
$\partial F$. The line bundle $\Phi^*(T\Sigma)$ has a canonical
trivialization over $\partial F$, since $\Phi$ induces an isomorphism
$$TF \cong \Phi^*(T\Sigma),$$ and $TF$ is canonically trivialized near
$\partial F$ (using the outward normal orientation on $F$). We define
$\EulerMeasure({\mathcal P})$ to be the Euler number of $\Phi^*(T\Sigma)$,
relative to this trivialization at $\partial F$,
$$\EulerMeasure({\mathcal P})=\langle c_1(\Phi^*(T\Sigma),\partial F), F\rangle.$$

\begin{lemma}
\label{lemma:CalcEulerMeasure}
Writing a triply-periodic domain 
$${\mathcal P}=\sum n_i {\mathcal D}_i,$$
we can calculate the Euler measure by the formula
$$\EulerMeasure({\mathcal P})=\sum_i n_i \left(\chi ({\mathrm{int}} {\mathcal D}) -\frac{1}{4}(\#{\text
{corner points of $F$}})\right),$$
where the corner points are to be counted with multiplicity.
\end{lemma}

\begin{proof}
Endow $\Sigma$ with a Riemannian metric $g$ for which the
$\alphas$, $\betas$, and $\gammas$ are simultaneously geodesics, any
two of which meet at right angles, and let  $\Phi\colon F \longrightarrow
\Sigma$ be a branched cover representing 
${\mathcal P}$. The relative 
Chern class of $\Phi^*(T\Sigma)$ is calculated by $\int_F
\Phi^*(K_g)$, where $K_g$ is the Gaussian curvature of the metric
$g$; which in turn can be calculated locally over each ${\mathcal
D}$ using the Gauss-Bonet formula to give the formula stated above.
\end{proof}

Next, we set up the dual spider number of a Whitney triangle
$$u\colon \Delta \longrightarrow \Sym^g(\Sigma)$$
and a triply-periodic domain $\PerDom$. 

Note first that the orientation of $\Sigma$, and the orientations on
all the attaching circles $\alphas$, $\betas$, and $\gammas$ naturally
induce ``inward'' normal vector fields to the attaching circles
(i.e. if $\gamma\colon S^1\longrightarrow \Sigma$ is a unit speed
immersed curve, this inward normal vector is given by
$J\frac{d\gamma}{dt}$).  Let $\alpha_i'$, $\beta_i'$, and $\gamma_i'$
denote copies of the corresponding attaching circles $\alpha_i'$,
$\beta_i'$, and $\gamma_i'$, translated slightly in these normal
directions. Let $\alphaprs$, $\betaprs$, and $\gammaprs$ denote the
corresponding $g$-tuples, and $\TaPr$, $\TbPr$, and $\TcPr$ be the
corresponding tori in $\Sym^g(\Sigma)$. By construction, then,
$u(e_{\alpha})$ misses $\TaPr$, $u(e_{\beta})$ misses $\TbPr$, and
$u(e_{\gamma})$ misses $\TcPr$.

Let $x\in\Delta$ be a point in the interior, chosen in general
position, so that the $g$-tuple $u(x)$ misses all of $\alphaprs$,
$\betaprs$, and $\gammaprs$. Choose three paths $a$, $b$, and $c$ from
$x$ to $e_{0}$, $e_{1}$, and $e_{2}$ respectively. The central
point $x$ and the three paths $a$, $b$, and $c$ is called a {\em dual
spider}. We can think of the paths $a$, $b$, and $c$ as one-chains
in $\Sigma$.
Recall that $\partial\PerDom$ has three types of boundaries: the
$\alpha$, $\beta$, and $\gamma$ boundaries, which we denote
$\partial_\alpha \PerDom$, $\partial_\beta \PerDom$, and
$\partial_\gamma\PerDom$.
Let 
$\partial'_\alpha\PerDom$, $\partial'_\beta\PerDom$, and
$\partial'_\gamma\PerDom$ respectively denote the one-chains obtained
by translating the corresponding boundary components using the induced
normal vector fields.
The dual spider number of $u$ and $\PerDom$ is defined by
$$\Spider(u,\PerDom)=n_{u(x)}(\PerDom)+\#(a\cap\partial'_\alpha\PerDom)+\#(b\cap\partial'_\beta\PerDom)+\#(c\cap
\partial'_\gamma\PerDom).$$

It is elementary to verify that $\sigma(u,\PerDom)$ is
depends only on the homotopy class of $u$ and the periodic domain
$\PerDom$; in particular, it is independent of the choice of dual
spider used in its definition.

\begin{prop}
\label{prop:COneFormula}
Fix a Whitney triangle $u$, a base-point $z$, and a triply-periodic
domain $\PerDom$, whose boundary represents the two-dimensional
homology class $H(\PerDom)\in H_2(X;\Z)$. Then, 
\begin{equation}
\label{eq:COneFormula}
\langle c_1(\spinc_z(u)), H(\PerDom)\rangle
=
\EulerMeasure(\PerDom) + \#(\partial\PerDom)-2n_z(\PerDom)+2 \Spider(u,\PerDom).
\end{equation}
\end{prop}

\begin{proof}
To calculate the first Chern class, we adopt the notation from
Subsection~\ref{HolDiskTwo:subsec:TrianglesAndSpinCStructures}
of~\cite{HolDiskTwo}, where the $\SpinC$ structure belonging to a
triangle is constructed. The construction is done by constructing a
two-plane field ${\mathcal L}$ in $X$, which is defined away from a
finite collection of balls in the four-manifold $X$. We make use of
the complex decomposition $TX\cong {\mathcal L}\oplus {\mathcal
L}^\perp$ (which holds wherever ${\mathcal L}$ is defined), so that
$$c_1(TX)\cong c_1({\mathcal L}\oplus {\mathcal L}^\perp).$$

Consider the representative for $H(\PerDom)$ constructed earlier. We can
assume that this representative misses the locus in $X$ where the
two-plane field ${\mathcal L}$ is undefined. We think of its
projection to $\Delta$ is a dual spider for $u$, with arcs $a$, $b$,
and $c$.

Note that ${\mathcal L}^\perp$ has a
nowhere-vanishing vector field defined on the support of our
representative, so its first Chern class evaluates trivially. (Such a
vector field can be obtained by pulling back a vector field defined
over the triangle which does not vanish on the support of the dual spider.)
Thus,
$$\langle c_1(\spinc_z(u)),H(\PerDom)\rangle = \langle c_1({\mathcal L}),
H(\PerDom)\rangle.$$

Now, we perform the Chern class evaluation in three parts: the region over the center
point of the spider, where the representative is identified with the
periodic domain $\PerDom$, the region over the three legs $a$, $b$,
and $c$, where the representative is given by a number of cylinders,
and the regions in $U_i\times e_i$ (with $i=\alpha$, $\beta$, or
$\gamma$), where the representative is given by a collection of
disks. To justify this, recall that ${\mathcal L}|\partial \PerDom$ is
canonically identified with the tangent bundle to $F$, where it has a
canonical trivialization, given by the tangents along the
curve. Similarly, over the other three boundary points of $a$, $b$,
and $c$ respectively, ${\mathcal L}$ is identified with the tangent
bundle to the disk, and then calculate the first Chern number
of ${\mathcal L}$ by evaluating the three corresponding relative first
Chern numbers.

For the region over the central point $x$ in the spider, we have the
identification
${\mathcal
L}\cong T\Sigma$, away from tubular neighborhoods of the $g+1$ points
$u(\tau)$ and $z$. Indeed, by a local calculation around these $g+1$
special points (see Lemma~\ref{HolDiskTwo:lemma:RelativeCOne} of~\cite{HolDiskTwo}), it
follows that
\begin{eqnarray*}
\lefteqn{\langle e({\mathcal L},\partial F),F\rangle} \\
&=& \langle e(\Phi^*(T\Sigma),\partial F),F\rangle
+ 2 \#\{x\in F\big|\Phi(x)\in u(\tau)\}
-2 \#\{x\in F \big|\Phi(x)=z\} \\
&=& \EulerMeasure(\PerDom)+2 n_{u(\tau)}(\PerDom) - 2 n_z(\PerDom).
\end{eqnarray*}

Consider, next, the cylinders and then the caps which are added to $F$
to obtain the representative in $X$ -- these cylinders consist of some
number of copies of $\alpha_i'\times a\subset \Sigma\times \Delta$,
$\beta_i'\times b \subset \Sigma\times \Delta$, and 
$\gamma_i'\times c \subset \Sigma\times \Delta$, 
followed then by
some number of caps.  Over the cylinder $\alpha_i\times a$, the
restriction of ${\mathcal L}$ is also canonically identified with the
tangent bundle of the cylinder, except at $\tau\in a$ for which
$u(\tau)$ meets $\alpha_i'$.  In a neighborhood of each such crossing,
the tangent bundle to the cylinder and the actual two-plane field
${\mathcal L}$ differ by a relative Chern number of $2$ (this is the
same as the contribution of the points $x\in F$ mapping to
$u(\tau)$). The contribution of the $\beta_i'$ and $\gamma_i'$ works
the same way. Hence, the contributions from the cylinders is:
$$2 \#(a\cap \alpha_i') + 2 \#(b\cap \beta_i') + 2(c\cap \gamma_i').$$
Finally, each closing disk contributes a relative Chern
number of $1$, since over each disk, ${\mathcal L}$ is canonically
identified with the tangent bundle, but the trivialization on the
boundary is the one induced by the tangent to the boundary. Adding
this up, we
get a contribution
$$\sum_{i=1}^g |a_i|+|b_i|+|c_i|.$$

Adding up all the contributions, we obtain
Equation~\eqref{eq:COneFormula}.
\end{proof}

\subsection{Examples with $g=1$}
\label{subsec:CPtwoBar}

Lemma~\ref{lemma:CPtwoBar} is a straightforward application of 
Proposition~\ref{prop:COneFormula} at this point.

\vskip.3cm
\noindent{\bf{Proof of Lemma~\ref{lemma:CPtwoBar}}}

\begin{figure}
\mbox{\vbox{\epsfbox{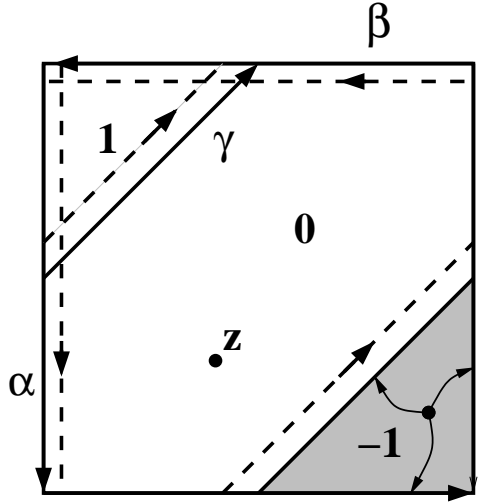}}}
\caption{\label{fig:cptwo}
Three curves in the torus. Note that the square is to be given the
usual edge identifications in this picture.  The dotted lines indicate
the translates of the $\alpha$, $\beta$, and $\gamma$ curves, oriented
as boundaries of the triply-periodic domain whose multiplicities are
shown. A dual spider is included for the small shaded triangle.  Since
its vertex meets the periodic domain with multiplicity $-1$, and its
legs are disjoint from the translated attaching circles, the dual
spider number is $-1$, so the evaluation $\langle c_1(\spinc),
E\rangle=+1$.}
\end{figure}

Consider the genus one Heegaard triple 
$(E,\{\alpha_0\},\{\beta_0\},\{\gamma_0\},z_0)$ for $\mCP$ (punctured in three points)
pictured in Figure~\ref{fig:cptwo}. 

In the picture the multiplicities illustrated represent
the triply-periodic domain associated to an embedded
two-sphere $E\subset \mCP$ with self-intersection number $-1$.

According to Proposition~\ref{prop:COneFormula}, the shaded triangle
-- $\psi^+_1$ in the notation at the beginning of the section --
represents a $\SpinC$ structure 
with $\langle
c_1(\spinc), [E]\rangle = 1$. Specifically, this periodic domain has
\begin{eqnarray*}
\EulerMeasure(\PerDom)&=& \frac{1}{4}-\frac{1}{4}=0 \\
\#(\partial\PerDom)&=& 3 \\
n_z(\PerDom)&=& 0 \\
\sigma(\psi^+_1,\PerDom)&=&-1.
\end{eqnarray*}

For $\psi^+_k$, the only quantity which changes is $\sigma(\psi^+_k,\PerDom)=k-1$,
which gives $\langle c_1(\psi^+_k),[E]\rangle = 2k+1$. 
The formula $\langle c_1(\psi^-_k),[E]\rangle = -2k-1$ follows symmetrically.
\qed

%

\section{Absolute gradings}
\label{sec:AbsGrade}

Let $Y$ be an oriented three-manifold, equipped with a torsion
$\SpinC$ structure (i.e. one for which $c_1(\spinct)$ is torsion).
We have seen that
$\HFc(Y,\spinct)$ is a relatively $\Z$-graded Abelian group: it is
generated by homogeneous elements ${\mathfrak A}$, on which there is a
relative grading function $$\gr\colon {\mathfrak A}\times {\mathfrak
A} \longrightarrow
\Z.$$

\begin{theorem}
\label{thm:AbsGrade}
Let $\spinct$ be a torsion $\SpinC$ structure.
Then, the homology groups $\HFc(Y,\spinct)$ can be endowed with an 
absolute grading 
$$\liftGr\colon {\mathfrak A}\longrightarrow \Q$$
satisfying the following properties:
\begin{itemize}
\item the homogeneous elements of least grading in $\HFp(S^3,\spinc_0)$ have absolute
grading zero
\item the absolute grading lifts the relative grading, in the sense that
if $\xi,\eta\in{\mathfrak A}$, then 
$$\gr(\xi,\eta)=\liftGr(\xi)-\liftGr(\eta)$$
\item the natural maps 
$\iota$ and $\pi$ in the long exact sequence
(Equation~\eqref{eq:HFinfExactSequence}) preserve the absolute
grading, while the coboundary map decreases absolute degree by one,
and the $U$ action decreases it by two.
\item if $W$ is a cobordism from $Y_1$ to $Y_2$ endowed with a $\SpinC$ 
structure whose restriction $\spinct_i$ to $Y_i$ is torsion for $i=1,2$, then
\begin{equation}
\label{eq:DimensionShiftFormula}
\liftGr(F_{W,\spinc}(\xi))-\liftGr(\xi) = 
\frac{c_1(\spinc)^2-2\chi(W)-3\sigma(W)}{4},
\end{equation}
where $\spinct_i=\spinc|Y_i$ for $i=1,2$
\end{itemize}
\end{theorem}

To define $\liftGr$, we present $Y$ as a surgery on a link
$\Link\subset S^3$, so that $\spinct$ is the restriction of a $\SpinC$
structure $\spinc$ over the induced cobordism $W(S^3,\Link)$ from
$S^3$ to $Y$.  Let $(\Sigma,\alphas,\betas,\gammas,z)$ be a Heegaard
triple subordinate to some bouquet for the link 
(in the sense of Defintion~\ref{def:HeegTrip}), so that
$Y_{\alpha,\beta}\cong S^3$, $Y_{\beta,\gamma}\cong \#^{n}(S^1\times S^2)$, and
$Y_{\alpha,\gamma}\cong Y$). Fix intersection points
$\x_0\in\Ta\cap\Tb$, $\x_1\in \Tb\cap\Tc$ so that $\x_0$
resp. $\x_1$ are in the same degree as the highest non-zero 
generators of $\HFa(S^3,\spinct_0)$ 
and $\HFa(\#^n(S^2\times S^1),\spinct_0)$ respectively.
(We say simply that the intersection points $\x_0$ and $\x_1$ lie {\em in the canonical degree}.)

We define the absolute grading on $\CFa(Y,\spinct)$, 
by letting (for each $\y\in\Ta\cap\Tc$)
\begin{equation}
\label{eq:DefAbsGrade}
\liftGr(\y)=-\Mas(\psi)+2n_z(\psi)+\frac{c_1(\spinc)^2-2\chi(W)-3\sigma(W)}{4}
\end{equation}
where $W=W(\Link)$, and $\psi\in\pi_2(\x_0,\x_1,\y)$ is a homotopy
class whose $\SpinC$ structure $\spinc_z(\psi)=\spinc$. This induces
an absolute grading on $\CFinf(Y,\spinct)$ by
$$\liftGr[\y,i]=2i+\liftGr(\y),$$
and hence on the sub- and quotient-complexes
$\CFm(Y,\spinct)$ and $\CFp(Y,\spinct)$ (so that the inclusion and
projection preserve grading).

\subsection{Invariance of the absolute grading}
We show that the grading defined above is well-defined (depending only
on the three-manifold and $\SpinC$ structure), and satisfies the
requirements of Theorem~\ref{thm:AbsGrade}, in a sequence of steps
which are reminiscent of the link invariant constructed in
Subsection~\ref{subsec:TwoHandles}. Commutative diagrams coming from 
associativity are replaced
by (more elementary) index statements.

The definition of the absolute grading as defined in
Equation~\eqref{eq:DefAbsGrade} depends on a link $\Link\subset S^3$,
a $\SpinC$ structure $\spinc\in\SpinC(W(\Link))$ extending $\spinct$,
and a Heegaard triple $(\Sigma,\alphas,\betas,\gammas,z)$ subordinate
to some bouquet for the link.

\begin{prop}
\label{prop:AbsGradeDepLink}
The absolute grading defined in Equation~\eqref{eq:DefAbsGrade} is
independent of the bouquet $B(\Link)$ for the link, and the
subordinate Heegaard triple. Specifically, if
$(\Sigma_1,\alphas_1,\betas_1,\gammas_1,z_1)$ and
$(\Sigma_2,\alphas_2,\betas_2,\gammas_2,z_2)$ are a pair of Heegaard
triples then the absolute grading induced by first Heegaard triple on
a homogeneous $\xi\in\HFc(\alphas_1,\gammas_1,\spinct)$ agrees with the
absolute grading induced by the second Heegaard triple on
$\Psi(\xi)\in\HFc(\alphas_2,\gammas_2,\spinct)$, where $$\Psi\colon
\HFc(\alphas_1,\gammas_1,\spinct)\longrightarrow
\HFc(\alphas_2,\gammas_2,\spinct)$$ is the isomorphism 
induced by the equivalence of Heegaard diagrams.
\end{prop}

This is shown in a sequence of steps.

\begin{lemma}
\label{lemma:Baspoints}
For a fixed Heegaard triple $(\Sigma,\alphas,\betas,\gammas,z)$, 
the absolute grading is independent of the choices of $\x_0$ and $\x_1$ 
in the canonical degree, and independent of the particular triangle
$\psi\in\pi_2(\x_0,\x_1,\y)$ used in its definition.
\end{lemma}

\begin{proof}
First, we show independence of the particular triangle: i.e. fix a
Heegaard triple subordinate to some bouquet for the framed link
$\Link\subset S^3$, and intersection points $\x_0$, $\x_1$ as
before. Suppose that $\psi,\psi'\in\pi_2(\x_0,\x_1,\y)$ are a pair of
triangles with $\spinc_z(\psi)=\spinc_z(\psi')=\spinc$. Then,
$$\psi'=\psi+\phi_{\alpha,\beta}+\phi_{\beta,\gamma}+\phi_{\alpha,\gamma}+\ell
[S],$$ where $\phi_{\xi,\eta}$ are periodic domains for
$Y_{\xi,\eta}$, and $S$ is the generator of
$\pi_2(\Sym^g(\Sigma))$. By additivity of the Maslov index, the fact
that the restriction of $\spinc$ to the three boundary components is
torsion, and the fact that $\langle c_1(T\Sigma), [S]\rangle = 1$, it
follows that $$\Mas(\psi')=\Mas(\psi)+2\ell.$$ Moreover, it is
immediate that $n_z(\psi')=n_z(\psi)+\ell$. Since all the other terms
of Equation~\eqref{eq:DefAbsGrade} remain unchanged, it follows that
$\liftGr$ is independent of $\psi\in\pi_2(\x_0,\x_1,\y)$. Independence
of the choice of $\x_0\in\Ta\cap\Tb$ and $\x_1\in\Tb\cap\Tc$ follows
similarly, with the observation that any two intersection points
$\x_0, \x_0'$ for $\Ta\cap\Tb$ both of which lie in the canonical degree
can be connected by a Whitney disk $\phi\in\pi_2(\x_0,\x_0')$ with 
$\Mas(\phi)=n_z(\phi)=0$.
\end{proof}

\begin{lemma}
\label{lemma:AbsGradeStab}
The absolute grading is invariant under stabilizations of the Heegaard
triple.
\end{lemma}

\begin{proof}
We stabilize the triple $(\Sigma,\alphas,\betas,\gammas,z)$ by forming
the connected sum with the standard genus one diagram
$(E,\alpha_{g+1},\beta_{g+1},\gamma_{g+1},z_0)$ encountered in
Lemma~\ref{lemma:Stabilize}. If $\psi\in\pi_2(\x_0,\x_1,\y)$, then its
stabilization $\psi'=\psi\#\psi_0$ is formed by splicing a standard triangle
$\psi_0\in\pi_2(x_0,y_0,w_0)$ in the genus one
surface. Now, $\psi'\in\pi_2(\x_0',\x_1',\y')$ (where $\y'$ is the
stabilization of $\y$, and $\x_0'$ and $\x_1'$ are corresponding intersection points
in the canonical degrees for the stabilized diagram) satisfies
$\Mas(\psi')=\Mas(\psi)$ and $n_{z'}(\psi')=n_{z}(\psi)$.
Of course, the topological terms in the definition of the absolute
grading remain unchanged.  Thus, $\liftGr(\y)=\liftGr(\y')$.
\end{proof}

\begin{lemma}
\label{lemma:StrongEquivAbsGrade}
If $(\Sigma,\alphas,\betas,\gammas,z)$ and
$(\Sigma',\alphas',\betas',z')$ are strongly equivalent Heegaard
triples, then the absolute gradings are identified.
\end{lemma}

\begin{proof}
First, we observe that we are free to change the $\betas$ by an
isotopy without changing the relative grading. To see this, observe
that if $\Psi_t$ is an isotopy carrying $\betas$ to $\betas'$, then we
can find intersection points $\x_0\in\Ta\cap\Tb$, $\x_1\in\Tb\cap\Tc$,
$\x_0'\in\Ta\cap\TbPr$, $\x_1\in\Tb\cap\TbPr$ satisfying the
requirements of the definition of the absolute grading, and also
Whitney disks with dynamic boundary conditions
$\phi_0\in\pi_2^{\Psi_t}(\x_0',\x_0)$ and
$\phi_1\in\pi_2^{\Psi_t}(\x_1',\x_1)$ with 
$\Mas(\phi_0)= \Mas(\phi_1)=0$ and $n_z(\phi_0)=n_z(\phi_1)=0$.
If $\psi\in\pi_2(\x_0,\x_1,\y)$ is a triangle representing $\spinc$,
then, we claim that the triangle with dynamic boundary conditions
$\psi+\phi_0+\phi_1$ is homotopic to a triangle
$\psi'\in\pi_2(\x_0',\x_1',\y)$ with (stationary boundary conditions) for
the Heegaard triple $(\Sigma,\alphas,\betas',\gammas)$. 
By the homotopy invariance of the Maslov index and the local
multiplicity, 
$\Mas(\psi')=\Mas(\psi)$ and $n_z(\psi')=\Mas(\psi)$.
A similar argument allows us to modify $\alphas$ and 
$\gammas$ by isotopies.

Next, we claim that if the $\gammas'$ are obtained from
the $\gammas$ by a sequence of handleslides and isotopies
(so that $(\Sigma,\betas,\gammas',z)$ is still admissible for
$\#^n(S^2\times S^1)$, and of course, $(\Sigma,\alphas,\gammas',z)$ is
still admissible for $Y$), then the induced gradings on $Y$ are identified. 
This follows from the index analogue of the associativity of triangles. 
Specifically, suppose $\xi\in\HFa(\alphas,\gammas)$ is a homology class, and
$\Phi(\xi)$ is its image under the strong equivalence map. 
Then, there are intersection points $\y\in\Ta\cap\Tc$ and $\y'\in
\Ta\cap\TcPr$ in the same degrees as $\xi$ and 
$\Phi(\xi)$ respectively, and there is also a triangle
$\psi_0\in\pi_2(\y,\x_2,\y')$ with
$n_z(\psi_0)=\Mas(\psi)=0$, where $\x_2\in\Ta\cap\TcPr$ is
an intersection point in the canonical degree. Now, juxtaposing the
triangles $\psi\in\pi_2(\x_0,\x_1,\y)$ and
$\psi_0\in\pi_2(\y,\x_2,\y')$, we obtain a square
$\pi_2(\x_0,\x_1,\x_2,\y')$. We claim that there is also a triangle
$\psi_0'\in\pi_2(\x_1,\x_1',\x_2)$ with
$n_z(\psi_0')=\Mas(\psi_0')=0$, for some $\x_1'\in\Tb\cap\TcPr$; and then 
also a representative $\psi\in\pi_2(\x_0,\x_1',\y')$ with $\spinc_z(\psi)=\spinc$.
Now the juxtapositions of $\psi+\psi_0$ and $\psi'+\psi_0'$ represent
the same $\SpinC$ structure for the Heegaard quadruple so, up to 
adding doubly-periodic domains (which leaves $\Mas$ and $n_z$ unchanged), we 
see that the two squares become homotopic, and hence that
$\Mas(\psi_0)+\Mas(\psi)=\Mas(\psi_0')+\Mas(\psi')$ and also
$n_z(\psi_0)+n_z(\psi)=n_z(\psi_0')+n_z(\psi')$. It follows that 
the $(\Sigma,\alphas,\betas,\gammas,z)$-induced 
grading of $\y$ (which is the grading of $\xi$) coincides
with the $(\Sigma,\alphas,\betas,\gammas',z)$-induced grading of $\y'$ (which 
is the grading of $\Phi(\xi)$). 
We can apply the same arguments to allow handleslides and isotopies amongst
$\alphas$ and $\betas$ as well.
\end{proof}

\vskip.2cm
\noindent{\bf{Proof of Proposition~\ref{prop:AbsGradeDepLink}.}}
Independence of the Heegaard triple subordinate to a given bouquet now
follows immediately from Lemmas~\ref{lemma:AbsGradeStab} and
\ref{lemma:StrongEquivAbsGrade}, in view of
Lemma~\ref{lemma:LinkInvariance}.  In fact, independence of the
bouquet follows as in Lemma~\ref{lemma:IndepOfBouquet}. In that lemma,
we have seen that difference choices of bouquet correspond to
handleslides amongst the $\gammas$, and the proof of
Lemma~\ref{lemma:StrongEquivAbsGrade} actually 
shows that in fact these more general moves leave the
absolute grading unchanged.
\vskip.2cm

Since we have not yet established that $\liftGr$ is
independent of the link and $\SpinC$ structure, we include them in the
notation, writing $\liftGr_{\Link,\spinc}([\y,i])$.

\begin{lemma}
\label{lemma:Composite}
Let $Y=S^3(\Link_1)=S^3(\Link_1')$, and let $\spinc_1$
resp. $\spinc_1'$ be $\SpinC$ structures over $W(\Link_1)$, resp
$W(\Link_1')$ whose restrictions to $Y$ agree (and are torsion). Then
for any other link $\Link_2\subset S^3$ and $\SpinC$ structure
$\spinc_2\in
\SpinC(W(Y,\Link_2))$, we have that
$$\liftGr_{\Link_1\cup\Link_2,\spinc_1\cup \spinc_2}-
\liftGr_{\Link_1'\cup\Link_2,\spinc_1'\cup \spinc_2}
=\liftGr_{\Link_1,\spinc_1}-\liftGr_{\Link_1',\spinc_2'}.$$
\end{lemma}

\begin{proof}
This is the analogue of the composition law for the link invariant
(Proposition~\ref{prop:CompositeLink}).

Fix a Heegaard triple $(\Sigma,\alphas,\betas,\gammas,z)$ subordinate to a link, and let
$x_0\in\Ta\cap\Tb$, $x_1\in\Tb\cap\Tc$, $\x_2\in\Tc\cap\Td$ be elements in the same degree as the
corresponding top-dimensional non-zero homology in $\HFleq(\#^{\ell}(S^1\times S^2))$.
We define
$$\liftGr^\circ(\y)=-\Mas(\psi)+2n_z(\psi).$$ 

Let $\psi_1\in\pi_2(\x_0,\x_1,\y)$ represent $\spinc_1$, and
$\psi_2\in\pi_2(\y,\x_2,\w)$ represent $\spinc_2$. Up to addition of
doubly-periodic domains (which do not change the Maslov index), we can
decompose the square obtained by juxtaposing $\psi_1$ and $\psi_2$ as
another juxtaposition, of $\psi_3\in\pi_2(\x_0,\x_3,\y)$ and
$\psi_4\in\pi_2(\x_1,\x_3,\y)$, for $\x_3\in\Tb\cap\Td$, which we can
also assume lies in the canonical degree. It follows that
$-\Mas(\psi_4)+2n_z(\psi_4)=0$. Now, since the squares are homotopic,
we get that
$$-\Mas(\psi_1)+2n_z(\psi_1)-\Mas(\psi_2)+2n_z(\psi_2)=-\Mas(\psi_3)+2n_z(\psi_3);$$
i.e.
$$\liftGr^{\circ}(\y)-\Mas(\psi_2)+2n_z(\psi_2)=\liftGr^{\circ}(\w),$$
where $\psi_2\in\pi_2(\y,\x_2,\w)$. As in the proof of
Proposition~\ref{prop:AbsGradeDepLink} above, it follows that
$-\Mas(\psi_2)+2n_z(\psi_2)$ depends only on the link $\Link_2$
(thought of as a link in $Y$), the $\SpinC$ structure $\spinc_2$, and,
of course, the gradings of $\y$ and $\w$.  In particular, it is
independent of $\Link_1$ and $\spinc_1$. 

The corresponding excisive property of the other term in $\liftGr$ is obvious.
The lemma follows.
\end{proof}

The three-sphere $S^3$ has a standard absolute grading for which
$\HFa(S^3)$ is supported in degree zero.

\begin{lemma}
\label{lemma:BlowUp}
Let $\Link\subset S^3$ be the unknot with framing $-1$.  Then the
induced absolute grading on $S^3\cong S^3(\Link)$, $\liftGr_{\Link,\spinc}$
induced from any
$\SpinC$ structure over $W(\Link)$ agrees with the standard grading 
of $\HFp(S^3)$.
\end{lemma}

\begin{proof}
A Heegaard triple
subordinate to the unknot with framing $-1$ is
the standard genus one diagram with
three curves $\alpha$, $\beta$ and $\gamma$ each pairwise intersecting
in one point pictured in Figure~\ref{fig:cptwo}.  The $\SpinC$ structures
over $W(\Link)$ are represented by triangles $\psi_k^\pm$ indexed by a
sign $\pm$ and a non-negative integer $k$, as in
Proposition~\ref{HolDiskTwo:prop:HoClassesCancel}
of~\cite{HolDiskTwo}.  These triangles have $\Mas(\psi_k)=0$ and
$n_z(\psi_k)=\frac{k(k-1)}{2}$. According to Lemma~\ref{lemma:CPtwoBar},
$\langle c_1(\spinc_z(\psi^\pm_k)),
[E]\rangle=\pm(2k-1)$.  It follows that $\liftGr[\y,j]=2j$,
independent of the $\SpinC$ structure used over the cobordism (where
here $\y$ is the unique intersection point between $\alpha$ and
$\gamma$).
\end{proof}

\begin{lemma}
\label{lemma:AntiBlowUp}
Let $\Link\subset S^3$ be the unknot with framing $+1$ (observe the sign here).  
Then the
induced absolute grading on $S^3\cong S^3(\Link)$, $\liftGr_{\Link,\spinc}$
induced from any
$\SpinC$ structure over $W(\Link)$ agrees with the standard grading 
of $\HFp(S^3)$.
\end{lemma}

\begin{proof}
Let $W(\Link)$ be the cobordism, and let $H\in H_2(W(\Link);\Z)\cong
\Z$ be a generator. Fix a $\SpinC$ structure $\spinc\in\SpinC(W(\Link))$ with
$\langle c_1(\spinc),H\rangle = 2k+1$. Let $\Link'$ be the new link
obtained from $\Link$ by adding another unknot, this one with framing
$-1$, and let $E\in H_2(W(\Link');\Z)$ be the new homology class
introduced by this unknot. We endow $W(\Link')$ with the $\SpinC$
structure $\spinc'$ with $\langle c_1(\spinc), E \rangle =
-2k-1$. After handlesliding the circle with framing $+1$ over the
circle with framing $-1$, we obtain a new link $\Link''$ consisting of
a pair of circles with linking number one and framings $0$ and
$-1$. Thus, the cobordism decomposes along $S^1\times S^2$. It is easy
to see that $c_1(\spinc')|_{S^1\times S^2}=0$.

We claim that the induced absolute grading on $S^3$ induced from
$\Link''$ is the standard grading. This follows from the fact that for
$\Link''$, we have that the induced grading by
$2i-\Mas(\psi)+2n_z(\psi)$ is shifted up by one (this is the shift
appearing in the the composite cobordism $S^3\Rightarrow (S^1\times
S^2)\Rightarrow S^3$, appearing in the surgery long exact sequence for
the unknot -- observe also that the relative grading here is an
integer, not an integer modulo two, since we are factoring through a
torsion $\SpinC$ structure on $S^1\times S^2$). On the other hand, it
is easy to see that $c_1(\spinc')=0$, so
$$\frac{c_1(\spinc')^2-2\chi(W(\Link''))-3\sigma(W(\Link''))}{4}=-1.$$

Since the absolute grading induced by a link is invariant under
handleslides (this can be seen by adapting Lemma~\ref{lemma:FourDHS} in 
the now familiar manner),
it follows that the induced grading on $S^3$ by $\Link'$ is also the
standard grading. Finally by Lemmas~\ref{lemma:Composite} and
\ref{lemma:BlowUp}, it follows that the absolute grading induced from 
the $\Link$ must be standard, as well. 
\end{proof}

\begin{lemma}
\label{lemma:FourDHSAbsGrade}
The grading on $Y$
induced by a link $\Link\subset S^3$ and a 
$\SpinC$ structure $\spinc\in W(\Link)$ is invariant under handleslides between
components of the link.
\end{lemma}

\begin{proof}
As in Lemma~\ref{lemma:FourDHS}, if $\Link$ and $\Link'$ differ by
four-dimensional handleslides, then we can find Heegaard triples
$(\Sigma,\alphas,\betas,\gammas,z)$ and
$(\Sigma,\alphas',\betas',\gammas',z)$ subordinate to a bouquet for
$\Link$ and $\Link'$, with the property that the $\alphas'$
(resp. $\betas'$ resp. $\gammas'$) are gotten by the $\alphas$
(resp. $\betas$ resp $\gammas$) by a sequence of
handleslides and isotopies. The fact that the gradings remain identified 
now follows from the proof of Lemma~\ref{lemma:StrongEquivAbsGrade}.
\end{proof}

\begin{lemma}
\label{lemma:SThree}
Let $\Link\subset S^3$ be any link with the property that
$S^3(\Link)\cong S^3$. Then for any $\SpinC$ structure $\spinc$ over
$W(\Link)$, the induced absolute grading on $S^3$ is the standard grading.
\end{lemma}

\begin{proof}
According to a theorem of Kirby (see~\cite{KirbyCalculus}), any two
links $\Link$ and $\Link'$ which give rise to $S^3$ can be connected
by a sequence of moves which either: introduce a new, disjoint unknot
with framing $\pm 1$, delete a disjoint unknot with framing $\pm 1$,
(four-dimensional) handleslides amongst the components of $\Link$, and
isotopies of the components of $\Link$. According to
Lemmas~\ref{lemma:BlowUp} and \ref{lemma:AntiBlowUp}, and
\ref{lemma:FourDHSAbsGrade}, it follows that absolute grading is
unchanged under all of these operations.
\end{proof}

\vskip.2cm
\noindent{\bf{Proof of Theorem~\ref{thm:AbsGrade}.}}
We put the above results together to show that the absolute grading we
have defined is independent of the link and $\SpinC$ structure used in
its definition.
Suppose that $\Link_1$ and $\Link_1'$ are a pair of
framed links in $S^3$ so that $S^3(\Link_1)\cong S^3(\Link_1')\cong
Y$, equipped with $\SpinC$ structures $\spinc_1$ and $\spinc_1'$ with
$\spinc_1|Y\cong
\spinct \cong \spinc_1'|Y$. Then, we fix a link $\Link_2\subset Y$
with $Y(\Link_2)\cong S^3$, and $\spinc_2$ so that $\spinc_2|Y\cong
\spinct$. Lemma~\ref{lemma:SThree}, the gradings on $\HFp(S^3)$
induced by the links $\Link_1\cup\Link_2$ (endowed with a $\SpinC$
structure $\spinc$ whose restriction to $W(\Link_1)=\spinc_1$, and
whose restriction to $W(Y,\Link_2)$ agrees with $\spinc_2$) and
$\Link_1'\cup\Link_2$ (endowed with a $\SpinC$ structure whose
restriction to $W(\Link_1)=\spinc_1$, and whose restriction to
$W(Y,\Link_2)$ agres with $\spinc_2$) coincide. The result then follows directly from
Lemma~\ref{lemma:Composite}.

The fact that the absolute grading is a lift of the relative grading
is a direct consequence of the additivity of the Maslov index and
$n_z$ of a triangle under juxtaposition with Whitney disks. In particular, 
the fact that the $H_1$ action decreases degree
by one and $U$ decreases degree by two follows.
Moreover, the degree is by definition 
preserved by $\iota$ and $\pi$, and is compatible with the 
usual grading of $\HFp(S^3,\spinc_0)$. 

It remains to verify Equation~\eqref{eq:DimensionShiftFormula}.  For cobordisms
composed of two-handles entirely, this equation follows from another
application of the index analogue of the associativity. Next, let $U$
be a cobordism from $Y$ to $Y\#(S^2\times S^1)$ consisting of a single
one-handle, and equipped with the $\SpinC$ structure 
${\widehat \spinct}$ whose restriction to $Y$ is $\spinct$. 
We claim that if $\Knot\subset Y$ is is an unknot and $W$
is the cobordism from $Y$ to $Y\#(S^2\times S^1)$ obtained by
zero-surgery on $\Knot$, then we claim that 
$$\pm \gamma\circ \Fc{U,{\widehat\spinct}} =
\Fc{W,\spinc},$$ where $\gamma$ denotes
the action of a generator 
$H_1(Y\#(S^2\times S^1);\Z)$
coming from the $S^1$ factor, and $\spinc$ is the $\SpinC$ structure
over $W$ which with $c_1(\spinc)=0$, and $\spinc|Y=\spinct$. This is
an immediate consequence of the Heegaard diagrams (note that $\Fc{U}$
is the map $\Gc{U,\spinc}$ as defined in
Subsection~\ref{subsec:OneThreeHandles}). Similarly, if $V$ is the
cobordism from $Y\#(S^2\times S^1)$ to $Y$ consisting of a single
three-handle addition (along a sphere $S$ from the $S^2\times S^1$
factor), and $W'$ is the cobordism from $Y\#(S^2\times S^1)$
consisting of a two-handle which cancels the two-sphere, then
$$\pm F_{V,{\widehat\spinct}}\circ \gamma = F_{W',\spinc}$$ for the
obvious choices of $\SpinC$ structure over $V$ and $W'$ (again,
observe that $F_{V,{\widehat\spinct}}$ is the map
$\Ec{V,{\widehat\spinct}}$ from
Subsection~\ref{subsec:OneThreeHandles}). It is an easy consequence of
these observations (and Equation~\eqref{eq:DimensionShiftFormula}, which we
already know holds for the cobordisms $W$ and $W'$) that both $F_V$
and $F_U$ increase degree by $\OneHalf$, verifying
Equation~\eqref{eq:DimensionShiftFormula} for one- and three-handle
additions. Thus, the equation follows in general.

\qed
\vskip.2cm

\subsection{Absolute gradings and duality}

When $\spinct$ is a torsion $\SpinC$ structure, 
there is a duality map (an isomorphism):
$$\Dual\colon \HFp(Y,\spinct)\longrightarrow
\coHFm(Y,\spinct)$$
which, on the chain level, is defined by
$\Dual[\x,i] = \langle [\x,i], \cdot \rangle$
(in the notation of Section~\ref{sec:Duality}); i.e.
$\Dual[\x,i]=-[\x,-i-1]^*$. This sets up an isomorphism of relatively
graded groups.  Indeed, we have the following absolutely graded
version:

\begin{prop}
\label{prop:GradedDuality}
The map 
$$\Dual\colon \HFp_i(Y,\spinct)\longrightarrow \coHFm^{-i-2}(Y,\spinct)$$
is an isomorphism.
\end{prop}

\begin{proof}
The following follows from the fact that the absolute grading on $S^3$ is well-defined.

Let $S^3=Y(\Link)$ for some link $\Link\subset Y$, and
$(\Sigma,\alphas,\gammas,\betas,z)$ be a corresponding  Heegaard triple.
Let
$\psi\in\pi_2(\y,\x_1,\x_0)$, with $\y\in\Ta\cap\Tc$,
$\x_1\in\Tc\cap\Tb$, and $\x_0\in\Ta\cap\Tb$ (where $\x_0$ and $\x_1$
are as before).  Then,
$$\liftGr(\y)=\Mas(\psi)-2n_z(\psi)-\frac{c_1(\spinc)^2-2\chi(W)-3\sigma(W)}{4}.$$

Now, if $(\Sigma,\alphas,\betas,\gammas,z)$ is subordinate to a
bouquet for a link $\Link\subset S^3$ (hence representing a cobordism
from $S^3$ to $Y$), then 
$$(-\Sigma,\alphas,\gammas,\betas,z)$$
represents a cobordism (with only two-handles) from $-Y$ to $-S^3$. We
have an identification $$\pi_2(\x_0,\x_1,\y)\cong
\pi_2(\y,\x_1,\x_0)$$ (where the first group is taken with respect to
the first Heegaard triple, the second to the second Heegaard triple)
by precomposition with the reflection $R$ on the triangle which fixes one
vertex and switches the other two. Observe that
$$\Mas_\Sigma(\psi)=\Mas_{-\Sigma}(\psi\circ R),$$
and
$$n_z^{\Sigma}(\psi)=n_z^{-\Sigma}(\psi\circ R).$$
It follows that
$\liftGr_Y(\y)=-\liftGr_{-Y}(\y)$.

Thus,
$$\liftGr_{-Y} \Dual[\x,i]=\liftGr_{-Y}[\x,-i-1]^*=-2i-2-\liftGr(\x)=-\liftGr_Y([\x,i])-2.$$
\end{proof}

\section{Mixed invariants}
\label{sec:DefMixed}

Let $W$ be a cobordism with $b_2^+(W)>1$, then the cobordism has a 
refined invariant
$$\Fmix{W,\spinc}\colon \HFm(Y_1,\spinc|Y_1)\longrightarrow
\HFp(Y_2,\spinc|Y_2).$$

To define this, we need the following basic results:

\begin{lemma}
Let $K\subset Y$ be a null-homologous knot in an oriented
three-manifold.  Then the map induced on $\HFinf$ by the cobordism $W$
from $Y$ to $Y_n$ (where $n$ is any positive integer; i.e.  the
cobordism has $b_2^+(W)=1$) is trivial.
\end{lemma}

\begin{proof}
In~\cite{HolDiskTwo}, we defined an absolute $\Zmod{2}$ grading on the Floer homology groups,
by observing that the $\Zmod{2}$-graded version of 
$\uHFinf(Y)$ is non-trivial in one degree (which we declare to be even). We claim that 
the map induced by the cobordism $W$ shifts this $\Zmod{2}$ degree. Once we establish this,
it will follow that the map
$$\uFinf{W,\spinc}\colon \uHFinf(Y,\spinc|Y) \longrightarrow \uHFinf(Y_n,\spinc|Y_n)$$
is trivial, and hence, by the universal coefficients spectral sequence, it follows that for any
module $M$ over $\Z[H^1(Y;\Z)]$, the induced map
$$\uFinf{W,\spinc}\colon \uHFinf(Y,\spinc|Y,M) \longrightarrow \uHFinf(Y_n,\spinc|Y_n,M)$$
(and, in particular, the untwisted version) is trivial.

Consider the surgery long exact
sequence relating $Y$, $Y_0$, and $Y_n$: 
$$
\begin{CD}
... @>>> \HFp(Y) @>{F_1}>> \HFp(Y_0) @>{F_2}>> \HF(Y_n) @>>>...
\end{CD}
$$ 

Observe that $F_1$ and $F_2$ are sums of maps induced by
cobordisms. The component of $F_1$ which lands in the torsion $\SpinC$
structure of $Y_0$ induces an isomorphism on $\HFinf$ (with twisted
coefficients), so it follows that it must preserve the absolute
$\Zmod{2}$ degree.  The map $F_2$ shifts absolute $\Zmod{2}$ gradings;
this can be seen by twisting the homology of $Y_0$ by the
$\Z[H^1(Y_0;\Z)]$-module, $\Z[H^1(Y;\Z)]$ (which, in the torsion
$\SpinC$ structure, has homology $\Z$ in each dimension). Thus, the
composite $F_2\circ F_1$ shifts absolute $\Zmod{2}$ degrees.

On the other hand, by the degree shift formula Equation~\eqref{eq:DimensionShiftFormula}
from Theorem~\ref{thm:AbsGrade}, 
the component of $F_2\circ F_1$ which factors
through the torsion $\SpinC$ structure shifts the absolute
grading down by one. This composite corresponds to the map obtained by
a blowup of $W$. Clearly, blowing down the two-sphere does not affect the 
parity of the absolute degree shift (again, by the degree shift formula);
thus, $F_W$ shifts absolute degree down by one as well, modulo two. Thus, it, too, must shift the
absolute $\Zmod{2}$ degree. 
\end{proof}

More generally, we have the following:

\begin{lemma}
\label{lemma:BTwoPlusLemma}
Let $W$ be a cobordism with $b_2^+(W)>0$. Then the induced map
$$\Finf{W,\spinc}\colon \HFinf(Y_1)\longrightarrow \HFinf(Y_2)$$
vanishes. Indeed, for any module $M$ over $\Z[H^1(Y_1;\Z)]$, the map
$\uFinf{W,\spinc}$ vanishes as well.
\end{lemma}

\begin{proof}
We can find an embedded surface $\Sigma\subset W$ with $\Sigma\cdot
\Sigma>0$. Let $M$ be the boundary of the tubular neighborhood of
$\Sigma$. By fixing a path joining $Y_1$ to $\Sigma$, and taking a
regular neighborhood, we break the cobordism apart into a piece from
$Y_1$ to $Y_1\# Q$, and then from $Y_1\# Q$ to $Y_2$, where here $Q$
is the boundary of a tubular neighborhood of $\Sigma$. Since $\delta
H^1(Y_1\# Q)\subset H^2(W)$ is trivial, the map $\Finf{W,\spinc}$
factors through the map induced by the standard cobordism from
$Y_1$ to $Y_1\#Q$. Thus, the result follows 
once we show that the map on $\HFinf$ vanishes for this cobordism. 

We can further break the cobordism from $Y_1$ to $Y_1\# Q$ to 
cobordism from $Y_1$ to 
$Y_1'=Y_1\#^{2g}(S^2\times S^1)$ followed by a single two-handle addition $W'$, to give
$Y_2'=Y_1\# Q$, along a null-homologous knot, with positive framing.
The previous lemma applies to show that the induced map on $\HFinf$ 
is trivial.
(indeed, for any possibly twisting) The result then follows from the
composition law (Theorem~\ref{thm:Composition}).
\end{proof}

\begin{defn}
\label{def:AdmissibleCut}
Fix a $\SpinC$ structure $\spinc$ over $W$.  An {\em admissible cut}
of a cobordism $W$ is a three-manifold $N\subset W$ which divides $W$
into two pieces $W_1$ and $W_2$ with $b_2^+(W_1), b_2^+(W_2)>0$, and
$\delta H^1(N;\Z)=0$ in $H^2(W,\partial W;\Z)$.
\end{defn}

\begin{example}
Suppose that $b_2^+(W)>1$.
Let $\Sigma\subset W$ be an embedded surface with $\Sigma\cm
\Sigma>0$, and let $Q$ denote the boundary of its tubular
neighborhood. Then, both $Y_1\# Q$ and $Q\# Y_2$ determine 
admissible cuts for $W$. 
\end{example}

Thanks to  Lemma~\ref{lemma:BTwoPlusLemma},
the image of the map
$$\Fm{W_1,\spinc|W_1}\colon \HFm(Y_1,\spinc|_{Y_1})\longrightarrow
\HFm(N,\spinc|N)$$
lies in $\HFmred(N,\spinc|N)$. Moreover (by the same lemma), the
map 
$$\Fp{W_2,\spinc|{W_2}}\colon \HFp(N,\spinc|{N})
\longrightarrow \HFp(Y_2,\spinc|{Y_2}).$$ factors through
the projection of $\HFp(N,\spinc|N)$ to $\HFpred(N,\spinc|N)$.
Thus, we can 
define
$$\Fmix{W,N,\spinc}\colon \HFm(Y_1,\spinc|Y_1)\otimes_\Z
\Wedge^* \left(H_1(W;\Z)/\Tors(W)\right)
\longrightarrow 
\HFp(Y_2,\spinc|Y_2)$$
to be the composite:
$$\Fp{W_2,\spinc|W_2}\circ \tau^{-1}\circ \Fm{W_1,\spinc|W_1},$$
where
$$\tau\colon \HFpred(N,\spinc|N)\longrightarrow
\HFmred(N,\spinc|N)$$
is the natural isomorphism induced by the coboundary map 
for the exact sequence of Equation~\eqref{eq:HFinfExactSequence}.

\begin{theorem}
\label{thm:MixedInvariant}
Let  $W$ be a cobordism
with $b_2^+(W)>1$, and 
let $N$ and $N'$ be a pair of admissible cuts. Then the mixed invariants
defined by these cuts coincide.
\end{theorem}

In the proof of this result, we will make repeated use of the following:

\begin{lemma}
\label{lemma:SubCuts}
If we have a pair of admissible cuts $N, N'\subset W$ which are disjoint,
then the invariants $\Fmix{W,N,\spinc}$ and $\Fmix{W,N',\spinc}$ agree.
\end{lemma}

\begin{proof}
This follows from the following commutative diagram
\begin{equation}
\label{eq:SubCuts}
\begin{CD}
 	& & \HFpRed(N,\spinc|N) @>>> \HFpRed(N',\spinc|N')
 	@>>> \HFp(Y_2,\spinc|Y_2) \\
	& & @V{\tau_1}VV		@V{\tau_2}VV	\\
\HFm(Y_1,\spinc|Y_1) @>>> \HFmRed(N,\spinc|N) @>>> \HFmRed(N',\spinc|N'),
\end{CD}
\end{equation}
which in turn follows from the naturality of the long exact sequences 
of Theorem~\ref{thm:CobordismInvariant} (see Remark~\ref{rmk:Naturality}).
\end{proof}

\vskip.2cm
\noindent{\bf{Proof of Theorem~\ref{thm:MixedInvariant}.}}
We will repeatedly make use of the following trick: if $c_1, c_2\in
H_2(X;\Z)$ are a pair of homology classes with $c_1\cm c_2=0$, then
there are smoothly embedded representatives $\Sigma_1$ and $\Sigma_2$
which are disjoint. To find these, take any pair of smoothly embedded
representatives $\Sigma_1$ and $\Sigma_0$ for $c_1$ and $c_2$
respectively, chosen so that they meet transversally. Since $c_1\cm
c_2=0$, the intersection points between $\Sigma_1$ and $\Sigma_0$ come
in canceling pairs. $\Sigma_2$ is obtained from $\Sigma_0$ by
surgering out the pairs of canceling intersection points, and
attaching cylinders which are disjoint from $\Sigma_1$ (we can find
such cylinders in a small tubular neighborhood of $\Sigma_1$).

Now, when $b_2^+(W)\geq 3$, we proceed as follows. Consider cuts $N$ and
$N'$ which divide $W$ into $W_1\cup_{N} W_2$ and $W_1'\cup_{N'} W_2'$
respectively.  Let $\Sigma_1\subset W_1$ and $\Sigma_2\in W_1'$ be a
pair of disjoint surfaces with positive square. Since $b_2^+(W)\geq 3$, it
follows that $b_2^+(W-\Sigma_1-\Sigma_2)\geq 1$, so we can find
another surface $\Sigma_3\in W-\Sigma_1\cup\Sigma_2$ with positive
square. Let $Q_1$, $Q_2$, and $Q_3$ be the boundaries of the
tubular neighborhoods. Then, by Lemma~\ref{lemma:SubCuts}
$Y_1\# Q_1$ and
$Y_1\# Q_2$ give cuts calculating $\Fmix{W,N,\spinc}$ and
$\Fmix{W,N',\spinc}$ respectively. Applying the lemma once more, we
see that both of these invariants agree with the invariant calculated
using the cut $Q_3\# Y_2$.

When $b_2^+(W)=2$, we make use of the blowup formula, as follows.  Fix
cuts $Q_1$ and $Q_2$ and corresponding surfaces $\Sigma_1\subset W_1$
and $\Sigma_2\subset W_1'$ as above. Blowing up near $Y_2$, the
blowup formula gives
\begin{eqnarray*}
\Fmix{W\#\mCP,N,{\widehat \spinc}}=\Fmix{W,N,\spinc} &{\text{and}}&
\Fmix{W\#\mCP,N',{\widehat \spinc}}=\Fmix{W,N',\spinc}
\end{eqnarray*}
where ${\widehat \spinc}$ is the $\SpinC$ structure over $W\#\mCP$
which agrees with $\spinc$ in a complement of $\mCP$, and whose first
Chern class $c_1(\widehat \spinc)$ evaluates as $+1$ on an exceptional
sphere $E$ inside $\mCP$.  We can assume that $\Sigma_1\cm
\Sigma_2\geq 0$.  In the family $n([\Sigma_1]+[\Sigma_2])+k E$, we can
find a homology class $w$ with
\begin{eqnarray*}
w^2 &>& 0 \\
{[\Sigma_1]}^2 + w^2 - ([\Sigma_1]\cm w)^2 &<& 0 \\
{[\Sigma_2]}^2 + w^2 - ([\Sigma_2]\cm w)^2 &<& 0.
\end{eqnarray*}
Let $\Sigma_3$ be an embedded surface representing $w$. The above
conditions ensure that $b_2^+(W\#\mCP-\Sigma_1-\Sigma_3)=1$, and
$b_2^+(W\#\mCP-\Sigma_2-\Sigma_3)=1$. Applying
Lemma~\ref{lemma:SubCuts}, it follows that the invariants of $W\#
\mCP$ calculated using the cuts corresponding to $\Sigma_1$ and
$\Sigma_2$ agree with that calculated using $\Sigma_3$. 
Thus, the result follows.
\qed

In view of Theorem~\ref{thm:MixedInvariant}, we will drop the cut for
the notation of the mixed invariant.

We also have the following result:

\begin{prop}
Let $W$ be a cobordism from $Y_1$ to $Y_2$ with $b_2^+(W)>1$, and fix
$\SpinC$ structures $\spinct_1$ and $\spinct_2$ over $Y_1$ and $Y_2$
respectively. Then, there are only finitely many $\SpinC$ structures over $W$ for 
which $\Fmix{W,\spinc}$ is non-trivial. 
\end{prop}

\begin{proof}
This follows from Theorem~\ref{thm:Finiteness}. Fix an admissible cut
$W=W_1\#_N W_2$. Since the kernel of the map from
$\HFm(Y)$ to $\HFinf(Y)$ (where we add up all $\SpinC$ structures) is
$\HFmRed(Y)\cong \HFpRed(Y)$, which  is a finitely generated $\Z$-module, we
can find an integer $d$ large enough that $U^d \HFm(N,\spinct)$
injects into $\HFinf(N,\spinct)$ for every $\SpinC$ structure.  
Let $\xi_1,...,\xi_n$ be the generators for
$\HFm(Y_1,\spinct_1)$ as a $\Z[U]$ module.  According to the
finiteness theorem, there is a finite subset ${\mathfrak S}_2\subset
\SpinC(W_2)$ with the property that for $\spinc_2\in{\mathfrak S}_2$,
the map $F^+_{W_2,\spinc_2}$ is non-trivial.  Moreover, according to
the same theorem, there is a finite subset ${\mathfrak S}_1\subset
\SpinC(W_1)$ consisting of elements $\spinc_1\in{\mathfrak S}_1$ for
which there is some element $\xi_i$ with
$F^-_{W_1,\spinc_1}(\xi_i)\not\in U^d \HFm(N,\spinc|N)$. Of course, if
$\spinc$ is any $\SpinC$ structure over $W$ whose restriction to $W_1$
is not in ${\mathfrak S}_1$, the composite
$$\begin{CD}
\HFm(Y_1,\spinc|Y_1)@>{F^-_{W_1,\spinc|W_1}}>>
\HFm(N,\spinc|N) @>>> \HFmRed(N,\spinc|N)
\end{CD}
$$ is trivial, so that $\Fmix{W,\spinc}$ is trivial.  Thus, if
$\Fmix{W,\spinc}$ is non-trivial, its restrictions to $W_1$ and $W_2$
are constrained to lie in the finite sets ${\mathfrak S}_1$ and
${\mathfrak S}_2$ respectively. Since $\delta H^1(N;\Z)=0$, the
$\SpinC$ structure is uniquely determined by its restrictions to $W_1$
and $W_2$.
\end{proof}

\subsection{Other cuts}
When we drop the hypothesis that $\delta H^1(N;\Z)\subset H^2(W;\Z)$
is trivial, then we can still get information about sums of mixed
invariants.

Now we have the following formula for a sum of invariants:

\begin{lemma}
\label{lemma:OtherCut}
Suppose that $W$ is separated by $N$ into a pair of cobordisms
$W_1\cup_N W_2$ with $b_2^+(W_i)>0$. Then, we can still form the
composite $$\Fp{W_2,\spinc|W_2}\circ \tau^{-1}\circ
\Fm{W_1,\spinc|W_1},$$ where $$\tau\colon
\HFpred(N,\spinc|N)\longrightarrow \HFmred(N,\spinc|N)$$ is the map as
before. The composite can be expressed as a sum:
$$\Fp{W_2,\spinc|W_2}\circ \tau^{-1}\circ
\Fm{W_1,\spinc|W_1}=
\sum_{\{\spinc\in\SpinC(W)|\spinc|W_1=\spinc_1,\spinc|W_2=\spinc_2\}}
\Fmix{W,\spinc}.$$
\end{lemma}

\begin{proof}
Since $b_2^+(W_i)>0$, the maps on $\HFinf$ induced by these cobordisms vanish,
according to Lemma~\ref{lemma:BTwoPlusLemma}, we can find
$\Fp{W_2,\spinc|W_2}\circ \tau^{-1}\circ
\Fm{W_1,\spinc|W_1}$.
We find a further subdivision of $W_1=W_0 \cup_{N'} W_1'$ so that $N'$
is an admissible cut for $W$. This can be done by letting $W_0$ be a
neighborhood of an embedded surface in $W_1$ with positive square, and
$W_1'$ be its complement. The result then follows from the arguments from 
Commutative Diagram~\eqref{eq:SubCuts}, together with the composition law, 
Theorem~\ref{thm:Composition}.
\end{proof}

\section{Closed four-manifold invariants}
\label{sec:ClosedInvariant}

Let $X$ be a closed four-manifold with $b_2^+(X)\geq 2$. Then, by deleting a pair of 
disjoint balls from $X$, we can view it as a cobordism $W$ from $S^3$ to $S^3$.

We can now define the {\em absolute invariant of $X$} to be the map
$$\Phi_{X,\spinc}\colon \Z[U]\otimes \Wedge^*(H_1(X)/\Tors) 
\longrightarrow \Z/\pm 1$$
given by $\Phi_{X,\spinc}(U^n\otimes \zeta)$ is the coefficient of
$\Theta^+\in\HFp(S^3,\spinc)$ in the expression
$\Fmix{W,\spinc}(U^n\cm \Theta_-\otimes \zeta)$, where
$\Theta_-\in\HFm(S^3,\spinc_0)$ resp.  $\Theta_+\in\HFp(S^3,\spinc_0)$
is a generator whose degree is maximal resp. minimal.  Of course,
$\Phi_{X,\spinc}$ is vanishes on those homogeneous elements whose
degree is different from
$$d(\spinc)=\frac{c_1(\spinc)^2-2\chi(X)-3\sigma(X)}{4}.$$

\begin{theorem}
The map $$\Phi_{X,\spinc}\colon \Z[U]\otimes \Wedge^*(H_1(X)/\Tors)
\longrightarrow \Z/\pm 1$$ is a smooth, oriented four-manifold
invariant. In particular, if $f\colon X_1\longrightarrow X_2$ is an
orientation-preserving diffeomorphism, then
$$\Phi_{X_1,f^*(\spinc)}(U^n\otimes \zeta)=\Phi_{X_2,\spinc}(U^n\otimes f_*(\zeta)).$$
\end{theorem}

\begin{proof}
This is an immediate consequence of Theorem~\ref{thm:MixedInvariant}.
\end{proof}

\section{Properties of the closed four-manifold invariant}
\label{sec:BasicProperties}

We now return to the closed four-manifold invariant introduced in
Section~\ref{sec:ClosedInvariant}. We first state principle which is
obvious from the definition of the invariant, which implies
Theorem~\ref{intro:VanishingTheorem}. We then turn to the adjunction
inequality stated in Theorem~\ref{intro:Adjunction}

\subsection{Vanishing theorems}

\begin{theorem}
\label{thm:Vanish}
Let $Y$ be a rational homology three-sphere and $\spinc\in\SpinC(Y)$
with $\HFred(Y,\spinct)=0$. Suppose that 
$X$ is a smooth, closed, oriented four-manifold
which admits a decomposition $X=X_1\#_Y X_2$, with $b_2^+(X_1), b_2^+(X_2)>0$. 
Then, for each $\spinc\in\SpinC(X)$ with $\spinc|Y=\spinct$, we have that
$\Phi_{X,\spinc}=0$.
\end{theorem}

\begin{proof}
The hypotheses guarantee that the cut along $Y$ is admissible in
the sense of Definition~\ref{def:AdmissibleCut}. In particular, the map $\Fmix{}$
which is used to define $\Phi_{X,\spinc}$ (i.e. the mixed invariant
associated to the cobordism obtained by puncturing $X$) factors through
$\HFred(Y,\spinct)=0$, so it must vanish.
\end{proof}

Since $\HFred(S^3)=0$, the above theorem gives Theorem~\ref{intro:VanishingTheorem}.

Observe that the blow-up formula for the cobordism invariants can be
rephrased for the absolute invariant, as well. The following is a restatement of
Theorem~\ref{intro:BlowUp}:

\begin{theorem}
\label{thm:BlowUpClosed}
Let $X$ be a closed, smooth, four-manifold with $b_2^+(X)>1$,
and let ${\widehat X}=X\#\mCP$ be its blowup.
Then, for each $\SpinC$ structure ${\widehat \spinc}\in\SpinC({\widehat X})$, with
$d({\widehat X},{\widehat \spinc})\geq 0$ we have the relation
$$\Phi_{{\widehat X},{\widehat \spinc}}(U^{\frac{\ell(\ell+1)}{2}}\cm \xi)=\Phi_{X,\spinc}(\xi),$$
where $\spinc$ is the $\SpinC$ structure over $X$ which agrees over $X-B^4$ with the restriction of
${\widehat\spinc}$, $\xi\in\Z[U]\otimes\Wedge^*(H_1(X)/\Tors)$ is any element of degree 
$d(X,\spinc)$, and $\ell$ is calculated by $\langle c_1({\widehat \spinc}),[E]\rangle = \pm (2\ell+1)$,
where $E\subset{\widehat X}$ is the exceptional sphere.
\end{theorem}

\begin{proof}
Suppose that $N$ is an allowable cut for $X=X_1\#_N X_2$.  Then, we
can decompose 
${\widehat X}=X_1\#_N X_2 \#_{S^3} \mCP$, and
still use $N$ as the cut.  The theorem then follows from the
composition law for the cobordism $X_2\#(\mCP-B^4)$, together with the
blowup formula for the maps of cobordisms (Theorem~\ref{thm:BlowUp}).
\end{proof}

\subsection{Adjunction inequalities}

We turn to a proof of Theorem~\ref{intro:Adjunction}. 

\vskip.3cm
\noindent{\bf{Proof of Theorem~\ref{intro:Adjunction}}.}
As is standard practice, we can reduce to the case where $\Sigma\cm\Sigma=0$ with the
help of the blowup-formula. 

Since we have assumed that $b_2^+(X)>1$, we can find an admissible cut
$N$ which is disjoint from $\Sigma$. Specifically, we can find another
smoothly embedded surface $T\subset X$ with positive self-intersection
number, which is disjoint from $\Sigma$, and use its tubular
neighborhood as the admissible cut. Thus, letting $W$ be complement in
$X$ of a pair of balls, we $W=W_1\cup W_2$, where $W_2$ contains the
embedded surface $\Sigma$. Now, the map $$\Fp{W_2,\spinc|W_2}\colon
\HFp(N,\spinc|N)\longrightarrow \HFp(S^3)$$ factors through
$\HFp(S^1\times \Sigma,\spinc|S^1\times \Sigma)$, where we take
$S^1\times \Sigma\subset X_2$ to be the tubular neighborhood of
$\Sigma$. According to the adjunction inequality for the
three-manifold $S^1\times \Sigma$, this group is zero -- forcing
$\Phi_{X,\spinc}=0$ -- unless 
$$\big|\langle c_1(\spinc),[\Sigma]\rangle\big|
\leq 2g-2.$$
\qed

\commentable{
\bibliographystyle{plain}
\bibliography{biblio}
}

\end{document}